\documentclass{ucthesis}

\ssp

\newcommand{\pathtotrunk}{./}

\usepackage{amsmath,amssymb,amsfonts,amsthm}
\usepackage{ifpdf}
\usepackage{ifthen}

 \usepackage[plainpages=false,hypertexnames=false,pdfpagelabels]{hyperref}
 \newcommand{\arxiv}[1]{\href{http://arxiv.org/abs/#1}{\tt arXiv:\nolinkurl{#1}}}
 \newcommand{\doi}[1]{\href{http://dx.doi.org/#1}{{\tt DOI:#1}}}
 
\ifpdf
\usepackage[pdftex,all,color]{xy}
\else
\usepackage[all,color]{xy}
\fi
\SelectTips{cm}{}%

\theoremstyle{plain}
\newtheorem*{conj}{Conjecture}
\newtheorem{prop}{Proposition}[section]
\newtheorem{thm}[prop]{Theorem}
\newtheorem{lem}[prop]{Lemma}
\newtheorem*{lem*}{Lemma}
\newtheorem{cor}[prop]{Corollary}
\newtheorem*{cor*}{Corollary}
\theoremstyle{definition}
\newtheorem{defn}[prop]{Definition}         
\newtheorem*{defn*}{Definition}             

\newenvironment{rem}{\noindent\textsl{Remark.}}{}  
\numberwithin{equation}{section}

\newcounter{comment}
\newcommand{\noop}[1]{}

\ifpdf
\newcommand{\pathtodiagrams}{\pathtotrunk diagrams/pdf/}
\else
\newcommand{\pathtodiagrams}{\pathtotrunk diagrams/eps/}
\fi

\newcommand{\mathfig}[2]{{\hspace{-3pt}\begin{array}{c}%
  \raisebox{-2.5pt}{\includegraphics[width=#1\textwidth]{\pathtodiagrams #2}}%
\end{array}\hspace{-3pt}}}

\newcommand{\Natural}{\mathbb N}
\newcommand{\Integer}{\mathbb Z}
\newcommand{\Rational}{\mathbb Q}
\newcommand{\Complex}{\mathbb C}
\newcommand{\Real}{\mathbb R}

\newcommand{\abs}[1]{\left\vert#1\right\vert}
\newcommand{\vect}[1]{\overrightarrow{#1}}
\newcommand{\numberof}[1]{\sharp\left(#1\right)}

\newcommand{\bdy}{\partial}
\newcommand{\eset}{\emptyset}

\newcommand{\Id}{\boldsymbol{1}}
\renewcommand{\imath}{\mathfrak{i}}
\renewcommand{\jmath}{\mathfrak{j}}

\newcommand{\ssum}[1]{\Sigma#1}

\newcommand{\csl}[1]{\mathfrak{sl}_{#1}}
\newcommand{\csa}{\mathfrak{h}}
\newcommand{\usl}[1]{U\!\left(\csl{#1}\right)}
\newcommand{\uqsl}[1]{U_q\!\left(\csl{#1}\right)}

\newcommand{\V}[2]{V_{#1<#2}}
\newcommand{\im}{\imath_{-1}}
\newcommand{\iz}{\imath_0}

\newcommand{\A}{\mathcal{A}}
\newcommand{\qRing}{\Integer[q,q^{-1}]}
\newcommand{\qi}[2][q]{\left[#2\right]_{#1}}
\newcommand{\qBinomial}[3][q]{\genfrac{[}{]}{0pt}{}{#2}{#3}_{#1}}
\newcommand{\qPoch}[3]{\left(#1;#2\right)_{#3}}

\newcommand{\sumhat}{\overline{\sum}}
\newcommand{\sumtah}{\underline{\sum}}

\newcommand{\comult}{\Delta}
\newcommand{\counit}{\varepsilon}
\newcommand{\antipode}{S}

\newcommand{\pairing}{p}
\newcommand{\Spairing}{\pairing^{\mathcal{S}}}%
\newcommand{\copairing}{c}
\newcommand{\vout}[2]{v^{#1}_{#2}}
\newcommand{\vin}[2]{w^{#1}_{#2}}

\newcommand{\FreeCat}{\mathcal{T}}
\newcommand{\PivCat}{{\mathcal{P}iv}}
\newcommand{\SymCat}{{\mathcal{S}ym}}
\newcommand{\Draw}{\mathbf{\operatorname{Draw}}}

\newcommand{\GT}[1][]{\mathcal{GT}\parenthesise{#1}}
\newcommand{\dGT}[1][]{d\GT\parenthesise{#1}}
\newcommand{\dGTe}[1][]{\dGT_{\eset}\parenthesise{#1}}
\newcommand{\dGTa}[1][]{\dGT_{\bdy}\parenthesise{#1}}

\renewcommand{\SS}{\mathcal{SS}}
\newcommand{\APR}{\mathcal{APR}}
\newcommand{\AP}{\mathcal{AP}}
\renewcommand{\P}[3]{{\mathcal{P}^{#1}_{#2;#3}}}
\newcommand{\AQR}{\mathcal{AQR}}
\newcommand{\AQ}{\mathcal{AQ}}
\newcommand{\Q}[3]{{\mathcal{Q}^{#1}_{#2;#3}}}
\newcommand{\Label}{\mathcal{L}}

\newcommand{\rotl}[1][]{\operatorname{rotl}\parenthesise{#1}}
\newcommand{\rotr}[1][]{\operatorname{rotr}\parenthesise{#1}}
\newcommand{\drotl}[1][]{\operatorname{drotl}\parenthesise{#1}}

\newcommand{\Braids}{\mathfrak{B}}

\newcommand{\Kekule}{Kekul\'{e} }

\newcommand{\Rep}{\protect\operatorname{\mathbf{Rep}}}
\newcommand{\SRep}{\mathcal{S}\Rep}%
\newcommand{\FundRep}{\mathbf{FundRep}}
\newcommand{\SFundRep}{\mathcal{S}\!\FundRep}
\newcommand{\Mat}[1][]{\operatorname{Mat}\parenthesise{#1}}
\newcommand{\Kar}[1][]{\operatorname{Kar}\parenthesise{#1}}

\newcommand{\To}{\rightarrow}
\newcommand{\Into}{\hookrightarrow}
\newcommand{\Onto}{\twoheadrightarrow}
\newcommand{\Iso}{\cong}

\newcommand{\IsoTo}{\overset{\tiny\Iso}{\To}}

\newcommand{\restrict}[2]{#1{}_{\mid #2}{}}
\newcommand{\set}[1]{\left\{#1\right\}}
\newcommand{\setc}[2]{\left\{#1 \;\left| \; #2 \right. \right\}}

\newcommand{\Span}[2]{\operatorname{span}_{#1}\left\{\displaystyle{#2}\right\}}
\newcommand{\splitSpan}[3]{\notag \operatorname{span}_{#1}\left\{\displaystyle{#2}\vphantom{\displaystyle{#3}}\right. \\ & \qquad  \qquad \qquad \left.\displaystyle{#3}\vphantom{\displaystyle{#2}}\right\}}
\newcommand{\splitSpanThree}[4]{\notag \operatorname{span}_{#1}\left\{\displaystyle{#2}\vphantom{\displaystyle{#3#4}}\right. \\ \notag & \qquad  \qquad \qquad \left.\displaystyle{#3}\vphantom{\displaystyle{#2#4}}\right. \\ & \qquad  \qquad \qquad \left.\displaystyle{#4}\vphantom{\displaystyle{#2#3}}\right\} }

\newcommand{\compose}{\circ}
\newcommand{\inner}[2]{#1 \cdot #2}%

\newcommand{\psmallmatrix}[1]{\left(\begin{smallmatrix} #1 \end{smallmatrix}\right)}

\newcommand{\alt}{\wedge}
\newcommand{\directSum}{\oplus}
\newcommand{\DirectSum}{\bigoplus}
\newcommand{\tensor}{\otimes}
\newcommand{\Tensor}{\bigotimes}

\newcommand{\Hom}[3]{\operatorname{Hom}_{#1}\!\left(#2,#3\right)}

\newcommand{\parenthesise}[1]{\ifthenelse{\equal{#1}{}}{}{\left(#1\right)}}

\def\clap#1{\hbox to 0pt{\hss#1\hss}}

\ifpdf
\usepackage[usenames,dvipsnames]{color}
\usepackage[pdftex]{graphicx}
\else
\usepackage{breakurl}
\usepackage[usenames,dvips]{color}
\usepackage[dvips]{graphicx}
\fi

\begin{document}


\title{A Diagrammatic Category for the Representation Theory of $\uqsl{n}$}
\author{Kim Morrison}
\degreesemester{Spring}%
\degreeyear{2007}%
\degree{Doctor of Philosophy}%
\chair{Professor Vaughan Jones}%
\othermembers{Professor Richard Borcherds \\
  Professor Steve Evans}%
\numberofmembers{3}%
\prevdegrees{B.Sc. (Hons) (University of New South Wales) 2001}%
\field{Mathematics}%
\campus{Berkeley}%

\maketitle%
\copyrightpage%

\begin{abstract}
This thesis provides a partial answer to a question posed by Greg
Kuperberg in \cite{q-alg/9712003} and again by Justin Roberts as
problem 12.18 in \emph{Problems on invariants of knots and
3-manifolds} \cite{math.GT/0406190}, essentially:
\begin{quote}
Can one describe the category of representations of the quantum
group $\uqsl{n}$ (thought of as a spherical category) via generators
and relations?
\end{quote}

For each $n \geq 0$, I define a certain tensor category of trivalent
graphs, modulo isotopy, and construct a functor from this category
onto (a full subcategory of) the category of representations of the
quantum group $\uqsl{n}$. One would like to describe completely the
kernel of this functor, by providing generators for the tensor
category ideal. The resulting quotient of the diagrammatic category
would then be a category equivalent to the representation category
of $\uqsl{n}$.

I make significant progress towards this, describing certain
elements of the kernel, and some obstructions to further elements.
It remains a conjecture that these elements really generate the
kernel. The argument is essentially the following. Take some
trivalent graph in the diagrammatic category for some value of $n$,
and consider the morphism of $\uqsl{n}$ representations it is sent
too. Forgetting the full action of $\uqsl{n}$, keeping only a
$\uqsl{n-1}$ action, the source and target representations branch
into direct sums, and the morphism becomes a matrix of maps of
$\uqsl{n-1}$ representations. Arguing inductively now, we attempt to
write each such matrix entry as a linear combination of diagrams for
$n-1$. This gives a functor $\dGT$ between diagrammatic categories,
realising the forgetful functor at the representation theory level.
Now, if a certain linear combination of diagrams for $n$ is to be in
the kernel of the representation functor, each matrix entry of
$\dGT$ applied to that linear combination must already be in the
kernel of the representation functor one level down. This allows us
to perform inductive calculations, both establishing families of
elements of the kernel, and finding obstructions to other linear
combinations being in the kernel.

This thesis is available electronically from the arXiv, at \arxiv{0704.1503},
and at \url{http://tqft.net/thesis}.

\end{abstract}

\begin{frontmatter}
\setcounter{tocdepth}{1}%
\tableofcontents
\begin{acknowledgements}
First I'd like to thank Vaughan for being such a great advisor. I'm
glad he took a chance on me, and I hope he's not too disappointed by
the complete absence of subfactors in what follows! I'm grateful for
so many things; good advice, lots of mathematics, constant interest
and encouragement, gentle reminders to keep climbing and work in
balance, windsurfing lessons, and a great friendship.

Thanks also to Dror Bar-Natan and Kevin Walker. I've learnt a ton
from each of them, thoroughly enjoyed working with them, and look
forward to more! And further thanks to Greg Kuperberg, for
introducing me to the subject this thesis treats, and taking
interest in my work. Conversations with Joel Kamnitzer, Mikhail
Khovanov, Ari Nieh, Ben Webster, Noah Snyder and Justin Roberts
helped me along the way.

Thanks to my family, Pam, Graham and Adele Morrison, for endless
love and support. And finally, thanks to my friends Nina White,
Yossi Farjoun, Rahel Wachs, Erica Mikesh, and Carl Mautner.
\end{acknowledgements}

\end{frontmatter}

\chapter{Introduction}
\section{Summary}
The eventual goal is to provide a diagrammatic presentation of the
representation theory of $\uqsl{n}$. The present work describes a
category of diagrams, along with certain relations amongst these
diagrams, and a functor from these diagrams to the representation
theory. Further, I conjecture that the relations given are in fact
all of them---that this functor is an equivalence of categories.

We begin by defining a `freely generated category of diagrams', and show that there's a well-defined functor from this
category to the category of representations of $\uqsl{n}$. Essentially, this is a matter of realising that the representation
category is a \emph{pivotal category}, and, as a pivotal category, it is finitely generated. It's then a matter of trying to find the kernel
of this functor; if we could do this, the quotient by the kernel would give the desired diagrammatic category equivalent to the representation category.

This work extends Kuperberg's work \cite{q-alg/9712003} on $\Rep \uqsl{3}$, and agrees with a previously conjectured
\cite{math.QA/0310143} description of $\Rep \uqsl{4}$. Some, but not all,
of the relations have been presented previously in the context of
the quantum link invariants \cite{math.GT/0506403,MR2000a:57023}.
There's a detailed discussion of connections with
previous work in \S \ref{sec:relationships}.

For each $n \geq 0$ we have a category of (linear combinations of)
diagrams
\begin{equation*}
\PivCat_n = \left(\left.
    \mathfig{0.1}{webs/inwards_vertex_abc}, \qquad \mathfig{0.1}{webs/outwards_vertex_abc}, \qquad \mathfig{0.06}{duals/dtr_k},  \qquad \mathfig{0.06}{duals/dar_k}
    \right|
    \begin{array}{c}
    a+b+c=n \\
    1 \leq k \leq n-1
    \end{array}
\right)_{\text{pivotal}}.
\end{equation*}
with edges labelled by integers $1$ through $n-1$, generated by two
types of trivalent vertices, and orientation reversing `tags', as shown above. We allow arbitrary planar
isotopies of the diagrams. This category has no relations; it is a
free pivotal category.

We can construct a functor from this category into the
representation theory
\begin{equation*}
\Rep_n : \PivCat_n \To \Rep \uqsl{n}.
\end{equation*}
This functor is well-defined, in that isotopic diagrams give the
same maps between representations. Our primary goal is thus to
understand this functor, and to answer two questions:

\begin{enumerate}
\item Is $\Rep_n$ full? That is, do we obtain all morphisms between
representations?
\item What is the kernel of $\Rep_n$? When do different diagrams give
the same maps of representations? Can we describe a diagrammatic
quotient category which is equivalent to the representation theory?
\end{enumerate}

The first question has a relatively straightforward answer. We do
not get all of $\Rep \uqsl{n}$, but if we lower our expectations to
the subcategory containing only the fundamental representations, and
their tensor products, then the functor is in fact full. Kuperberg
gave a proof of this fact for $n=3$, by recognising the image of the
functor using a Tannaka-Krein type theorem. This argument continues
to work with only slight modifications for all $n$. I'll also give a
direct proof using quantum Schur-Weyl duality, in \S
\ref{sec:generators-for-fundrep}.

The second question has proved more difficult. Partial answers have
been known for some time. I will describe a new method for
discovering elements of the kernel, based on branching. This method also gives us a limited ability to find obstructions for further relations.

The core of the idea is that there is a forgetful functor $$\GT:\Rep
\uqsl{n} \To \Rep \uqsl{n-1},$$ which forgets the full $\uqsl{n}$
action but does not change the underlying linear maps, and that this
should be reflected somehow in the diagrams. A diagram in
$\PivCat_{n}$ `represents' some morphism in $\Rep \uqsl{n}$;
thinking of this as a morphism in $\Rep \uqsl{n-1}$ via $\GT$, we can hope to
represent it by diagrams in $\PivCat_{n-1}$. This hope is borne out---in
\S \ref{sec:dGT} I construct a functor $\dGT:\PivCat_{n} \To
\Mat(\PivCat_{n-1})$ (and explain what a `matrix category' is), in such
a way that the following diagram commutes:
\begin{equation*}
\xymatrix@R+5mm@C+15mm{
  \PivCat_{n} \ar[d]^{\dGT} \ar[r]^{\Rep_{n}} & \Rep \uqsl{n} \ar[d]^{\GT} \\
  \Mat(\PivCat_{n-1}) \ar[r]^{\Mat(\Rep_{n-1})}           & \Rep \uqsl{n-1}
}
\end{equation*}
With this functor on hand, we can begin determining the kernel of
$\Rep_{n}$. In particular, given a morphism in
$\PivCat_{n}$ (that is, some linear combination of diagrams) we
can consider the image under $\dGT$. This is a matrix of (linear
combinations of) diagrams in $\PivCat_{n-1}$. Then the original diagrammatic morphism becomes
zero in the $\uqsl{n}$ representation theory exactly if each entry
in this matrix of diagrams is zero in the $\uqsl{n-1}$ representation
theory. Thus, if we understand the kernel of $\Rep_{n-1}$, we can
obtain quite strong restrictions on the kernel of $\Rep_{n}$. Of
course, the kernels of $\Rep_2$ and $\Rep_3$ are well known, given
by the relations in the Temperley-Lieb category and Kuperberg's
spider for $\csl{3}$. Moreover, the kernel of $\Rep_1$ is \emph{really} easy to describe. The method described allows us to work up from
these, to obtain relations for all $\uqsl{n}$.

In \S \ref{sec:kernel} I use this approach to find three families of
relations, and to show that these relations are the only ones of
certain types. It remains a conjecture that the proposed relations
are in fact complete.

The first family are the $I=H$ relations, essentially correspond to $6-j$ symbols:
\begin{equation*}
 \mathfig{0.1}{IH/I_abcd_l} = (-1)^{(n+1)a} \mathfig{0.1}{IH/H_abcd_d}.
\end{equation*}
For a given boundary, there are two types of squares, and each can be written as a linear combination of the others. I call these relations the `square-switch' relations. When $n+\ssum{a}-\ssum{b} \geq 0$, for $\max{b} \leq l \leq \min{a}+n$ we have
\begin{align*}
\mathfig{0.17}{PQ/P_square_l} = \sum_{m=\max{a}}^{\min{b}} \qBinomial{n+\ssum{a}-\ssum{b}}{m+l-\ssum{b}} \mathfig{0.17}{PQ/Q_square_m} \\
\intertext{and when $n+\ssum{a}-\ssum{b} \leq 0$, for $\max{a} \leq l \leq \min{b}$ we have}
\mathfig{0.17}{PQ/Q_square_l} = \sum_{m=\max{b}}^{n+\min{a}} \qBinomial{\ssum{b}-n-\ssum{a}}{m+l-\ssum{a}-n} \mathfig{0.17}{PQ/P_square_m}.
\end{align*}
Finally, there are relations amongst polygons of arbitrarily large size (but with a cutoff for each $n$), called the `\Kekule' relations. For each $\ssum{b} \leq j \leq \ssum{a}+n-1$,
\begin{equation*}
 \sum_{k=-\sumhat{b}}^{-\sumtah{a} + 1} (-1)^{j+k} \qBinomial{j+k-\max{b}}{j-\ssum{b}} \qBinomial{\min{a} + n -j -k}{\ssum{a} + n -1 -j} \mathfig{0.2}{PQ/octagon_jk} = 0.
\end{equation*}

\section{The Temperley-Lieb algebras}
The $n=2$ part of this story has, unsurprisingly, been understood
for a long time. The Temperley-Lieb category gives a diagrammatic
presentation of the morphisms between tensor powers of the standard
representation of $\uqsl{2}$. The objects of this category are
natural numbers, and the morphisms from $n$ to $m$ are $\qRing$-linear combinations of diagrams drawn in a horizontal strip
consisting of non-intersecting arcs,
with $n$ arc endpoints on the bottom edge of the strip, and $m$ on
the top edge. (Notice, in particular, that I'm an optimist, not a
pessimist; time goes up the page.) Composition of morphisms is achieved by gluing diagrams, removing each closed circle in exchange for a factor of $\qi{2}$.

\section{Kuperberg's spiders}
The $n=3$ story, dates back to around Kuperberg's paper \cite{q-alg/9712003}. There he defines the notion of a `spider' (in this work, we use the parallel
notion of a pivotal category%
), and constructs the spiders for each of the rank $2$ Lie algebras $A_2 = \mathfrak{su}(3)$, $B_2 = \mathfrak{sp}(4)$ and $G_2$ and their quantum analogues.
Translated into a category, his $A_2$ spider has objects words in $(+,-)$, and morphisms (linear combinations of) oriented trivalent graphs drawn in a horizontal strip, with
orientations of boundary points along the top and bottom edges coinciding with the target and source word objects, and each trivalent vertex either `oriented inwards' or `oriented outwards',
subject to the relations
\begin{align}
\mathfig{0.06}{webs/clockwise_circle} & = \qi{3} = q^2 + 1 + q^{-2} \label{eq:kuperberg-loop}  \\
\mathfig{0.065}{webs/bigon} & = -\qi{2} \mathfig{0.015}{webs/strand} \label{eq:kuperberg-bigon} \\
\intertext{and}
\mathfig{0.125}{webs/sl_3/oriented_square} & = \mathfig{0.125}{webs/sl_3/two_strands_horizontal} + \mathfig{0.125}{webs/sl_3/two_strands_vertical}. \label{eq:kuperberg-square}
\end{align}
(Note that there's a `typo' in Equation (2) of \cite{q-alg/9712003}, corresponding to Equation \eqref{eq:kuperberg-square} above; the $\mathfig{0.05}{webs/sl_3/two_strands_horizontal}$ term
has been replaced by another copy of the $\mathfig{0.05}{webs/sl_3/two_strands_vertical}$ term.)

Kuperberg proves that this category is equivalent to a full subcategory of the category of representations of the quantum group $\uqsl{3}$; the subcategory with objects arbitrary tensor products of the two $3$-dimensional representations.
It is essentially an equinumeration proof, showing that the number of diagrams (modulo the above relations) with a given boundary agrees with the dimension of the appropriate $\uqsl{3}$ invariant space.
I'm unable to give an analogous equinumeration argument in what follows.

Note that the $n=3$ special case of my construction will not quite reproduce Kuperberg's relations above; the bigon relation will involve a $+\qi{2}$, not a $-\qi{2}$. This is just a normalisation issue, resolved by multiplying each vertex by $\sqrt{-1}$.

\chapter{The `diagrammatic' category $\SymCat_n$}
\label{sec:diagrams}%
Just as permutations form groups, planar diagrams up to planar
isotopy form pivotal categories. In what follows, we'll define a
certain `free (strict) pivotal category', $\PivCat_n$ along with a
slight modification called $\SymCat_n$ obtained by adding some
symmetries and some relations for degenerate cases. Essentially,
$\PivCat_n$ will be the category of trivalent graphs, with edges
carrying both orientations and labels $1$ through $n-1$, up to planar
isotopy.

For lack of a better place, I'll introduce the notion of a matrix
category here; given any category $\mathcal{C}$ in which the
$\operatorname{Hom}$ spaces are guaranteed to be abelian groups, we
can form a new category $\Mat[\mathcal{C}]$, whose objects are
formal finite direct sums of objects in $\mathcal{C}$, and whose
morphisms are matrices of appropriate morphisms in $\mathcal{C}$.
Composition of morphisms is just matrix multiplication (here's where
we need the abelian group structure on $\operatorname{Hom}$ spaces).
If the category already had direct sums, then there's a natural
isomorphism $\mathcal{C} \Iso \Mat[\mathcal{C}]$.

\section{Pivotal categories}
\label{sec:pivotal}%
I'll use the formalism of pivotal categories in the following. This
formalism is essentially interchangeable with that of spiders, due
to Kuperberg \cite{q-alg/9712003}, or of planar algebras, due to Jones
\cite{math.QA/9909027}.\footnote{These alternatives are perhaps more
desirable, as the pivotal category view forces us to make a
artificial distinction between the domain and codomain of a
morphism. I'll have to keep reminding you this distinction doesn't matter, in what follows. On the other hand, the categorical setup allows us to more
easily incorporate the notion of direct sum.}

A pivotal category is a monoidal category\footnote{For our purposes,
we need only consider strict monoidal categories, where, amongst
other things, the tensor product is associative on the nose, not
just up to an isomorphism. The definitions given here must be
modified for non-strict monoidal categories.} $\mathcal{C}$ equipped
with
\begin{enumerate}
\item a cofunctor $*:\mathcal{C}^{op} \To \mathcal{C}$, called the
dual,
\item a natural isomorphism $\tau:1_{\mathcal{C}} \To **$,
\item a natural isomorphism $\gamma:\tensor \compose (* \times *)
\To * \compose \tensor^{op}$,
\item an isomorphism $e \To e^*$, where $e$ is the neutral object
for tensor product, and
\item for each object $c \in \mathcal{C}$, a `pairing' morphism
$p_c : c^* \tensor c \To e$.
\end{enumerate}
with the natural isomorphisms satisfying certain coherence
conditions, and the pairing morphisms certain axioms, all given in
\cite{MR1686423}.\footnote{Actually, \cite{MR1020583} is an earlier
reference for the strict case (see below), and they cite a preprint
of \cite{MR1250465} for the full version; however, the published
version of that paper ended up introducing a slightly different
notion, of `autonomy' for duals.} The primary example of a pivotal
category is the category of representations of an involutory or ribbon Hopf algebra
\cite{MR1686423}.
A pivotal functor between two pivotal categories should intertwine
the duality cofunctors, commute with the isomorphisms $e \To e^*$,
and take pairing morphisms to pairing morphisms. It should also
intertwine the natural isomorphisms $\tau$; that is,
for a functor $\mathcal{F}:\mathcal{C} \To \mathcal{D}$, we require
$\tau^{\mathcal{D}}_{\mathcal{F}(a)} =
\mathcal{F}(\tau^{\mathcal{C}}_a)$. There's a similar condition for
$\gamma$.

In the case that $\tau$ and $\gamma$ are simply the identity on each
object (and $e=e^*$), we say the pivotal category is strict.
Unfortunately, representations of a Hopf algebra do not generally
form a strict pivotal category; $\tau$ cannot be the identity. However, every pivotal category is
equivalent to a strict pivotal category \cite{MR1686423}. We'll
take advantage of this later!

In a pivotal category with $e=e^*$ and $\gamma$ the identity (but
not necessarily also $\tau$), the axioms satisfied by the natural
isomorphism $\tau$ simplify to the requirements that it's a tensor
natural transformation: $\tau_{a \tensor b} = \tau_a \tensor
\tau_b$, and that $\tau_a^* : a^{***} \To a^*$ and $\tau_{a^*} : a^*
\To a^{***}$ are inverse morphisms.\footnote{This condition might also
be stated as $\tau$ commuting with the cofunctor $*$.} In this
same situation, the axioms for the pairing morphisms simplify to:
\begin{enumerate}
\item For each object $a$,
\begin{equation}
\label{eq:pivotal-straightening}%
\mathfig{0.135}{morphisms/p_cp} = \mathfig{0.036}{morphisms/id_a_star}.
\end{equation}
\item For each morphism $f:a \To b$,
\begin{equation}
\label{eq:duals-by-pairing-maps}%
\mathfig{0.1}{morphisms/pairing_f} = \mathfig{0.1}{morphisms/f_pairing}.
\end{equation}
\item For all objects $a$ and $b$,
\begin{equation}
\label{eq:duals-ab}
\mathfig{0.15}{morphisms/pairing_ab} = \mathfig{0.15}{morphisms/pairing_ab_2}.
\end{equation}
\end{enumerate}
From these, we can derive
\begin{lem}
\label{lem:dual-formula}%
The dual of a morphism can be written in terms of
the original morphism, components of $\tau$, and the pairing morphisms as
\begin{equation*}
\mathfig{0.04}{morphisms/f_star} = \mathfig{0.1}{morphisms/f_star_formula}
\end{equation*}
\end{lem}
\begin{proof}
\begin{align*}
\mathfig{0.1}{morphisms/f_star_formula} & = \mathfig{0.1}{morphisms/f_star_formula_1} \\
\intertext{using by the naturality of $\tau$, which becomes}
    & = \mathfig{0.1}{morphisms/f_star_formula_2} \\
\intertext{by Equation \eqref{eq:duals-by-pairing-maps}, and finally}
    & = \mathfig{0.04}{morphisms/f_star}
\end{align*}
by Equation \eqref{eq:pivotal-straightening}.
\end{proof}

Note that the category of matrices over a pivotal category is still pivotal, in an essentially obvious way.

\section{Quotients of a free tensor category}
\label{sec:quotients}%
To begin with, let's just define a free (strict) monoidal
category on some generating morphisms, which we'll call
$\FreeCat_n$. We'll then add some relations implementing planar
isotopy to obtain $\PivCat_n$, and some more relations to obtain
$\SymCat_n$.
The objects of $\FreeCat_n$ form a monoid under tensor product, with
neutral object $0$ (sometimes also called $n$), generated by the set $\{1, \dotsc, n-1, 1^*,
\dotsc, (n-1)^*\}$.
The `generating morphisms' are diagrams
\begin{equation*}
\mathfig{0.072}{pairings/pl_a}, \qquad \mathfig{0.072}{pairings/pr_a}, \qquad \mathfig{0.072}{pairings/cpl_a}, \qquad \mathfig{0.072}{pairings/cpr_a}
\end{equation*}
for each $a = 1, \dotsc, n-1$ (but not for the `dual integers' $1^*,\dotsc,(n-1)^*$), along with
\begin{equation}
\label{eq:duals}
\mathfig{0.078}{duals/dtr}, \qquad
\mathfig{0.078}{duals/dal}, \qquad
\mathfig{0.078}{duals/dtl}, \qquad \text{and} \qquad
\mathfig{0.078}{duals/dar}
\end{equation}
again for each $a = 1, \dotsc, n-1$, and finally
\begin{equation}
\label{eq:vertices}
\mathfig{0.1}{vertices/v_outup_0} \qquad \text{and} \qquad
\mathfig{0.1}{vertices/v_inup_0}
\end{equation}
for each $a,b,c = 0, 1, \dotsc, n-1$ such that $a+b+c=n$. We'll
sometimes speak of the `type' of a vertex; the first, outgoing,
vertex here is of $+$-type, the second, incoming, vertex is of
$-$-type. We say a vertex is `degenerate' if one of its edges is
labelled with $0$. Notice there are no nondegenerate trivalent
vertices for $n=2$, and exactly one of each type for $n=3$. To read
off the source of such a morphism, you read across the lower
boundary of the diagram; each endpoint of an arc labelled $a$ gives
a tensor factor of the source object, either $a$ if the arc is
oriented upwards, or $a^*$ if the arc is oriented downwards. To read
the target, simply read across the upper boundary. Thus the source
of $\mathfig{0.05}{vertices/v_outup_0}$ is $0$, and the target is $a
\tensor b \tensor c$. All morphisms are then generated from these,
by formal tensor product and composition, subject only to the usual
identities of a tensor category.

Next $\PivCat_n$. This category has exactly the same objects. The
morphisms, however, are arbitrary trivalent graphs drawn in a strip,
which look locally like one of the pictures above, up to planar
isotopy fixing the boundary of the strip. (Any boundary points of
the graph must lie on the boundary of the strip.) Thus the graphs
are oriented, with each edge carrying a label $1$ through $n$, and
edges only ever meet bivalently as in Equation \eqref{eq:duals} or
trivalently as in Equation \eqref{eq:vertices}. Being a little more
careful, we should ask that the diagrams have product structure near
the boundary, that this is preserved throughout the isotopies, and
that small discs around the trivalent vertices are carried around
rigidly by the isotopies, so we can always see the ordering (not
just the cyclic ordering) of the three edges incident at a vertex.
(Note, though, that we're not excited about being able to see this
ordering; we're going to quotient it out in a moment.) The source
and target of such a graph can be read off from the graph exactly as
described for the generators of $\FreeCat_n$ above.

Next, we make this category into a strict pivotal category. For
this, we need to define a duality functor, specify the evaluation
morphisms, and then check the axioms of \S \ref{sec:pivotal}. The
duality functor on objects is defined by $0^* = 0$ and otherwise
$(k)^* = k^*, (k^*)^* = k$. On morphisms, it's a $\pi$ rotation of
the strip the graph is drawn in. Clearly the double dual functor $**$ is the identity on the nose. The evaluation morphisms for $a =
1, \dotsc, n-1$ are `leftwards-oriented' cap diagrams; the
evaluation morphisms for $a = 1^*, \dotsc, n-1^*$ are
`rightwards-oriented'. The axioms in Equation \eqref{eq:pivotal-straightening} and \eqref{eq:duals-by-pairing-maps} are then satisfied automatically, because we allow isotopy of diagrams. The evaluation morphisms for iterated tensor
products are just the nested cap diagrams, with the unique
orientations and labels matching the required source object. This definition ensures the axiom of Equation \eqref{eq:duals-ab} is satisfied.

There's a (tensor) functor from $\FreeCat_n$ to $\PivCat_n$, which I'll call $\Draw$. Simply take a morphism in $\FreeCat_n$, which can
be written as a composition of tensor products of generating morphisms, and \emph{draw} the corresponding diagram, using the usual rules
of stacking boxes to represent composition, and juxtaposing boxes side by side to represent tensor product.
The resulting diagram can then be interpreted as a morphism in $\PivCat_n$.
This is well-defined by the usual nonsense of \cite{MR1113284}, that
the identities relating tensor product and composition in a tensor
category correspond to `rigid' isotopies (that is, isotopies which
do not rotate boxes). The functor is obviously full; or at least,
obviously modulo some Morse theory. The kernel of this functor is
generated by the extra isotopies we allow in $\PivCat_n$. Thus, as a
tensor ideal of $\FreeCat_n$, $\ker{\Draw}$ is generated by
\begin{align}
\label{eq:rotate-vertex-2pi}
\mathfig{0.27}{isotopy/v_in_2pi_rotated}  & = \mathfig{0.11}{isotopy/v_in} &
\mathfig{0.27}{isotopy/v_out_2pi_rotated} & = \mathfig{0.11}{isotopy/v_out}, \\
\label{eq:straighten-strand}
\mathfig{0.084}{isotopy/strand_left} & = \mathfig{0.024}{isotopy/strand} &
\mathfig{0.024}{isotopy/strand} & = \mathfig{0.084}{isotopy/strand_right}, \\
\notag
\mathfig{0.084}{isotopy/strand_down_left} & = \mathfig{0.024}{isotopy/strand_down} &
\mathfig{0.024}{isotopy/strand_down} & = \mathfig{0.084}{isotopy/strand_down_right}, \\
\label{eq:tag-near-maxima}
\intertext{and}
\mathfig{0.1}{pairings/cp_tag_left} & = \mathfig{0.1}{pairings/cp_tag_right} &
\mathfig{0.1}{pairings/p_tag_left} & = \mathfig{0.1}{pairings/p_tag_right}.
\end{align}
There are other obvious variations of Equations
\eqref{eq:rotate-vertex-2pi} (rotating the vertex other way) and
\eqref{eq:tag-near-maxima} (tags pointing the other way), but these
follow easily from the ones given here. Whenever we want to define a
functor on $\PivCat_n$ by defining it on generators, we need to
check these morphisms are in the kernel.

Finally, we can define the category we're really interested in,
which I'll call $\SymCat_n$. We'll add just a few more relations to
$\PivCat_n$; these will be motivated shortly when we define a
functor from $\FreeCat_n$ to the representations category of $\uqsl{n}$. This functor will
descend to the quotient $\PivCat_n$, and then to the quotient
$\SymCat_n$, and the relations we add from $\PivCat_n$ to
$\SymCat_n$ will be precisely the parts of the kernel of this
functor which only involve a single generator. The real work of this
thesis is, of course, understanding the rest of that kernel! We add
relations insisting that the trivalent vertices are rotationally
symmetric
\begin{align}
\label{eq:rotate-vertex}
 \mathfig{0.16}{isotopy/v_in_rotated} & = \mathfig{0.105}{isotopy/v_in}, & \mathfig{0.16}{isotopy/v_out_rotated} & = \mathfig{0.105}{isotopy/v_out},
\end{align}
that opposite tags cancel
\begin{equation}
\label{eq:tags-cancel}
\mathfig{0.052}{duals/two_tags} = \mathfig{0.026}{isotopy/strand},
\end{equation}
that dual of a tag is a $\pm 1$ multiple of a tag
\begin{align}
\label{eq:dual-of-a-tag}
 \mathfig{0.078}{duals/dtr} & = (-1)^{(n+1)a} \mathfig{0.078}{duals/dtl}, & \mathfig{0.078}{duals/dar} & = (-1)^{(n+1)a} \mathfig{0.078}{duals/dal},
\end{align}
and that trivalent vertices `degenerate' to tags
\begin{align}
\label{eq:degeneration}
 \mathfig{0.105}{vertices/v_inup_degenerate}  & = \mathfig{0.105}{pairings/p_tag_left} &
 \mathfig{0.105}{vertices/v_outup_degenerate} & = \mathfig{0.105}{pairings/cp_tag_right}.
\end{align}
Notice here we're implicitly using the canonical identifications between the objects $0 \tensor a$, $a$, and $a \tensor 0$, available because our
tensor categories are strict.

Clearly the element of $\ker{\Draw}$ in Equation
\eqref{eq:rotate-vertex-2pi} can be constructed by tensor product
and composition out of the briefer rotations in Equation
\eqref{eq:rotate-vertex}, and so in checking the well-definedness of
a functor on $\SymCat_n$, we only need to worry about the latter.

\section{Flow vertices}
\label{sec:flow-vertices}%
We'll now introduce two new types of vertices. You could add them as
diagrammatic generators, then impose as relations the formulas
below, but it's less cumbersome to just think of them as a
convenient notation. In each of these vertices, there will be some
`incoming' and some 'outgoing' edges, and the sum of the incoming
edges will be the same as the sum of the outgoing edges.
\begin{defn}
The `flow vertices' are
\begin{align*}
\mathfig{0.12}{vertices/v_merge} & = \mathfig{0.12}{vertices/v_merge_def} \\
\intertext{and}
\mathfig{0.12}{vertices/v_split} & = \mathfig{0.12}{vertices/v_split_def}.
\end{align*}
\end{defn}
The convention here is that the `hidden tag' lies on the `thick' edge, and points counterclockwise. These extra vertices will be convenient in what
follows, hiding a profusion of tags. `Splitting' vertices are of $+$-type, `merging' vertices are of $-$-type.

\section{Polygonal webs}
\label{sec:polygons}
To specify the kernel of the representation functor, in \S \ref{sec:kernel}, we'll need to introduce some
notations for `polygonal webs'. These webs will come in two families,
the `$\mathcal{P}$' family and the `$\mathcal{Q}$' family. In each family, the vertices around the polygon will alternate in type. A boundary edge which is
connected to a $+$-vertex in a $\mathcal{P}$-polygon will be connected to a $-$-vertex in a $\mathcal{Q}$-polygon.

For $a,b \in \Integer^k$ define\footnote{I realise this definition is `backwards', or at least easier to read from right to left than from left to right. Sorry---I only realised too late.} the boundary label pattern
$$\Label(a,b) =
(b_k-a_1,-) \tensor (b_k-a_k,+) \tensor \dotsb \tensor (b_2-a_2,+) \tensor (b_1-a_2,-) \tensor (b_1-a_1,+).$$
We'll now define some elements of $\Hom{\SymCat^n}{\eset}{\Label(a,b)}$,
$\P{n}{a,b}{l}$ for $\max{b} \leq l \leq \min{a} + n$ and
$\Q{n}{a,b}{l}$ for $\max{a} \leq l \leq \min{b}$ by
\begin{align}
\P{n}{a,b}{l} & = \mathfig{0.65}{PQ/P_flows} \notag \\
              & = \mathfig{0.65}{PQ/P_labels} \label{eq:P-defn} \\
\intertext{and}
\Q{n}{a,b}{l} & = \mathfig{0.65}{PQ/Q_flows} \notag \\
              & = \mathfig{0.65}{PQ/Q_labels}. \label{eq:Q-defn}
\end{align}
(The diagrams are for $k=3$, but you should understand the obvious generalisation for any $k \in \Natural$.)
You should consider the first of each pair of diagrams simply as notation for the second. Each edge label is
a signed sum of the `flows labels' on either side, determining signs by relative orientations. It's trivial\footnote{Actually, perhaps only trivial after acknowledging that the disk in which the diagrams are drawn is simply connected.} to
see that for web diagrams with only `2 in, 1 out' and `1 in, 2 out'
vertices, it's always possible to pick a set of flow labels
corresponding to an allowable set of edge labels. Not every set of
flow labels, however, gives admissible edge labels, because the edge
labels must be between $0$ and $n$. The allowable flow labels for the $\mathcal{P}$- and $\mathcal{Q}$-polygons are
exactly those for which $a_i, a_{i+1} \leq b_i \leq n+a_i,
n+a_{i+1}$. 
Further, there's a $\Integer$ redundancy in flow labels; adding a
constant to every flow label in a diagram doesn't actually change
anything. Taking this into account, there's a finite set of pairs $a,b$ for each $n$ and $k$. The inequalities on the `internal flow label' $l$ for both
$\P{n}{a,b}{l}$ and $\Q{n}{a,b}{l}$ simply demand that all the
internal edges have labels between $0$ and $n$, inclusive.

We denote the subspace of $\Hom{\SymCat^n}{\eset}{\Label(a,b)}$ spanned
by all the $\mathcal{P}$-type polygons by $\AP^{n}_{a,b}$, and the
subspace spanned by the $\mathcal{Q}$-type polygons by
$\AQ^{n}_{a,b}$. The space $\AP^{n}_{a,b}$ is $\min{a} - \max{b} + n
- 1$ dimensional (or $0$ dimensional when this quantity is
negative), and the space $\AQ^{n}_{a,b}$ is $\max{a} - \min{b} + 1$
dimensional.

A word of warning; $a$ and $b$ each having $k$ elements does not
necessarily mean that $\P{n}{a,b}{l}$ or $\Q{n}{a,b}{l}$ are honest
$2k$-gons. This can fail in two ways. First of all, if some $a_i =
b_i$, $a_i + n = b_i$, $a_{i+1} = b_i$ or $a_{i+1} + n = b_i$, then
one of the external
edges carries a trivial label.
Further, when $l$ takes on one of its extremal allowed values, at
least one of the internal edges of the polygon becomes trivial, and
the web becomes a tree, or a disjoint union of trees and arcs. For
example (ignoring the distinction between source and target of
morphisms; strictly speaking these should all be drawn with all
boundary points at the top of the diagram),
\begin{align*}
\P{4}{(0,0,0),(1,1,1)}{1} & = \mathfig{0.1}{webs/sl_4/kekule0} &
\P{4}{(0,0,0),(1,1,1)}{2} & = \mathfig{0.1}{webs/sl_4/kekule1} \\
\P{4}{(0,0,0),(1,1,1)}{3} & = \mathfig{0.1}{webs/sl_4/kekule2} &
\P{4}{(0,0,0),(1,1,1)}{4} & = \mathfig{0.1}{webs/sl_4/kekule3} \\
\intertext{and}
\P{5}{(0, 1, 1),(2, 2, 2)}{2} & = \mathfig{0.1}{webs/sl_5/P_011_222_2} &
\P{5}{(0, 1, 1),(2, 2, 2)}{3} & = \mathfig{0.1}{webs/sl_5/P_011_222_3} \\
\P{5}{(0, 1, 1),(2, 2, 2)}{4} & = \mathfig{0.1}{webs/sl_5/P_011_222_4} &
\P{5}{(0, 1, 1),(2, 2, 2)}{5} & = \mathfig{0.1}{webs/sl_5/P_011_222_5}
\end{align*}

Finally, it's actually possible for a $\mathcal{P}$-type polygon and
a $\mathcal{Q}$-type polygon to be equal in $\SymCat_n$. This can only
happen in the case that $a$ and $b$ each have length $2$, or at any length, when either $a$ or $b$ is constant. This only involves polygons with extreme values of the
internal flow label $l$. Specifically
\begin{lem}
\label{lem:P=Q}%
If $a$ and $b$ are each of length $2$,
\begin{align*}
\P{n}{a,b}{\max{b}} & = \Q{n}{a,b}{\min{b}} \\
\intertext{and}
\P{n}{a,b}{\min{a}+n} & = \Q{n}{a,b}{\max{a}}.
\end{align*}
Further, even if $a$ and $b$ have length greater than $2$, when $a$ is a constant vector $a = \vect{a}$
$$\P{n}{\vect{a},b}{a+n} =  \Q{n}{\vect{a},b}{a}$$
and when $b$ is a constant vector $b = \vect{b}$
$$\P{n}{a,\vect{b}}{b} =  \Q{n}{a,\vect{b}}{b}$$
in $\SymCat$. Otherwise, the $\mathcal{P}$- and $\mathcal{Q}$-polygons
are linearly independent in $\SymCat_n$. In particular, $\dim(\AP^n_{a,b} \cap
\AQ^n_{a,b})$ is $0$, $1$ or $2$ dimensional, depending on whether
neither $a$ nor $b$ are constant, one is, or either both are or $a$ and $b$ have length $2$.
\end{lem}
For example
\begin{align*}
\P{3}{(0,0),(1,1)}{1} & = \mathfig{0.1}{webs/sl_3/two_strands_horizontal} = \Q{3}{(0,0),(1,1)}{1}, \\
\P{3}{(0,0),(1,1)}{2} & = \mathfig{0.1}{webs/sl_3/oriented_square}, \\
\intertext{and}
\P{3}{(0,0),(1,1)}{3} & = \mathfig{0.1}{webs/sl_3/two_strands_vertical} = \Q{3}{(0,0),(1,1)}{0}.
\end{align*}

\subsection{Rotations}
\label{sec:rotations}%
It's important to point out that the two types of polygons described
above are actually closely related. In fact, by a `rotation', and
adding some tags, we can get from one to the other. This will be
important in our later descriptions of the kernel of the
representation functor. Unsurprisingly, this symmetry between the
types of polygons is reflected in the kernel, and we'll save half
the effort in each proof, essentially only having to deal with one
of the two types.

First, by $\rotl:\Integer^k \To \Integer^k$ we just mean `rotate left': $$\rotl(a_1,a_2,\dotsc,a_k) = (a_2,\dotsc,a_k,a_1),$$ and by $\rotr:\Integer^k \To \Integer^k$ `rotate right'.
\begin{lem}
Writing the identity out by explicit tensor products and compositions would be tedious; much easier is to write it diagrammatically first off:
\begin{equation}
\label{eq:rotations}
 \mathfig{0.32}{PQ/P_rotl} = \mathfig{0.44}{PQ/Q_rotated},
\end{equation}
although keep in mind that we intend the obvious generalisation to
arbitrary size polygons, not just squares.
\end{lem}
\begin{proof}
It's just a matter of using isotopy, the definition of the `flow vertices' in \S \ref{sec:flow-vertices} and cancelling tags by Equation \eqref{eq:tags-cancel}.
\begin{align*}
\mathfig{0.28}{PQ/rotation_1} & = \mathfig{0.28}{PQ/rotation_2} = \mathfig{0.24}{PQ/rotation_3} \\
                            & = \mathfig{0.24}{PQ/rotation_4} = \mathfig{0.24}{PQ/rotation_5}. \qedhere
\end{align*}
\end{proof}
Notice that it's actually
not important which way the tags point in Equation \eqref{eq:rotations}; we could reverse them, since by Equation
\eqref{eq:dual-of-a-tag} we'd just pick up an overall sign of
$(-1)^{(n+1)(2\ssum{b}-2\ssum{a})} = 1$. We'll sometimes write Equation \eqref{eq:rotations}
as
\begin{align*}
\P{n}{b-n,\rotl[a]}{l} & = \drotl[\Q{n}{a,b}{l}] \\
\intertext{or equivalently, re-indexing}
\P{n}{a,b}{l} & = \drotl[\Q{n}{\rotr(b),a+n}{l}]
\end{align*}
where the operation $\drotl$ is implicitly defined by Equation
\eqref{eq:rotations}. (You should read $\drotl$ as `diagrammatic
rotation', and perhaps have the `$d$' prefix remind you that we're
not just rotating, but also adding a tag, called $d$ in the
representation theory, to each external edge.) The corresponding
identity writing a $\mathcal{Q}$-polygon in terms of a $\mathcal{P}$-polygon is
\begin{align*}
\Q{n}{\rotr(b),a+n}{l} & = \drotl[\P{n}{a,b}{l}] \\
\intertext{or}
\Q{n}{a,b}{l}& = \drotl[\P{n}{b-n,\rotl(a)}{l}].
\end{align*}

\chapter{Just enough representation theory}
In this chapter, we'll describe quite a bit of representation
theory, but hopefully only what's necessary for later! There are no
new results in this section, although possibly Proposition
\ref{prop:elementary-morphisms-generate}, explaining roughly that
`the fundamental representation theory of $\uqsl{n}$ is generated by
triple invariants' has never really been written down.

\section{The Lie algebra $\csl{n}$}
The Lie algebra $\csl{n}$ of traceless, $n$-by-$n$ complex matrices
contains the commutative Cartan subalgebra $\csa$ of diagonal
matrices, spanned by $H_i = E_{i,i}-E_{i+1,i+1}$ for
$i=1,\ldots,n-1$. (Here $E_{i,j}$ is simply the matrix with a single
nonzero entry, a $1$ in the $(i,j)$ position.) When $\csl{n}$ acts
on a vector space $V$, the action of $\csa$ splits $V$ into
eigenspaces called weight spaces, with eigenvalues in $\csa^*$. The
points
$$\Lambda = \setc{\lambda \in \csa^*}{\lambda(H_i) \in \Integer}
\Iso \Integer^{n-1}$$ form the weight lattice. We call $\csa^*_+ =
\setc{\lambda \in \csa^*}{\lambda(H_i) \geq 0}$ the positive Weyl
chamber, and the lattice points in it, $\Lambda_+ = \Lambda \cap
\csa^*_+$, the dominant weights. The positive Weyl chamber looks
like $\Real_{\geq 0}^{n-1}$, and the dominant weights like
$\Natural^{n-1}$, generated by certain fundamental weights,
$\lambda_1, \ldots, \lambda_{n-1}$, the dual basis to $\set{H_i}
\subset \csa$.

Under the adjoint action of $\csl{n}$ on itself, $\csl{n}$ splits
into $\csa$, as the $0$-weight space, and one dimensional root
spaces, each spanned by a root vector $E_{i \neq j}$ with weight
$[H_k, E_{i,j}] =
(\delta_{i,k}+\delta_{j,k+1}-\delta_{i,k+1}-\delta_{j,k}) E_{i,j}$.
In particular, the simple roots are defined $E_i^+ = E_{i,i+1}$, for
$i=1,\ldots,n-1$. Similarly define $E_i^- = E_{i+1,i}$. Together,
$E_i^+, E_i^-$ and $H_i$ for $i=1,\ldots,n-1$ generate $\csl{n}$ as
a Lie algebra.

The universal enveloping algebra $\usl{n}$ is a Hopf algebra
generated as an associative algebra by symbols $E_i^+, E_i^-$ and
$H_i$ for $i=1,\ldots,n-1$, subject only to the relations that the
commutator of two symbols agrees with the Lie bracket in $\csl{n}$. (We won't bother specify the other Hopf algebra structure, the comultiplication, counit or antipode. See the next section for all the details for quantum $\csl{n}$.)
Trivially, $\csl{n}$ and $\usl{n}$ have the same representation
theory. We denote the subalgebra generated by $H_i$ and $E_i^\pm$ for
$i=1,\ldots,n-1$ by $U^\pm\,\left(\csl{n}\right)$.

\section{The quantum groups $\uqsl{n}$}
We now recall the $q$-deformation of the Hopf algebra $\usl{n}$.
Unfortunately we can't straightforwardly recover $\usl{n}$ by
setting $q=1$, but we will see in \S \ref{ssec:representations} that
this is the case at the level of representation theory.

Let $\A = \Rational(q)$ be the field of rational functions in an
indeterminate $q$, with coefficients in $\Rational$. The quantum
group $\uqsl{n}$ is a Hopf algebra over $\A$, generated as an
associative algebra by symbols $X_i^+, X_i^-, K_i$ and $K_i^{-1}$,
for $i = 1,\ldots, n-1$, subject to the relations
\begin{alignat*}{2}
K_i K_i^{-1} & = 1 = K_i^{-1} K_i & \\
K_i K_j & = K_j K_i & \\
K_i X_j^\pm & = q^{\pm(i,j)} X_j^\pm K_i & \qquad & \text{where
$(i,j) = \begin{cases}2 & \text{ if $i=j$,} \\
-1 & \text{ if $\abs{i-j} = 1$, or} \\ 0 & \text{
if $\abs{i-j} \geq 2$}\end{cases}$ } \\
 X_i^+ X_j^- - X_j^- X_i^+ & = \delta_{i j} \frac{K_i - K_i^{-1}}{q-q^{-1}}\\
 X_i^\pm X_j^\pm & = X_j^\pm X_i^\pm && \text{when $\abs{i-j} \geq 2$}
\end{alignat*}
and the quantum Serre relations
\begin{equation*}
 {X_i^\pm}^2 X_j^\pm - \qi{2} X_i^\pm X_j^\pm X_i^\pm + X_j^\pm
 {X_i^\pm}^2 = 0 \qquad \text{ when $\abs{i-j} = 1$.}
\end{equation*}
The Hopf algebra structure is given by the following formulas for
the comultiplication $\comult$, counit $\counit$ and antipode
$\antipode$.
\begin{align*}
\comult(K_i) & = K_i \tensor K_i, \\
\comult(X_i^+) & = X_i^+ \tensor K_i + 1 \tensor X_i^+, \\
\comult(X_i^-) & = X_i^- \tensor 1 + K_i^{-1} \tensor X_i^-,
\end{align*}
\begin{align*}
\counit(K_i) & = 1 & \counit(X_i^\pm) & = 0,
\end{align*}
\begin{align*}
\antipode(K_i) & = K_i^{-1} & \antipode(X_i^+) & = - X_i^+ K_i^{-1}
& \antipode(X_i^i) & = - K_i X_i^-.
\end{align*}
(See, for example, \cite[\S 9.1]{MR1300632} for verification that
these do indeed fit together to form a Hopf algebra.) There are
Hopf subalgebras $\uqsl{n}^\pm$, leaving out the
generators $X_i^\mp$.

\section{Representations}\label{ssec:representations}
The finite-dimensional irreducible representations of both $\usl{n}$
and $\uqsl{n}$ can be succinctly catalogued.

\hyphenation{eigen-spaces}%
A representation $V$ of $\usl{n}$ is said
to be a highest weight representation if there is a weight vector $v
\in V$ such that $V = \usl{n}^-\,(v)$. The representation is said to
have weight $\lambda$ if $v$ has weight $\lambda$. As with
$\csl{n}$, representations of $\uqsl{n}$ split into the eigenspaces
of the action of the commutative subalgebra generated by the $K_i$
and $K_i^{-1}$. These eigenspaces are again called weight spaces. A
representation $V$ of $\uqsl{n}$ is called a high weight
representation if it contains some weight vector $v$ so $V =
\uqsl{n}^{-}\,(v)$. The finite-dimensional irreducible
representations (henceforth called irreps) are then classified by:

\begin{prop}
Every irrep of $\usl{n}$ or of $\uqsl{n}$ is a highest weight
representation, and there is precisely one irrep with weight
$\lambda$ for each $\lambda \in \Lambda_+$, the set of dominant
weights.
\end{prop}
\begin{proof}
See \cite[\S 23]{MR1153249} for the classical case, \cite[\S 10.1]{MR1300632} for the quantum case.
\end{proof}

Further, the finite dimensional representation theories are tensor
categories and the tensor product of two irreps decomposes uniquely as
a direct sum of other irreps. The combinatorics of the the tensor
product structures in $\Rep \usl{n}$ and $\Rep \uqsl{n}$ agree. A
nice way to say this is that the Grothendieck rings of $\Rep
\usl{n}$ and $\Rep \uqsl{n}$ are isomorphic, identifying irreps with
the same highest weight \cite[\S 10.1]{MR1300632}.

\subsection{Fundamental representations}
\label{sec:fundamental}%
Amongst the finite dimensional irreducible representations, there
are some particularly simple ones, whose highest weights are the
fundamental weights. These
are called the fundamental representations, and there are $n-1$ of
them for $\usl{n}$ or $\uqsl{n}$.

We'll write, in either case, $\V{a}{n}$ for the fundamental
representation with weight $\lambda_a$. For $\usl{n}$ this is just
the representation $\alt^a \Complex^n$. It will be quite convenient
to agree that $\V{0}{n}$ and $\V{n}{n}$ (poor notation, I admit!)
both denote the trivial representation.

We can now define $\FundRep \uqsl{n}$; it's the full subcategory of
$\Rep \uqsl{n}$, whose objects are generated by tensor product and
duality from the trivial representation $\A$, and the fundamental
representations $\V{a}{n}$. Note that there are no direct sums in
$\FundRep \uqsl{n}$.

All dominant weights are additively generated by fundamental
weights, and this is reflected in the representation theory; the
irrep with high weight $\lambda = \sum_a m_a \lambda_a$ is
contained, with multiplicity one, as a direct summand in the tensor
product $\tensor_a \V{a}{n}^{\tensor m_a}$.

This observation explains, to some extent, why it's satisfactory to
simply study $\FundRep \uqsl{n}$, instead of the full representation
theory. Every irrep, while not necessarily allowed as an object of
$\FundRep \uqsl{n}$, reappears in the Karoubi envelope
\cite{MR0358759, wiki:Karoubi-Envelope}, since any irrep appears as
a subrepresentation of some tensor product of fundamental
representations. In fact, there's a canonical equivalence of
categories $$\Kar[\FundRep \uqsl{n}] \Iso \Rep \uqsl{n}.$$

\subsection{The Gel`fand-Tsetlin basis}
\label{sec:gelfand-tsetlin}
We'll now define the Gel`fand-Tsetlin basis \cite{MR0034763,MR0035774}, 
a canonical basis which arises from the nice multiplicity free
splitting rules for $\uqsl{n} \Into \uqsl{n+1}$. (Here, and
hereafter, this inclusion is just $K_i \mapsto K_i, X_i^\pm \mapsto
X_i^\pm$.)

\begin{lem}
If $V$ is an irrep of $\uqsl{n+1}$, then as a representation of
$\uqsl{n}$ each irrep appearing in $V$ appears exactly once.

More specifically, if $V$ is the $\uqsl{n+1}$ irrep with highest
weight $(\lambda_1, \lambda_2, \ldots, \lambda_{n+1})$, then an
$\uqsl{n}$ irrep $W$ of weight $(\mu_1, \mu_2, \ldots, \mu_{n})$
appears (with multiplicity one) if and only if `$\mu$ fits inside
$\lambda$', that is $$\lambda_1 \leq \mu_1 \leq \lambda_2 \leq \mu_2
\leq \dotsb \leq \lambda_n \leq \mu_n \leq \lambda_{n+1}.$$
In particular, the fundamental irreps $\V{a}{n}$ break up as
$$\V{a}{n} \Iso \V{a-1}{n-1} \directSum \V{a}{n-1}.$$
\end{lem}
\begin{proof}
See \cite[\S 14.1.A]{MR1300632}, \cite{knutson-GT-notes}.
\end{proof}

This allows us to inductively define ordered bases for $\uqsl{n}$
irreps, at least projectively. This was first done in the quantum
case in \cite{MR841713}. 
Choose, without much effort, a basis for the only irreducible representation of
$\uqsl{1}$, the trivial representation $\Complex(q)$. Now for any
representation $V$ of $\uqsl{n+1}$, decompose $V$ over $\uqsl{n}$,
as $V \Iso \DirectSum_{\alpha} W_{\alpha}$, ordered
lexicographically by highest weight, and define the Gel`fand-Tsetlin
basis of $V$ to be the concatenation of the inclusions of the bases
for each $W_{\alpha}$ into $V$. These inclusions are unique up to
complex multiples, and we thus obtain a canonical projective basis.

We'll call this forgetful functor the `Gel`fand-Tsetlin' functor,
$\GT : \Rep\uqsl{n} \To \Rep\uqsl{n-1}$. Restricted to the
fundamental part of the representation theory (see \S
\ref{sec:fundamental}), it becomes a functor $\GT : \FundRep \uqsl{n}
\To \Mat[\FundRep \uqsl{n}]$.

Next, I'll describe in gory detail the action of $\uqsl{n}$ on each
of its fundamental representations $\V{a}{n}$, using the
Gel`fand-Tsetlin decomposition. We introduce maps $p_{-1}$ and
$p_0$, the $\uqsl{n-1}$-linear projections of $\V{a}{n} \Onto
\V{a-1}{n-1}$ and $\V{a}{n} \Onto \V{a}{n-1}$. We also introduce the
inclusions $\im:\V{a-1}{n-1} \Into \V{a}{n}$ and $\iz:\V{a}{n-1}
\Into \V{a}{n}$.

\begin{prop}
We can describe the action of $\uqsl{n}$ on $\V{a}{n}$ recursively
as follows. On $\V{0}{n}$ and $\V{n}{n}$, $\uqsl{n}$ acts trivially:
$X_i^\pm$ by $0$, and $K_i$ by $1$. On the non-trivial
representations, we have
\begin{align}
\restrict{X_{n-1}^+}{\V{a}{n}} & = \im \iz p_{-1} p_0, \label{eq:top-Xp} \\
\restrict{X_{n-1}^-}{\V{a}{n}} & = \iz \im p_0 p_{-1}, \notag \\
\restrict{K_{n-1}}{\V{a}{n}} & = \im(\im p_{-1} + q \iz p_0) p_{-1} + \iz(q^{-1} \im p_{-1} +\iz p_0)p_0, \label{eq:top-K} \\
\intertext{and for $Z \in \uqsl{n-1}$}
\restrict{Z}{\V{a}{n}} & = \im \restrict{Z}{\V{a-1}{n-1}} p_{-1} + \iz \restrict{Z}{\V{a}{n-1}} p_0. \label{eq:lower-generators}
\end{align}
\end{prop}
\begin{rem}
Relative to the direct sum decomposition $\V{a}{n} \Iso
\left(\V{a-2}{n-2} \directSum \V{a-1}{n-2}\right) \directSum
\left(\V{a-1}{n-2} \directSum \V{a}{n-2}\right)$ under $\uqsl{n-2}$,
we can write these as matrices, as
\begin{align*}
X_{n-1}^+ & = \begin{pmatrix} 0 & 0 & 0 & 0 \\ 0 & 0 & 1 & 0 \\0 & 0
& 0 & 0 \\ 0 & 0 & 0 & 0 \end{pmatrix}, &
X_{n-1}^- & = \begin{pmatrix} 0 & 0 & 0 & 0 \\ 0 & 0 & 0 & 0 \\0 & 1 & 0 & 0 \\ 0 & 0 & 0 & 0 \end{pmatrix}, \\
K_{n-1}   & = \begin{pmatrix} 1 & 0 & 0 & 0 \\ 0 & q & 0 & 0 \\0 & 0
& q^{-1} & 0 \\ 0 & 0 & 0 & 1 \end{pmatrix}, & Z         & =
\begin{pmatrix} Z & 0 & 0 & 0 \\ 0 & Z & 0 & 0 \\ 0 & 0 & Z & 0 \\ 0 & 0 & 0 & Z \end{pmatrix}.
\end{align*}
\end{rem}
\begin{proof}
We need to check the relations in $\uqsl{n}$ involving $X_{n-1}^\pm$
and $K_{n-1}$. We're inductively assuming the others hold. To begin,
the $K_i$'s automatically commute, as they're diagonal with respect
to the direct sum decomposition. Further, it's clear from the
definitions that $X_{n-1}^\pm$ and $K_{n-1}$ commute with the
subalgebra $\uqsl{n-2}$. We then need to check the relation $K_i
X_j^\pm = q^{\pm(i,j)} X_j^\pm K_i$, for $\abs{i-j} \leq 1$. Now
\begin{align*}
K_{n-1} X_{n-1}^+ & = (\im(\im p_{-1} + q \iz p_0) p_{-1} + \iz(q^{-1} \im p_{-1} +\iz p_0)p_0) \im \iz p_{-1} p_0 \\
                  & = q \im \iz p_{-1} p_0, \\
\intertext{while}
X_{n-1}^+ K_{n-1} & = \im \iz p_{-1} p_0 (\im(\im p_{-1} + q \iz p_0) p_{-1} + \iz(q^{-1} \im p_{-1} +\iz p_0)p_0) \\
                  & = q^{-1} \im \iz p_{-1} p_0, \\
\intertext{so $K_{n-1} X_{n-1}^+ = q^2 X_{n-1}^+ K_{n-1}$ and similarly $K_{n-1} X_{n-1}^- = q^{-2} X_{n-1}^- K_{n-1}$. Further}
K_{n_2} X_{n-1}^+ & = (\im(\im(\im p_{-1} + q \iz p_0) p_{-1} + \iz(q^{-1} \im p_{-1} +\iz p_0)p_0)p_{-1} + \\ & \qquad + \iz(\im(\im p_{-1} + q \iz p_0) p_{-1} + \iz(q^{-1} \im p_{-1} +\iz p_0)p_0)p_0) \im \iz p_{-1} p_0 \\
                  & = \im\iz(q^{-1} \im p_{-1} +\iz p_0) p_{-1} p_0, \\
\intertext{while}
X_{n-1}^+ K_{n-2} & = \im \iz p_{-1} p_0 (\im(\im(\im p_{-1} + q \iz p_0) p_{-1} + \iz(q^{-1} \im p_{-1} +\iz p_0)p_0)p_{-1} + \\ & \qquad \qquad + \iz(\im(\im p_{-1} + q \iz p_0) p_{-1} + \iz(q^{-1} \im p_{-1} +\iz p_0)p_0)p_0) \\
                  & = \im \iz (\im p_{-1} + q \iz p_0) p_{-1} p_0,
\end{align*}
so $K_{n-2} X_{n-1}^+ = q^{-1} X_{n-1}^+ K_{n-2}$, and by similar calculations $K_{n-2} X_{n-1}^- = q X_{n-1}^- K_{n-2}$, $K_{n-1} X_{n-2}^+ = q^{-1} X_{n-2}^+ K_{n-1}$ and $K_{n-1} X_{n-2}^- = q X_{n-2}^- K_{n-1}$.

Next, we check
$X_i^+ X_j^- - X_j^- X_i^+ = \delta_{i j} \frac{K_i - K_i^{-1}}{q-q^{-1}}$, for $i=j=n-1$, and for $i=n-1,j=n-2$.
\begin{align*}
X_{n-1}^+ X_{n-1}^- & = \im \iz p_{-1} p_0 \iz \im p_0 p_{-1} \\
                    & = \im \iz p_0 p_{-1} \\
X_{n-1}^- X_{n-1}^+ & = \iz \im p_0 p_{-1} \im \iz p_{-1} p_0, \\
                    & = \iz \im p_{-1} p_0 \\
\intertext{while}
K_{n-1} - K_{n-1}^{-1} & = \im(\im p_{-1} + q \iz p_0) p_{-1} + \iz(q^{-1} \im p_{-1} +\iz p_0)p_0 - \\
                       & \qquad - \im(\im p_{-1} + q^{-1} \iz p_0) p_{-1} - \iz(q \im p_{-1} +\iz p_0)p_0 \\
                       & = (q-q^{-1}) \im \iz p_0 p_{-1} + (q^{-1}-q) \iz \im p_{-1} p_0, \\
\intertext{so $X_{n-1}^+ X_{n-1}^- - X_{n-1}^- X_{n-1}^+ = \frac{K_{n-1} - K_{n-1}^{-1}}{q-q^{-1}}$, as desired, and}
X_{n-1}^+ X_{n-2}^- & = \im \iz p_{-1} p_0 (\im \iz \im p_0 p_{-1} p_{-1} + \iz \iz \im p_0 p_{-1} p_0) \\
                    & = \im \iz p_{-1} \iz \im p_0 p_{-1} p_0 = 0 \\
X_{n-2}^- X_{n-1}^+ & = (\im \iz \im p_0 p_{-1} p_{-1} + \iz \iz \im p_0 p_{-1} p_0) \im \iz p_{-1} p_0 \\
                    & = \im \iz \im p_0 p_{-1} \iz p_{-1} p_0 = 0
\end{align*}
so $X_{n-1}^+ X_{n-2}^- - X_{n-1}^+ X_{n-2}^- = 0$.

Finally, we don't need to check the Serre relations. By a result of
\cite{MR1104219}, if you defined a quantum group without the Serre
relations, you'd see them reappear in the ideal of elements acting
by zero on all finite dimensional representations.
\end{proof}

\begin{cor}
We'll collect here some formulas for the generators acting on the
dual representations. Abusing notation somewhat, we write $\iz:
\V{a}{n-1}^* \Into \V{a}{n}^*$, $\im: \V{a-1}{n-1}^* \Into
\V{a}{n}^*$, $p_0: \V{a}{n}^* \Onto \V{a}{n-1}^*$ and $p_{-1}:
\V{a}{n}^* \Onto \V{a-1}{n-1}^*$. With these conventions,
$(\restrict{\iz}{\V{a}{n-1}})^* = \restrict{p_0}{\V{a}{n}^*}$, and
so on. Then
\begin{align*}
\restrict{X_{n-1}^+}{\V{a}{n}^*} & = -q      \iz \im p_0 p_{-1}, \\
\restrict{X_{n-1}^-}{\V{a}{n}^*} & = -q^{-1} \im \iz p_{-1} p_0, \\
\restrict{K_{n-1}}{\V{a}{n}^*}   & = \im(\im p_{-1} + q^{-1} \iz p_0) p_{-1} + \iz(q^ \im p_{-1} +\iz p_0)p_0, \\
\intertext{and for $Z \in \uqsl{n-1}$} \restrict{Z}{\V{a}{n}^*} & =
\im \restrict{Z}{\V{a-1}{n-1}^*} p_{-1} + \iz \restrict{Z}{\V{a}{n-1}^*} p_0.
\end{align*}
\end{cor}

\section{Strictifying $\mathbf{Rep}\, \uqsl{n}$}
\label{sec:strictification}%
We'll now make something of a digression into abstract nonsense.
(But it's good abstract nonsense!) It's worth understanding at this
point that while $\Rep \uqsl{n}$ isn't a strict pivotal category, it
can be `strictified'. That is, the natural isomorphism $\tau : \Id
\To **$ isn't the identity, but $\Rep \uqsl{n}$ is equivalent to
another pivotal category in which that `pivotal isomorphism' $\tau$
\emph{is} the identity. This strictification will be a full
subcategory of $\Rep \uqsl{n}$ (as a tensor category), but with a
new, modified duality functor.

While this discussion may seem a little esoteric, it actually pays
off! Our eventual goal is a diagrammatic category equivalent to
$\Rep \uqsl{n}$, with diagrams which are only defined up to planar
isotopy. Such a pivotal category will automatically be strict. The
`strictification' we perform in this section bridges part of the
inevitable gap between such a diagrammatic category, and the
conventionally defined representation category.

Of course, it would be possible to have defined $\Rep \uqsl{n}$ in
the first place in a way which made it strict as a pivotal category,
but this would have required a strange and unmotivated definition of
the dual of a morphism. Hopefully in the strictification we describe
here, you'll see the exact origin of this unfortunate dual.

First of all, let's explicitly describe the pivotal category
structure on $\Rep \uqsl{n}$, in order to justify the statement that
it is not strict\footnote{Not strict as a pivotal category, that is.
We've defined it in such a way that it's strict as a tensor
category, meaning we don't bother explicitly reassociating tensor
products.}. We need to describe the duality functor $*$, pairing
maps, and the `pivotal' natural isomorphism $\tau: \Id \To **$.

The contravariant\footnote{In fact, it's doubly contravariant; both
with respect to composition, and tensor product.} duality functor
$*$ on $\Rep \uqsl{n}$ is given by the usual duality functor for
linear maps. We need to dress up duals of representations of
$\uqsl{n}$ as representations again, which we do via the the
antipode $S$, just as we use the comultiplication to turn tensor
products of representations into representations. Thus for $V$ some representation of $\uqsl{n}$, $Z \in
\uqsl{n}$, $f \in V^*$ and $v \in V$, we say $(Zf)(v) = f(S(Z)v)$.

The pairing maps $p_V : V^* \tensor V \To \A$ are just the usual
duality pairing maps for duals of vector spaces. That these are maps
of representations follows easily from the Hopf algebra axioms. Be
careful, however, to remember that the other vector space pairing
map $V \tensor V^* \To \A$ is not generally a map of representations
of the Hopf algebra. Similarly, the copairing $c_V : \A \To V
\tensor V^*$ is equivariant, while $\A \To V^* \tensor V$ is not.

The pivotal isomorphism is where things get interesting, and we make
use of the particular structure of the quantum group $\uqsl{n}$.

\begin{defn}
We'll define the element $\tau_n = K^\rho = \prod_{j=1}^{n-1}
K_j^{j(n-j)}$ of $\uqsl{n}$.
\end{defn}
\begin{lem}
\label{lem:tau}%
Abusing notation, this gives an isomorphism of $\uqsl{n}$
representations, $\tau_n:V \IsoTo V^{**}$, providing the components
of the pivotal natural isomorphism.
\end{lem}
\begin{proof}
This relies on the formulas for the antipode acting on $\uqsl{n}$;
we need to check $$\tau_n X_i^\pm = S^2(X_i^\pm) \tau_n.$$ Since
$S^2(X_i^\pm) = q^{\pm2} X_i^\pm$, this becomes the condition
$$ \tau_n X_i^\pm \tau_n^{-1} = q^{\pm2} X_i^\pm,$$ which follows
immediately from the definition of $\tau_n$ and the commutation
relations in the quantum group.

We also need to check that $\tau$ satisfies the axioms for a pivotal natural isomorphism.
Since the $K_i$ are group-like in $\uqsl{n}$ (i.e. $\comult(K_i) = K_i \tensor K_i$), $\tau$ is a tensor natural transformation.
Further
\begin{align*}
\restrict{\tau_V^*}{V^{***}} \compose \restrict{\tau_{V^*}}{V^*} & = \restrict{(K^\rho)^*}{V^{***}} \compose \restrict{K^\rho}{V^*} \\
    & = (K^\rho) \antipode(K^\rho) \\
    & = 1,
\end{align*}
so $\tau_V^*$ and $\tau_{V^*}$ are inverses, as required.
\end{proof}

To construct the strictification of $\Rep \uqsl{n}$, written as $\SRep \uqsl{n}$, we begin by defining it as a tensor category, taking the full tensor
subcategory of $\Rep \uqsl{n}$ whose objects are generated by tensor product (but not duality) from the
trivial representation $\A$, $\V{a}{n}$ and $\V{a}{n}^*$, for each $a
= 1, \dotsc, n-1$. Before describing the pivotal structure on $\SRep
\uqsl{n}$, we might as well specify the equivalence; already there's
only one sensible choice. We can include $\SRep \uqsl{n}$ into $\Rep
\uqsl{n}$ in one direction. For the other, we simply send
$\V{a}{n}^{*(n)}$ (that is, the $n$-th iterated dual) to $\V{a}{n}$
or $\V{a}{n}^*$, depending on whether $n$ is even or odd, and modify
each morphism by pre- and post-composing with the
unique\footnote{It's unique, given the axioms for $\tau$ described
in \S \ref{sec:pivotal}, and established for our particular $\tau$
in Lemma \ref{lem:tau}.} appropriate tensor product of compositions of $\tau$
and $\tau^{-1}$. Clearly the composition $\Rep \uqsl{n} \To \SRep
\uqsl{n} \To \Rep \uqsl{n}$ is not the identity functor, but equally
easily it's naturally isomorphic to the identity, again via
appropriate tensor products of compositions of $\tau$. Here, we're
doing nothing more than identifying objects in a category related by
a certain family of isomorphisms, and pointing out that the result
is equivalent to what we started with!

The interesting aspect comes in the pivotal structure on $\SRep
\uqsl{n}$. Of course, we define the duality cofunctor on objects of
$\SRep \uqsl{n}$ so it exchanges $\V{a}{n}$ and $\V{a}{n}^*$. It
turns out (inevitably, according to the strictification result of
\cite{MR1686423}, but explicitly here) that if we define pairing
maps on $\SRep \uqsl{n}$ by taking the image of the pairing maps for
$\Rep \uqsl{n}$, we make $\SRep \uqsl{n}$ a strict pivotal category
(and of course, we thus make the equivalence of categories from the
previous paragraph an equivalence of pivotal categories).
Explicitly, then, the pairing morphisms on $\SRep \uqsl{n}$ are
$\Spairing_a : \V{a}{n}^* \tensor \V{a}{n} \To \A = \pairing_a$ and
$\Spairing_{a^*} : \V{a}{n} \tensor \V{a}{n}^* \To \A =
\pairing_{a^*} \compose (\tau_a \tensor \Id_{a^*})$. We've yet to
define the duality cofunctor at the level of morphisms in $\SRep
\uqsl{n}$; there's in fact a unique definition forced on us by Lemma
\ref{lem:dual-formula}. Thus
$$f^* = (\Id_{a^*} \tensor \pairing_{b^*}) \compose (\Id_{a^*}
\tensor (\tau_b \compose f \compose \tau_a^{-1}) \tensor \Id_{b^*}))
\compose (\pairing_a^* \tensor \Id_{b^*}).$$

The processes of strictifying, and of limiting our attention to the fundamental part
of the representation theory, are independent; the discussion above
applies exactly to $\FundRep \uqsl{n}$. In that case, we obtain the
full subcategory $\SFundRep \uqsl{n}$.

\section{Generators for $\FundRep \uqsl{n}$}
\label{sec:generators-for-fundrep}%
In this section we'll define certain morphisms in $\FundRep
\uqsl{n}$, which we'll term `elementary' morphisms. Later, in
Proposition \ref{prop:elementary-morphisms-generate} we'll show that
they generate all of $\FundRep \uqsl{n}$ as a tensor category.

We already have the duality pairing and copairing maps, discussed
above, $\pairing_V : V^* \tensor V \To \A$, and $\copairing_V : \A
\To V \tensor V^*$, for each fundamental representation $V =
\V{a}{n}$, and for arbitrary iterated duals $V = \V{a}{n}^{*(k)}$.
We also have the non-identity isomorphisms $\tau$ identifying
objects with their double duals. Beyond those, we'll introduce some
more interesting maps:
\begin{itemize}
\item a map which identifies a fundamental representations with the dual of another fundamental representation, $d_{a,n}:\V{a}{n} \IsoTo \V{n-a}{n}^*$,
  along with the inverses of these maps,
\item a `triple invariant', living in the tensor product of three fundamental representations $\vout{n}{a,b,c}: \A \To \V{a}{n} \tensor \V{b}{n} \tensor \V{c}{n}$ with $a+b+c=n$, and
\item a `triple coinvariant', $\vin{n}{a,b,c} : \V{a}{n} \tensor \V{b}{n} \tensor \V{c}{n} \To \A$, again with $a+b+c=n$.
\end{itemize}

\begin{defn}
\label{defn:d}%
The map $d_{a,n}:\V{a}{n} \To \V{n-a}{n}^*$ is specified recursively
by
\begin{align*}
d_{0,n} & : 1 \mapsto 1^*, \\
d_{n,n} & : 1 \mapsto 1^*, \\
d_{a,n} & = \iz d_{a-1,n-1} p_{-1} + (-q)^{-a} \im d_{a,n-1} p_0.
\end{align*}
\end{defn}

\begin{lem}
\label{lem:d-map} The map $d_{a,n}$ is a map of $\uqsl{n}$
representations.
\end{lem}
\begin{proof}
Since $d_{a,n}$ is defined in terms of $\uqsl{n-1}$ equivariant
maps, we need only check $d_{a,n}$ commutes with $X_{n-1}^\pm$ and
$K_{n-1}$. We'll do one calculation explicitly:
\begin{align*}
d_{a,n} X_{n-1}^+ & = \left(\iz d_{a-1,n-1} p_{-1} + (-q)^{-a} \im d_{a,n-1} p_0\right) \im \iz p_{-1} p_0 \\
    & = \iz d_{a-1,n-1} \iz p_{-1} p_0 \\
    & = \iz \left(\iz d_{a-2,n-2} p_{-1} + (-q)^{1-a} \im d_{a-1,n-2} p_0\right) \iz p_{-1} p_0 \\
    & = (-q)^{1-a} \iz \im d_{a-1,n-2} p_{-1} p_0, \\
\intertext{while}
X_{n-1}^+ d_{a,n} & = -q \iz \im p_0 p_{-1} \left(\iz d_{a-1,n-1} p_{-1} + (-q)^{-a} \im d_{a,n-1} p_0\right) \\
    & = (-q)^{1-a} \iz \im p_0 d_{a,n-1} p_0) \\
    & = (-q)^{1-a} \iz \im p_0 \left( \iz d_{a-1,n-2} p_{-1} + (-q)^{-a} \im d_{a,n-2} p_0 \right) p_0) \\
    & = (-q)^{1-a} \iz \im d_{a-1,n-2} p_{-1} p_0).
\end{align*}
The other two cases (commuting with $X_{n-1}^-$ and $K_{n-1}$) are
pretty much the same.
\end{proof}

\begin{prop}
\label{prop:duals-of-d}%
The duals of the maps $d_{a,n}$ satisfy:
\begin{equation}
\label{eq:duals-of-d} d_{a,n}^* \tau_n = (-1)^{(n+1)a} d_{n-a,n}.
\end{equation}
\end{prop}

We need two lemmas before proving this.
\begin{lem}
\label{lem:product-of-Ks}
$$\prod_{j=1}^{n-1} \restrict{K_j^j}{\V{a}{n}} = q^{n-a} \im p_{-1} + q^{-a} \iz p_0.$$
\end{lem}
\begin{proof}
Arguing inductively, and using the formula for
$\restrict{K_{n-1}}{\V{a}{n}}$ from Equation \eqref{eq:top-K}, we
find
\begin{align*}
 \prod_{j=1}^{n-1} \restrict{K_j^j}{\V{a}{n}}
    & = \left(\im\left(\prod_{j=1}^{n-2} \restrict{K_j^j}{\V{a-1}{n-1}}\right) p_{-1} + \iz \left(\prod_{j=1}^{n-2} \restrict{K_j^j}{\V{a}{n-1}}\right) p_0\right) \restrict{K_{n-1}^{n-1}}{\V{a}{n}} \\
    & = \left(\im(q^{n-a} \im p_{-1} + q^{-a+1} \iz p_0)p_{-1} + \iz(q^{n-a-1} \im p_{-1} + q^{-a} \iz p_0)p_0\right) \times \\
    & \qquad \times \left(\im(\im p_{-1} + q^{n-1} \iz p_0) p_{-1} + \iz(q^{-n+1}\im p_{-1} +\iz p_0)p_0\right) \\
    & = q^{n-a} \im\im p_{-1} p_{-1} + q^{n-a} \im \iz p_0 p_{-1} + q^{-a} \iz\im p_{-1} p_0 + q^{-a} \iz\iz p_0 p_0 \\
    & = q^{n-a} \im p_{-1} + q^{-a} \iz p_0. \qedhere
\end{align*}
\end{proof}

\begin{lem}
\label{lem:tau-branching}
$$\restrict{\tau_n}{\V{a}{n}} = q^{n-a} \im \restrict{\tau_{n-1}}{\V{a-1}{n-1}} p_{-1}+ q^{-a} \iz \restrict{\tau_{n-1}}{\V{a}{n-1}} p_0$$
\end{lem}
\begin{proof}
Making use of Equation \eqref{eq:lower-generators} and Lemma
\ref{lem:product-of-Ks}:
\begin{align*}
\restrict{\tau_n}{\V{a}{n}} & = \restrict{\tau_{n-1}}{\V{a}{n}} \times \prod_{j=1}^{n-1} \restrict{K_j^j}{\V{a}{n}} \\
    & = \left( \im \restrict{\tau_{n-1}}{\V{a-1}{n-1}} p_{-1} + \iz \restrict{\tau_{n-1}}{\V{a}{n-1}} p_0 \right) \times \left( q^{n-a} \im p_{-1} + q^{-a}\iz p_0 \right) \\
    & = q^{n-a} \im \restrict{\tau_{n-1}}{\V{a-1}{n-1}} p_{-1}+ q^{-a} \iz \restrict{\tau_{n-1}}{\V{a}{n-1}} p_0. \qedhere
\end{align*}
\end{proof}

\begin{proof}[Proof of Proposition \ref{prop:duals-of-d}]
The proposition  certainly holds for $a=0$ or $a=n$, where all three
of the maps $d_{a,n}^*, \tau_n$ and $d_{n-a,n}$ are just the
identity, and the sign is $+1$. Otherwise, we proceed inductively.
First, write
$$
 d_{a,n}^* = \im d_{a-1,n-1}^* p_0 + (-q)^{-a} \iz d_{a,n-1}^* p_{-1}.
$$
Using Lemma \ref{lem:tau-branching} we then have
\begin{align*}
d_{a,n}^* \tau_n
    & = \left(\im d_{a-1,n-1}^* p_0 + (-q)^{-a} \iz d_{a,n-1}^* p_{-1}\right) \left( q^a \im \tau_{n-1} p_{-1}+ q^{a-n} \iz \tau_{n-1} p_0 \right) \\
    & = (-q)^{-a} q^a \iz d_{a,n-1}^* \tau_{n-1} p_{-1} + q^{a-n} \im d_{a-1,n-1}^* \tau_{n-1} p_0 \\
    & = (-1)^a \iz (-1)^{na} d_{n-1-a,n-1} p_{-1} + q^{a-n} \im (-1)^{n(a-1)} d_{n-a,n-1} p_0 \\
    & = (-1)^{(n+1)a} \iz d_{n-1-a,n-1} p_{-1} + q^{a-n} (-1)^{n(a-1)} \im d_{n-a,n-1} p_0, \\
\intertext{while} (-1)^{(n+1)a} d_{n-a,n}
    & = (-1)^{(n+1)a} (\iz d_{n-a-1,n-1} p_{-1} + (-q)^{a-n} \im d_{n-a,n-1} p_0) \\
    & = (-1)^{(n+1)a} \iz d_{n-a-1,n-1} p_{-1} + q^{a-n} (-1)^{n(a-1)} \im d_{n-a,n-1} p_0). \qedhere
\end{align*}
\end{proof}

Finally, we need to define the triple invariants $\vout{n}{a,b,c}$
and $\vin{n}{a,b,c}$.

\begin{defn}
\label{defn:trivalent-vertices}
When $a+b+c=m$, we define $\vout{n}{a,b,c} \in \V{a}{n} \tensor
\V{b}{n} \tensor \V{c}{n}$ and $\vin{n}{a,b,c} \in \V{a}{n}^*
\tensor \V{b}{n}^* \tensor \V{c}{n}^*$ by the formulas
\begin{align}
 \vout{0}{0,0,0} & = 1 \tensor 1 \tensor 1, \notag \\
 \vin{0}{0,0,0} & = 1^* \tensor 1^* \tensor 1^*, \notag \\
 \vout{n}{a,b,c} & = (-1)^c q^{b+c}                 (\im \tensor \iz \tensor \iz)(\vout{n-1}{a-1,b,c}) + {} \label{eq:vout} \\
                     & \phantom{=} (-1)^a q^{c\phantom{+b}} (\iz \tensor \im \tensor \iz)(\vout{n-1}{a,b-1,c}) + {} \notag \\
                     & \phantom{=} (-1)^b \phantom{q^{b+c}} (\iz \tensor \iz \tensor \im)(\vout{n-1}{a,b,c-1}), \notag \\
\intertext{and}
 \vin{n}{a,b,c} & = (-1)^c \phantom{q^{-a-b}}   (\vin{n-1}{a-1,b,c})(p_{-1} \tensor p_0 \tensor p_0) + {} \label{eq:vin} \\
                    & \phantom{=} (-1)^a q^{-a\phantom{-b}}   (\vin{n-1}{a,b-1,c})(p_0 \tensor p_{-1} \tensor p_0) + {} \notag \\
                    & \phantom{=}(-1)^b q^{-a-b}             (\vin{n-1}{a,b,c-1})(p_0 \tensor p_0 \tensor p_{-1}). \notag
\end{align}
\end{defn}
\begin{lem}
The maps $\vout{n}{a,b,c}$ and $\vin{n}{a,b,c}$ are maps of
$\uqsl{n}$ representations.
\end{lem}
\begin{proof}
As in Lemma \ref{lem:d-map}, we just need to check that these maps
commute with $X_{n-1}^\pm$ and $K_{n-1}$. We'll do the explicit
calculation for $X_{n-1}^+$; here we need to check that
$$\restrict{X_{n-1}^+}{\V{a}{b} \tensor \V{b}{n} \tensor \V{c}{b}}
\vout{n}{a,b,c} = 0$$ and $$\vin{n}{a,b,c}
\restrict{X_{n-1}^+}{\V{a}{b} \tensor \V{b}{n} \tensor \V{c}{b}} =
0.$$ First, we need to know how $X_{n-1}^+$ acts on the tensor
product of three representations, via
\begin{equation*}
\comult^{(2)}(X_{n-1}^+) = X_{n-1}^+ \tensor K_{n-1} \tensor K_{n-1}
+ 1 \tensor X_{n-1}^+ \tensor K_{n-1} + 1 \tensor 1 \tensor
X_{n-1}^+.
\end{equation*}
Next, let's use two steps of the inductive definition of
$\vout{n}{a,b,c}$ to write
\begin{align*}
 \vout{n}{a,b,c} & =  q^{2b+2c}           (\im\im \tensor \iz\iz \tensor \iz\iz)(\vout{n-2}{a-2,b,c}) + {} \\
                   & \quad q^{2c\phantom{+2b}} (\iz\iz \tensor \im\im \tensor \iz\iz)(\vout{n-2}{a,b-2,c}) + {} \\
                   & \quad \phantom{q^{2b+2c}} (\iz\iz \tensor \iz\iz \tensor \im\im)(\vout{n-2}{a,b,c-2}) + {} \\
                   & \quad (-1)^{a+b} q^c      (- \iz\iz \tensor \im\iz \tensor \iz\im + q^{-1} \iz\iz \tensor \iz\im \tensor \im\iz) (\vout{n-2}{a,b-1,c-1}) + {} \\
                   & \quad (-1)^{b+c} q^{b+c}  (  \im\iz \tensor \iz\iz \tensor \iz\im - q^{-1} \iz\im \tensor \iz\iz \tensor \im\iz) (\vout{n-2}{a-1,b,c-1}) + {} \\
                   & \quad (-1)^{a+c} q^{b+2c} (- \im\iz \tensor \iz\im \tensor \iz\iz + q^{-1} \iz\im \tensor \im\iz \tensor \iz\iz) (\vout{n-2}{a-1,b-1,c}).
\end{align*}
We'll now show $X_{n-1}^+$ kills each term (meaning each line, as
displayed above) separately. In each case, it follows immediately
from the formulas for $X_{n-1}^+$ and $K_{n-1}$ in Equations
\eqref{eq:top-Xp} and \eqref{eq:top-K}. In the first three terms, we
use $X_{n-1}^+ \im \im = X_{n_1}^+ \iz \iz = 0$. In the fourth term,
we see
\begin{align*}
 \comult^{(2)}(X_{n-1}^+)(- \iz\iz & \tensor \im\iz \tensor \iz\im + q^{-1} \iz\iz \tensor \iz\im \tensor \im\iz) \\
    & = -\iz\iz \tensor \im\iz \tensor \im\iz + q^{-1} q \iz\iz \tensor \im\iz \tensor \im\iz  \\
    & = 0.
\end{align*}
Here only one of the three terms of $\comult^{(2)}(X_{n-1}^+)$ acts
nontrivially on each of the two terms. The fifth and sixth terms are
exactly analogous.

For $\vin{n}{a,b,c}$ we use the same trick; write it in terms of
$\vin{n-2}{a-2,b,c}$, $\vin{n-2}{a,b-2,c}$, $\vin{n-2}{a,b,c-2}$,
$\vin{n-2}{a,b-1,c-1}$, $\vin{n-2}{a-1,b,c-1}$ and
$\vin{n-2}{a-1,b-1,c}$, and show that each of these terms multiplied
by $X_{n-1}^+$ gives zero separately.
\end{proof}
\begin{rem}
Going through this proof carefully, you'll see that it would still
work with some variation allowed in constants the definitions of
$\vout{n}{a,b,c}$ and of $\vin{n}{a,b,c}$ in Equations
\eqref{eq:vout} and \eqref{eq:vin}. However, given the normalisation
for $d_{a,n}$ that we've chosen in Definition \ref{defn:d}, these
normalisation constants are pinned down by Lemmas
\ref{lem:dGT-rotating} and \ref{lem:dGT-degeneration} below.
\end{rem}

Later, we'll discuss relationships between these elementary
morphisms, but for now we want to justify our interest in them.

\begin{prop}
\label{prop:elementary-morphisms-generate}%
The elementary morphisms generate, via tensor product,
composition and linear combination, all the morphisms in $\FundRep \uqsl{n}$.
\end{prop}
\begin{rem}
Certainly, they can't generate all of $\Rep \uqsl{n}$, simply
because the sources and targets of elementary morphisms are all in
$\FundRep \uqsl{n}$!
\end{rem}

\begin{rem}
A fairly abstract proof of this fact for $n=3$ has been given by
Kuperberg \cite{q-alg/9712003}.\footnote{His other results
additionally give a direct proof for $n=3$, although quite different
from the one here.} Briefly, he specialises to $q=1$, and considers
the subcategory of $\FundRep \usl{3} \cong \FundRep SU(3)$ generated
by the elementary morphisms. He extends this by formally adding
kernels and cokernels of morphisms (that is, by taking the Karoubi
envelope). This extension
being (equivalent to) all of $\Rep SU(3)$ is enough to obtain the
result. To see this, he notes that the extension is the sort of
category to which an appropriate Tannaka-Krein theorem applies,
allowing him to say that it is the representation category of some
compact Lie group. Some arguments about the symmetries of the
`triple invariant' morphisms allow him to conclude that this group
must in fact be $SU(3)$. While I presume this proof can be extended
to cover all $n$, I prefer to give a more direct proof, based on the
quantum version of Frobenius-Schur duality.
\end{rem}
\begin{proof}
It's a sort of bootstrap argument. First, we notice that the action
of the braid group $\Braids_m$ on tensor powers of the standard
representation $\V{1}{n}$ can be written in terms of elementary
morphisms. Next, we recall that the braid group action generates all
the endomorphisms of $\V{1}{n}^{\tensor m}$, and finally, we show
how to map arbitrary tensor products of fundamental representations
into a tensor power of the standard representation, using only
elementary morphisms.

These ideas are encapsulated in the following four lemmas.

\begin{lem}
The action of the braid group $\Braids_m$ on $\V{1}{n}^{\tensor m}$,
given by $R$-matrices, can be written in terms of elementary
morphisms.
\end{lem}
\begin{proof}
We only need to prove the result for $\V{1}{n}^{\tensor 2}$. There,
we can be particularly lazy, taking advantage of the fact that $\dim
\Hom{\uqsl{n}}{\V{1}{n}^{\tensor 2}}{\V{1}{n}^{\tensor 2}} = 2$.
(There are many paths to seeing this, the path of least effort perhaps being that
$\Hom{\uqsl{n}}{\V{1}{n}^{\tensor 2}}{\V{1}{n}^{\tensor 2}}$ is isomorphic to
$\Hom{\uqsl{n}}{\V{1}{n} \tensor \V{1}{n}^*}{\V{1}{n}^* \tensor \V{1}{n}}$, and that both
the source and target representations there decompose into the direct sum of the adjoint representation and the trivial representation.)
The map $$p = (\vin{n}{1,1,n-2} \tensor \Id_{\V{1}{n}} \tensor
\Id_{\V{1}{n}}) \compose (\Id_{\V{1}{n}} \tensor \Id_{\V{1}{n}}
\tensor \vout{n}{n-2,1,1})$$ 
is not a multiple of the identity, 
so every endomorphism of $\V{1}{n}^{\tensor 2}$ is a linear
combination of compositions of elementary morphisms, in particular
the braiding.
\end{proof}
\begin{rem}
In fact, the map associated to a positive crossing is $q^{n-1}
\Id_{\V{1}{n}^{\tensor 2}} - q^n p$. 
It's easy to prove that this, along with the negative crossing (obtained by replacing $q$ with $q^{-1}$), is the only linear combination of $\Id$
and $p$ which satisfies the braid relation. One could presumably
also check that this agrees with the explicit formulas given in
\cite[\S 8.3 and \S 10.1]{MR1300632}\footnote{Although be careful
there --- the formula for the universal $R$-matrix in \S 8.3.C is
incorrect, although it shouldn't matter for such a small representation. The order of the product is backwards
\cite{chari-pressley-correction}.}.
\end{rem}

\begin{lem}
The image of $\Braids_m$ linearly spans
$\Hom{\uqsl{n}}{\V{1}{n}^{\tensor m}}{\V{1}{n}^{\tensor m}}$.
\end{lem}
\begin{proof}
This is the quantum version of Schur-Weyl duality. See \cite[\S 10.2B]{MR1300632} and \cite{MR841713}.
\end{proof}

\begin{lem}
For any tensor product of fundamental representations $\Tensor_i V_{\alpha_i}$, there's some natural number $m$ and a pair of morphisms constructed out of elementary ones
$\iota:\Tensor_i V_{\alpha_i} \To V_1^{\tensor m}$ and $\pi: V_1^{\tensor m} \To \Tensor_i V_{\alpha_i}$ such that $\pi \compose \iota = \Id_{\Tensor_i V_{\alpha_i}}$.
\end{lem}
\begin{proof}
We just need to do this for a single fundamental representation, then tensor together those morphisms. For a single fundamental representation, first
define
$\iota' : V_a \To V_1 \tensor V_{a-1} = (\pairing_{V_{n-a}} \tensor \Id_{V_1} \tensor \Id_{V_{a-1}}) \compose (d_a \tensor \vout{n}{n-a,1,a-1})$ and
$\pi' : V_1 \tensor V_{a-1} \To V_a = x (\vin{n}{1,a-1,n-a} \tensor d_a) \compose (\Id_{V_1} \tensor \Id_{V_{a-1}} \tensor \copairing_{V_{n-a}})$.
The composition $\pi' \compose \iota'$ is an endomorphism of the irreducible $V_a$, which I claim is nonzero, and so for some choice of the coefficient $x$
is the identity. Now build the maps $\iota : V_a \To V_1^{\tensor a}$ and $\pi : \V_1^{\tensor a} \To V_a$ as iterated compositions of these maps.
\end{proof}

\begin{lem}
In fact, given two such tensor products of fundamental representations $\Tensor_i V_{\alpha_i}$ and $\Tensor_j V_{\beta_j}$, such that there
is some nonzero $\uqsl{n}$ map between them, it's possible to choose $\iota_\alpha, \pi_\alpha, \iota_\beta$ and $\pi_\beta$ as in the previous lemma,
with the same value of $m$.
\end{lem}
\begin{proof}
(This argument is due to Ben Webster. Thanks!)
First just pick $\iota_\alpha, \pi_\alpha, \iota'_\beta$ and $\pi'_\beta$ as in the previous lemma, with possibly different values of $m$, say $m_\alpha$
and $m_\beta$. Now, if there's some $\uqsl{n}$ map between $\Tensor_i V_{\alpha_i}$ and $\Tensor_j V_{\beta_j}$, then there's also a $\uqsl{n}$
map between $V_1^{\tensor m_\alpha}$ and $V_1^{\tensor m_\beta}$. This can happen if and only if there's some $\csl{n}$ map, and in fact some $SL(n)$ map, between the corresponding
classical representations, ${\Complex^n}^{\tensor m_\alpha}$ and ${\Complex^n}^{\tensor m_\beta}$. This $SL(n)$ intertwiner guarantees that
the center of $SL(n)$ acts in the same way on each representation; the element $\sqrt[n]{1} I$ acts by $1^{m_\alpha/n}$ and $1^{m_\beta/n}$
on the two representations, so $m_\alpha$ and $m_\beta$ must be congruent mod $n$. Assuming, without loss of generality, that $m_\beta \leq m_\alpha$ we now
define $\iota_\beta = \iota'_\beta \tensor \iota^{\tensor (m_\alpha - m_\beta)/n}$,
and $\pi_\beta = \pi'_\beta \tensor \pi^{\tensor (m_\alpha - m_\beta)/n}$, where $\iota : \A \To V_1^{\tensor n}$ and $\pi : V_1^{\tensor n} \To \A$
are diagrammatic morphisms satisfying $\pi \compose \iota = \Id_\A$. (These certainly exist!)
\end{proof}

Putting all this together, we take an arbitrary morphism $\phi: \Tensor_i V_{\alpha_i}
\To \Tensor_j V_{\beta_j}$ in $\FundRep \uqsl{n}$, and pick $\iota_\alpha, \pi_\alpha, \iota_\beta$ and $\pi_\beta$ as in
the last lemma. Then
$$\phi = \pi_\beta \compose \iota_\beta \compose \phi \compose \pi_\alpha \compose \iota_\alpha.$$
The middle of this composition, $\iota_\beta \compose \phi \compose \pi_\alpha$, is an endomorphism of $V_1^{\tensor m}$, so can be written
as a linear combination of braids, and hence as a linear combination of morphisms built from elementary morphisms.
\end{proof}

\section{The representation functor}
\label{sec:representation-functor}%
Finally, everything is in place for a definition of the
representation functor $\Rep : \SymCat_n \To \SFundRep \uqsl{n}$.
(Notice that we're defining the functor with target the strictification
$\SFundRep \uqsl{n}$, as discussed in \S \ref{sec:strictification}.
To obtain the functor to the more usual representation theory
$\FundRep \uqsl{n}$, you need to compose with the equivalences of
categories described in that section; of course, this is just the
inclusion of a full subcategory!) First define $\Rep'$ on the free
tensor category $\FreeCat_n$, by
\begin{align*}
\Rep'(a \in \{1, \dotsc, n-1\})  & = \V{a}{n}, &
\Rep'(a^*)                       & = \V{a}{n}^*,
\end{align*}
on objects, and
\begin{align*}
\Rep'\left( \mathfig{0.06}{pairings/pl_a}  \right) & = p_a, &
\Rep'\left( \mathfig{0.06}{pairings/pr_a}  \right) & = p_{a^*} \compose (\tau_n \tensor \Id), \\
\Rep'\left( \mathfig{0.06}{pairings/cpl_a} \right) & = c_a, &
\Rep'\left( \mathfig{0.06}{pairings/cpr_a} \right) & = (Id \tensor \tau_n^{-1}) \compose c_{a^*},
\end{align*}
\begin{align*}
\Rep'\left( \mathfig{0.065}{duals/dtr} \right) & = d_{a,n}, &
\Rep'\left( \mathfig{0.065}{duals/dal} \right) & = d_{a,n}^{-1}, \\
\Rep'\left( \mathfig{0.065}{duals/dtl} \right) & = (-1)^{(n+1)a} d_{a,n}, &
\Rep'\left( \mathfig{0.065}{duals/dar} \right) & = (-1)^{(n+1)a} d_{a,n}^{-1},
\end{align*}
and
\begin{align*}
\Rep'\left( \mathfig{0.1}{vertices/v_outup_0} \right) & = \vout{n}{a,b,c} &
\Rep'\left( \mathfig{0.1}{vertices/v_inup_0} \right)  & = \vin{n}{a,b,c}
\end{align*}
on morphisms.

\begin{prop}
\label{prop:rep-descends}
The functor $\Rep'$ descends to a functor defined on the quotient $\SymCat_n \To \SFundRep \uqsl{n}$.
\end{prop}
This is proved in the next chapter; we could do it now, but the proof will read more nicely once we have better diagrams available.

\chapter{The diagrammatic Gel`fand-Tsetlin functor}
\label{sec:dGT}

It's now time to define the diagrammatic Gel`fand-Tsetlin functor,
$$\dGT : \SymCat_n \To \Mat[\SymCat_{n-1}].$$ The first incarnation of the diagrammatic
functor will be a functor defined (in \S \ref{sec:dGT-on-generators}) on the free version of
the diagrammatic category, $\dGT' : \FreeCat_n \To \Mat[\SymCat_{n-1}]$. We'll need to show that it descends to
the quotient $\SymCat_n$ (in \S \ref{sec:dGT-descends}). This
definition is made so that the perimeter of the diagram
\begin{equation}
\label{eq:the-diagram}
 \xymatrix@C+10mm{
    \FreeCat_n \ar[dr]^{\text{Draw}} \ar@/^/[drr]^{\Rep'} \ar@/_/[ddr]_{\dGT'} & & \\
    & \SymCat_n \ar[d]^{\dGT} \ar[r]^-{\Rep} & \SFundRep \uqsl{n} \ar[d]^{\GT} \\
    & \Mat[\SymCat_{n-1}] \ar[r]^-{\Mat[\Rep]} & \Mat[\SFundRep \uqsl{n-1}]
 }
\end{equation}
commutes.

Happily, as an easy consequence of this, we'll see that the
representation functor $\Rep : \FreeCat_n \To \SFundRep \uqsl{n}$ also
descends to the quotient $\SymCat_n$.

At the end of the chapter in \S \ref{sec:path-model} we describe how
to compute the diagrammatic Gel`fand-Tsetlin functor, and perform a few small calculations which we'll need later.

\section{Definition on generators}
\label{sec:dGT-on-generators}%
The target category for the diagrammatic Gel`fand-Tsetlin functor is
the matrix category over $\SymCat_{n-1}$. We thus need to send each
object of $\FreeCat_n$ to a direct sum of objects of
$\SymCat_{n-1}$, and for each generating morphism of $\FreeCat_n$,
we need to pick an appropriate matrix of morphisms in $\SymCat_{n-1}$. On
objects, we use the obvious
\begin{align}
\label{eq:dGT-on-objects}%
\dGT'(a)   & = (a-1) \directSum a \\
\dGT'(a^*) & = (a-1)^* \directSum a^* \notag
\end{align}
(omitting the nonsensical direct summand in the case that $a = 0,
0^*, n$ or $n^*$), extending to tensor products by distributing over
direct sum. For morphisms we'll take advantage of the lack of
`multiplicities' in Equation \ref{eq:dGT-on-objects} to save on some
notation. An arbitrary morphism in $\Mat[\SymCat_{n-1}]$ has rows
and columns indexed by tensor products of the fundamental (and
trivial) objects and their duals, $0, 1, \dotsc, n-1, 0^*, 1^*,
\dotsc, (n-1)^*$. Notice, however, that morphisms in the image of
$\dGT'$ have distinct labels on each row (and also on each column).
Moreover, a (diagrammatic) morphism in $\SymCat_{n-1}$ explicitly
encodes its own source and target (reading across the top and bottom
boundary points). We'll thus abuse notation, and write a sum of
matrix entries, instead of the actual matrix, safe in the knowledge
that we can unambiguously reconstruct the matrix, working out which term should sit in each matrix entry. For example, if
$f:a \To b$ is a morphism in $\FreeCat_n$, then $\dGT'(f)$ is a
matrix $(a-1) \directSum a \To (b-1) \directSum b$, which we ought
to write as $\psmallmatrix{f_{11} & f_{12} \\ f_{21} & f_{22}}$, but
will simply write as $f_{11} + f_{12} + f_{21} + f_{22}$. This
notational abuse is simply for the sake of brevity when writing down
matrices with many zero entries. When composing matrices written
this way, you simply distribute the composition over summation,
ignoring any non-composable terms.

That said, we now define $\dGT'$ on trivalent vertices by
\begin{align*}
 \dGT'\left(\mathfig{0.1}{vertices/v_outup_0}\right) & =
     (-1)^c q^{b+c} \mathfig{0.1}{vertices/v_outup_1} + (-1)^a q^c \mathfig{0.09}{vertices/v_outup_2} + (-1)^b \mathfig{0.1}{vertices/v_outup_3}, \\
 \dGT'\left(\mathfig{0.1}{vertices/v_inup_0}\right) & =
     (-1)^c \mathfig{0.1}{vertices/v_inup_1} + (-1)^a q^{-a} \mathfig{0.095}{vertices/v_inup_2} + (-1)^b q^{-a-b} \mathfig{0.1}{vertices/v_inup_3} \\
\end{align*}
(the right hand sides of these equations are secretly $8 \times 1$ and $1 \times 8$ matrices, respectively), on the cups and caps by
\begin{align}
 \dGT'\left(\mathfig{0.066}{pairings/pl_a}\right) & = \mathfig{0.066}{pairings/pl_a} + \mathfig{0.099}{pairings/pl_am} &
 \dGT'\left(\mathfig{0.066}{pairings/pr_a}\right) & = q^{n-a} \mathfig{0.066}{pairings/pr_a} + q^{-a} \mathfig{0.099}{pairings/pr_am} \notag \\
 \dGT'\left(\mathfig{0.066}{pairings/cpl_a}\right) & = \mathfig{0.066}{pairings/cpl_a} + \mathfig{0.099}{pairings/cpl_am} &
 \dGT'\left(\mathfig{0.066}{pairings/cpr_a}\right) & = q^{a-n} \mathfig{0.066}{pairings/cpr_a} + q^a \mathfig{0.099}{pairings/cpr_am},
 \label{eq:dGT-on-pairings}%
\end{align}
and on the tags by
\begin{align}
 \dGT'\left(\mathfig{0.065}{duals/dtr}\right) & = \mathfig{0.076}{duals/dtrb} + (-1)^a q^{-a} \mathfig{0.088}{duals/dtrt} \notag \\
 \dGT'\left(\mathfig{0.065}{duals/dtl}\right) & = (-1)^{n+a} \mathfig{0.076}{duals/dtlb} + q^{-a} \mathfig{0.088}{duals/dtlt} \notag \displaybreak[1] \\
 \dGT'\left(\mathfig{0.065}{duals/dar}\right) & = q^a \mathfig{0.088}{duals/darb} + (-1)^{n+a} \mathfig{0.076}{duals/dart} \notag \\
 \dGT'\left(\mathfig{0.065}{duals/dal}\right) & = (-1)^a q^a \mathfig{0.088}{duals/dalb} + \mathfig{0.076}{duals/dalt}. \label{eq:dGT-on-tags}%
\end{align}
The functor $\dGT'$ then extends to all of $\FreeCat_n$ as a tensor functor.

\begin{prop}
\label{prop:outer-square-commutes}%
The outer square of Equation \eqref{eq:the-diagram} commutes.
\end{prop}
\begin{proof}
Straight from the definitions; compare the definition above with the
definition of $\Rep'$ from \S \ref{sec:representation-functor}, and
use the inductive Definitions \ref{defn:d} and
\ref{defn:trivalent-vertices} of the generating morphisms over in
the representation theory.
\end{proof}

\section{Descent to the quotient}
\label{sec:dGT-descends}

\begin{prop}
The functor $\dGT'$ descends from $\FreeCat_n$ to the quotient
$\SymCat_{n}$, where we call it simply $\dGT$.
\end{prop}
\begin{proof}
Read the following Lemmas \ref{lem:dGT-straightening},
\ref{lem:dGT-tag-near-maxima}, \ref{lem:dGT-rotating},
\ref{lem:dGT-flipping-tags} and \ref{lem:dGT-degeneration}.
\end{proof}

\begin{lem}
\label{lem:dGT-straightening}%
We first check the relation in Equation
\eqref{eq:straighten-strand}.
\begin{align*}
\dGT'\left(\mathfig{0.1}{isotopy/strand_left}\right) & = \dGT'\left(\mathfig{0.036}{isotopy/strand}\right)  = \dGT'\left(\mathfig{0.1}{isotopy/strand_right}\right). \\
\dGT'\left(\mathfig{0.1}{isotopy/strand_down_left}\right) & = \dGT'\left(\mathfig{0.036}{isotopy/strand_down}\right) = \dGT'\left(\mathfig{0.1}{isotopy/strand_down_right}\right).
\end{align*}
\end{lem}
\begin{proof}
This is direct from the definitions in Equation \eqref{eq:dGT-on-pairings}, and using the fact that we can straighten strands in the target category $\SymCat_{n-1}$.
\end{proof}

\begin{lem}
\label{lem:dGT-tag-near-maxima}%
Next, the relations in Equation \eqref{eq:tag-near-maxima}.
\begin{align*}
 \dGT'\left(\mathfig{0.1}{pairings/cp_tag_left}\right) & = \dGT'\left(\mathfig{0.1}{pairings/cp_tag_right}\right). \\
 \dGT'\left(\mathfig{0.1}{pairings/p_tag_left}\right)  & = \dGT'\left(\mathfig{0.1}{pairings/p_tag_right}\right).
\end{align*}
\end{lem}
\begin{proof}
\begin{align*}
 \dGT'\left(\mathfig{0.1}{pairings/cp_tag_left}\right)
    & = q^{n-a} \mathfig{0.1}{pairings/cp_tag_left_a} + (-1)^a \mathfig{0.1}{pairings/cp_tag_left_b} \\
    & = q^{n-a} \mathfig{0.1}{pairings/cp_tag_right_a} + (-1)^a \mathfig{0.1}{pairings/cp_tag_right_b} \\
    & = \dGT'\left(\mathfig{0.1}{pairings/cp_tag_right}\right) \\
\intertext{making use of the definitions in Equations
\eqref{eq:dGT-on-pairings} and \eqref{eq:dGT-on-tags}, and the
identity appearing in Equation \eqref{eq:tag-near-maxima} for
$\SymCat_{n-1}$, and}
 \dGT'\left(\mathfig{0.1}{pairings/p_tag_left}\right)
    & = \mathfig{0.1}{pairings/p_tag_left_a} + (-1)^a q^{-a} \mathfig{0.1}{pairings/p_tag_left_b} \\
    & = \mathfig{0.1}{pairings/p_tag_right_a} + (-1)^a q^{-a} \mathfig{0.1}{pairings/p_tag_right_b} \\
    & = \dGT'\left(\mathfig{0.1}{pairings/p_tag_right}\right). \qedhere
\end{align*}
\end{proof}

\begin{lem}
\label{lem:dGT-rotating}%
We check the `rotation relation' from Equation
\eqref{eq:rotate-vertex} is in the kernel of $\dGT'$, which of
course implies the `$2\pi$ rotation relation from Equations
\eqref{eq:rotate-vertex-2pi} is also in the kernel.
\begin{align}
\label{eq:dGT-v_in}%
\dGT'&\left(\mathfig{0.16}{isotopy/v_in_rotated}\right) \\
    & = (-1)^c \mathfig{0.16}{dGT/v_in_rotated_a} + (-1)^a q^{-a} \mathfig{0.16}{dGT/v_in_rotated_b} + (-1)^b q^{-a-b} \mathfig{0.16}{dGT/v_in_rotated_c} \notag \\
    & = \dGT'\left(\mathfig{0.1}{isotopy/v_in}\right). \notag \\
\dGT'&\left(\mathfig{0.16}{isotopy/v_out_rotated}\right) \\
    & = (-1)^c q^{b+c} \mathfig{0.16}{dGT/v_out_rotated_a} + (-1)^a q^c \mathfig{0.16}{dGT/v_out_rotated_b} + (-1)^b \mathfig{0.16}{dGT/v_out_rotated_c} \notag \\
    & = \dGT'\left(\mathfig{0.1}{isotopy/v_out}\right). \notag
\end{align}
\end{lem}
\begin{proof}
The statement of the lemma essentially contains the proof; we'll spell out Equation \eqref{eq:dGT-v_in} in gory detail:
\begin{align*}
\dGT'&\left(\mathfig{0.16}{isotopy/v_in_rotated}\right) \\
    & = \dGT'\left(\mathfig{0.06}{pairings/pl_c}\right) \compose
        \left(\Id_{c^*} \tensor \dGT'\left(\mathfig{0.1}{vertices/v_inup_cab_0}\right) \tensor \Id_c\right) \compose \left(\dGT'\left(\mathfig{0.06}{pairings/cpr_c}\right) \tensor \Id_{a \tensor b \tensor c}\right) \displaybreak[1] \\
    & = \left(\mathfig{0.06}{pairings/pl_c} + \mathfig{0.09}{pairings/pl_cm}\right) \compose \\
    & \qquad \compose \left(\Id_{c^*} \tensor \left((-1)^b \mathfig{0.1}{vertices/v_inup_cab_1} + (-1)^c q^{-c} \mathfig{0.1}{vertices/v_inup_cab_2} + \right. \right. \\
    & \qquad \qquad \qquad \qquad \qquad \qquad \qquad + \left. \left. (-1)^a q^{-c-a} \mathfig{0.1}{vertices/v_inup_cab_3}\right) \tensor \Id_c\right) \compose \\
    & \qquad \compose \left(\left(q^{c-n} \mathfig{0.06}{pairings/cpr_c} + q^c \mathfig{0.09}{pairings/cpr_cm}\right) \tensor \Id_{a \tensor b \tensor c}\right) \displaybreak[1] \\
    & = (-1)^c \mathfig{0.16}{dGT/v_in_rotated_a} + (-1)^a q^{-a} \mathfig{0.16}{dGT/v_in_rotated_b} + (-1)^b q^{-a-b} \mathfig{0.16}{dGT/v_in_rotated_c} \\
    & = \dGT'\left(\mathfig{0.1}{isotopy/v_in}\right).
\end{align*}
Going from the third last to the second last line, we simply throw out all non-composable cross terms.
\end{proof}

The next two lemmas are similarly direct.
\begin{lem}
\label{lem:dGT-flipping-tags}%
\begin{align*}
 \dGT'\left(\mathfig{0.065}{duals/dtr}\right) & = \mathfig{0.065}{duals/dtrb} + (-1)^a q^{-a}  \mathfig{0.075}{duals/dtrt} = (-1)^{(n+1)a} \dGT'\left(\mathfig{0.065}{duals/dtl}\right). \\
 \dGT'\left(\mathfig{0.065}{duals/dar}\right) & = q^a \mathfig{0.075}{duals/darb} + (-1)^{n+a} \mathfig{0.065}{duals/dart} = (-1)^{(n+1)a} \dGT'\left(\mathfig{0.065}{duals/dal}\right).
\end{align*}
\end{lem}

\begin{lem}
\label{lem:dGT-degeneration}%
\begin{align*}
 \dGT'\left(\mathfig{0.115}{vertices/v_outup_degenerate}\right) & = \dGT'\left(\mathfig{0.115}{pairings/cp_tag_right}\right), \\
 \dGT'\left(\mathfig{0.115}{vertices/v_inup_degenerate}\right)  & = \dGT'\left(\mathfig{0.115}{pairings/p_tag_left}\right).
\end{align*}
\end{lem}

We can now give the
\begin{proof}[Proof of Proposition \ref{prop:rep-descends}.]
Somewhat surprisingly, the fact that $\dGT'$ descends from
$\FreeCat_n$ to $\SymCat_n$ implies that $\Rep':\FreeCat_n \To
\SFundRep \uqsl{n}$ also descends to $\SymCat_n$, simply because the
outer square of Equation \eqref{eq:the-diagram} commutes, as pointed out in Proposition \ref{prop:outer-square-commutes}.
\end{proof}

\section{Calculations on small webs}
\label{sec:dGT-calculations}%
We'll begin with some calculations of
$\dGT$ on the flow vertices introduced in \S
\ref{sec:flow-vertices}.
\begin{align*}
 \dGT\left(\mathfig{0.12}{vertices/sideways_vertex1_0}\right) & = (-1)^a \left(
    \mathfig{0.12}{vertices/sideways_vertex1_0} +
    (-1)^{n+b} \mathfig{0.155}{vertices/sideways_vertex1_1} +
    q^{-a} \mathfig{0.135}{vertices/sideways_vertex1_2}
 \right) \displaybreak[1] \\
 \dGT\left(\mathfig{0.12}{vertices/sideways_vertex2_0}\right) & = (-1)^a \left(
    \mathfig{0.12}{vertices/sideways_vertex2_0} +
    (-1)^{n+b} q^b \mathfig{0.155}{vertices/sideways_vertex2_1} +
    \mathfig{0.145}{vertices/sideways_vertex2_2}
 \right)
\end{align*}

\begin{align*}
 \dGT\left(\mathfig{0.12}{vertices/sideways_vertex3_0}\right) & = (-1)^a \left(
    \mathfig{0.12}{vertices/sideways_vertex3_0} +
    (-1)^{n+b} q^b \mathfig{0.12}{vertices/sideways_vertex3_1} +
    \mathfig{0.12}{vertices/sideways_vertex3_2}
 \right) \displaybreak[1] \\
 \dGT\left(\mathfig{0.12}{vertices/sideways_vertex4_0}\right) & = (-1)^a \left(
    q^a \mathfig{0.12}{vertices/sideways_vertex4_0} +
    \mathfig{0.12}{vertices/sideways_vertex4_1} +
    (-1)^{n+b} q^{a-n} \mathfig{0.12}{vertices/sideways_vertex4_2}
 \right)
\end{align*}

We'll also need the following formulas for two vertex webs:
\begin{multline}
\dGT\left(\mathfig{0.25}{PQ/P_section_0}\right) =
    \mathfig{0.25}{PQ/P_section_0} + {} \\
  + q^{-a} \mathfig{0.26}{PQ/P_section_1} + \mathfig{0.26}{PQ/P_section_2} + {} \\
  + q^{-a} \mathfig{0.275}{PQ/P_section_3} + (-1)^{b+c} q^b \mathfig{0.275}{PQ/P_section_4} \label{eq:dGT-P-section}
\end{multline}
\begin{multline}
\dGT\left(\mathfig{0.18}{PQ/Q_section_0}\right) =
    q^b \mathfig{0.18}{PQ/Q_section_0} + {} \\
  + \mathfig{0.18}{PQ/Q_section_1} + \mathfig{0.18}{PQ/Q_section_2} + {} \\
  + q^b \mathfig{0.18}{PQ/Q_section_3} + (-1)^{a+c} q^{b+c-n} \mathfig{0.18}{PQ/Q_section_4}  \label{eq:dGT-Q-section}
\end{multline}

\section{A path model, and polygons.}
\label{sec:path-model}%
Call a diagram in $\SymCat_n$ in which all vertices are flow
vertices and there are no tags (except those implicitly hidden
inside the vertices) a \emph{flow diagram}. In particular, the flow
vertices themselves, and the $\mathcal{P}$- and
$\mathcal{Q}$-polygons described in \S \ref{sec:polygons} are flow
diagrams. A \emph{reduction path} on a flow diagram is some disjoint
union of cycles and paths beginning and ending at the boundary of
the diagram, always matching the orientations on the edges of the
diagram. We then notice that if $v$ is a flow vertex, the three terms of $\dGT(v)$ correspond to the three reduction paths on the vertex.
Specifically, for each reduction path, in the corresponding term of
$\dGT(v)$ the label on each traversed edge has been reduced by one.
Each term additionally has a coefficient $\pm q^k$ for some integer
$k$.

It's easy to see that this condition must also hold for larger flow
diagrams. For a reduction path $\pi$ on a diagram $D \in \SymCat_n$,
write $\pi(D)$ for the diagram in $\SymCat_{n-1}$ obtained from $D$
by reducing the label on each edge traversed by $\pi$. Further,
write $t(\pi)$ for the subset of the tensor factors in the target of
$D$ which $\pi$ reaches, and $s(\pi)$ for the subset of the tensor
factors in the source.

\begin{prop}
The result of applying $\dGT$ to a flow diagram $D$ is a sum over
the reduction paths on $D$:
\begin{equation}
\dGT(D) = \sum_{\substack{\text{reduction}\\\text{paths $\pi$}}} \pm q^\bullet \pi(D).
\end{equation}
Further, two terms in this sum lie in same matrix entry (remember,
$\dGT$ is a functor $\SymCat_n \To \Mat[\SymCat_{n-1}]$) iff the
boundaries of the two reduction paths agree.
\end{prop}
\begin{proof}
Certainly this is true for the generators which can appear in a flow
diagram; caps, cups, and flow vertices. The condition is also
preserved under tensor product and composition.

To see this, suppose we have two flow diagrams $D_1$ and $D_2$ in
$\SymCat_n$, with families of reduction paths $P_1$ and $P_2$, so
that $\dGT(D_i) = \sum_{\pi \in P_i} \pm q^\bullet \pi(D_i)$. Now
the reduction paths $P$ for $D_1 \tensor D_2$ are exactly those of
the form $\pi_1 \cup \pi_2$ for some $\pi_1 \in P_1$, $\pi_2 \in
P_2$ (remember tensor product of diagrams is side-by-side disjoint union). Thus
\begin{align*}
\dGT(D_1 \tensor D_2) & = \dGT(D_1) \tensor \dGT(D_2) \\
    & = \left(\sum_{\pi_1 \in P_1} (-1)^\bullet q^\bullet \pi_1(D_1)\right) \tensor \left(\sum_{\pi_2 \in P_2} (-1)^\bullet q^\bullet \pi_2(D_2)\right) \\
    & = \sum_{\substack{\pi_1 \in P_1 \\ \pi_2 \in P_2}} (-1)^\bullet q^\bullet (\pi_1 \cup \pi_2)(D_1 \tensor D_2) \displaybreak[1] \\
    & = \sum_{\pi \in P} (-1)^\bullet q^\bullet \pi(D_1 \tensor D_2).
\end{align*}
The argument for compositions is a tiny bit more complicated.
Supposing $D_1$ and $D_2$ are composable, the reduction paths $P$
for $D_1 \compose D_2$ are those of the form $\pi_1 \cup \pi_2$, for
some $\pi_1 \in P_1$, $\pi_2 \in P_2$ such that $s(\pi_1) =
t(\pi_2)$. Thus
\begin{align*}
\dGT(D_1 \compose D_2) & = \dGT(D_1) \compose \dGT(D_2) \\
    & = \left(\sum_{\pi_1 \in P_1} \pm q^\bullet \pi_1(D_1)\right) \compose \left(\sum_{\pi_2 \in P_2} \pm q^\bullet \pi_2(D_2)\right) \\
    & = \sum_{\substack{\pi_1 \in P_1 \\ \pi_2 \in P_2}} \pm q^\bullet \delta_{s(\pi_1), t(\pi_2)} (\pi_1 \cup \pi_2)(D_1 \compose D_2) \displaybreak[1] \\
    & = \sum_{\pi \in P} \pm q^\bullet \pi(D_1 \tensor D_2). \qedhere
\end{align*}
\end{proof}

We'll now calculate $\dGT$ on the $\mathcal{P}$- and
$\mathcal{Q}$-polygons. On a $\mathcal{P}$-polygon, there are two
reduction paths which do not traverse any boundary points; the empty
reduction path, and the path around the perimeter of the polygon.
We'll call these $\pi_{\mathcal{P};\eset}$ and $\pi_{\mathcal{P};\circ}$. Otherwise, there is
a unique reduction path for each (cyclic) subsequence of boundary
points which alternately includes incoming and outgoing edges.
We'll write $\pi_{\mathcal{P};s}$ for this, where $s$ is the subset
of $\{1,\dotsc,2k\}$ corresponding to the points on boundary
traversed by the path components. The reduction path has one
component for each such pair of incoming and outgoing edges, and
traverses the perimeter of the polygon counterclockwise between
them. On a $\mathcal{Q}$-polygon, there are two reductions paths
which traverse every boundary point; in one, each component enters
at an incoming edge, and traverse counterclockwise, departing at the
nearest outgoing edge, while in the other each component traverses
clockwise. We'll call these $\pi_{\mathcal{Q};\circlearrowleft}$ and $\pi_{\mathcal{Q};\circlearrowright}$.
Otherwise, there is a unique reduction path for each disjoint
collection of adjacent pairs of boundary edges; the path connects
each pair, traversing a single edge of the polygon between them, and
again, we'll write $\pi_{\mathcal{Q};s}$, where $s$ is the subset of
$\{1,\dotsc,2k\}$ corresponding to the subset of the boundary
consisting of the union of all the pairs. All these statements about
reduction paths should be immediately obvious, looking at the
definition of the $\mathcal{P}$- and $\mathcal{Q}$-polygons in
Equations \eqref{eq:P-defn} and \eqref{eq:Q-defn}.

Now,
\begin{align*}
\pi_{\mathcal{P};\eset}(\P{n}{a,b}{l}) & = \P{n-1}{a,b}{l} \\
\intertext{and}
\pi_{\mathcal{P};\circ}(\P{n}{a,b}{l}) & = \P{n-1}{a+\vect{1},b+\vect{1}}{l} = \P{n}{a,b}{l-1}. \\
\intertext{Further}
\pi_{\mathcal{P};s}(\P{n}{a,b}{l}) & = \P{n-1}{a+a',b+b'}{l},
\end{align*}
where $a',b' \in \{0,1\}^k$ is determined
from $s$ by $2i - 1 \in s$ iff $a'_i = 1$ and $b'_i = 0$, and $2i
\in s$ iff $a'_{i+1} = 1$ and $b'_i = 0$. The condition that $s$
contains alternately incoming and outgoing edges translates to $b'_i
\leq a'_i,a'_{i+1}$; we'll call such pairs $a', b'$
$\mathcal{P}$-admissible. Note that $\pi_{\mathcal{P};\eset}$ and $\pi_{\mathcal{P};\circ}$ correspond to two exceptional pairs, $a'=b'=\vect{0}$ and $a'=b'=\vect{1}$ respectively.

On the $\mathcal{Q}$-polygons, we have $\pi_{\mathcal{Q};\circlearrowleft}(\Q{n}{a,b}{l})
= \Q{n}{a,b+\vect{-1}}{l} = \Q{n}{a+\vect{1},b}{l+1}$, $\pi_{\mathcal{Q};\circlearrowright}(\Q{n}{a,b}{l}) =
\Q{n}{a+\vect{1},b}{l}$ and $\pi_{\mathcal{Q};s}(\Q{n}{a,b}{l}) =
\Q{n-1}{a+a',b+b'}{l}$, where $a' \in \{0,1,\}^k$ and $b' \in
\{-1,0\}^k$ are determined by $2i - 1 \in s$ iff $a_i = 1$ or
$b_i = -1$, and $2i \in s$ iff $a_{i+1} = 1$ or $b_i = -1$. The
condition that the points in $s$ come in disjoint adjacent pairs
translates to the condition $a'_i b'_i = 0 = a'_{i+1} b'_i$; such
pairs $a', b'$ are called $\mathcal{Q}$-admissible. Note that $\pi_{\mathcal{Q};\circlearrowleft}$ and $\pi_{\mathcal{Q};\circlearrowright}$ correspond to two exceptional pairs, $a'=\vect{0}, b'=\vect{-1}$ and $a'=\vect{1},b'=\vect{0}$ respectively.
Also, note that a pair $a',b'$ is both $\mathcal{P}$- and $\mathcal{Q}$-admissible exactly if $b'=\vect{0}$; any $a' \in \{0,1\}^k$ is allowed.

Thus we have
\begin{align*}
\dGT[\P{n}{a,b}{l}] & = \sum_{a',b' \text{ $\mathcal{P}$-admissible}} \dGT_{a',b'}\left(\P{n}{a,b}{l}\right) \\
\intertext{and}
\dGT[\Q{n}{a,b}{l}] & = \sum_{a',b' \text{ $\mathcal{Q}$-admissible}} \dGT_{a',b'}\left(\Q{n}{a,b}{l}\right)
\end{align*}
where
\begin{align}
\label{eq:dGTab-on-P}%
\dGT_{a',b'}\left(\P{n}{a,b}{l}\right) & = (-1)^{\inner{b'}{a+\rotl[a]}} q^{l(\ssum{b'}-\ssum{a'}+1)} q^{\inner{\rotl[a']}{b}-\inner{b'}{a}-a_1-n a_1'} \P{n-1}{a+a',b+b'}{l}, \\
\intertext{and}
\label{eq:dGTab-on-Q}%
\dGT_{a',b'}\left(\Q{n}{a,b}{l}\right) & = (-1)^{\inner{b'}{a+\rotl[a]}} q^{l(\ssum{a'}-\ssum{b'}-\frac{\abs{\bdy}}{2}+1)} q^{\ssum{b} - n\ssum{b'} +\inner{b'}{\rotl[a]}-\inner{a'}{b}-a_1-n a_1'} \Q{n-1}{a+a',b+b'}{l}.
\end{align}
The coefficients here are products of the coefficients appearing in the formulas for $\dGT$ on the two vertex webs in Equations \eqref{eq:dGT-P-section} and \eqref{eq:dGT-Q-section}.
We'll often also use the notation $\dGTe$ for the terms corresponding to reduction paths not traversing any boundary edges, and $\dGTa$ for the terms
corresponding to reduction paths traversing all boundary edges. Thus
\begin{align}
\label{eq:dGTe-P}%
\dGTe[\P{n}{a,b}{l}] & = \dGT_{\vect{0},\vect{0}}\left(\P{n}{a,b}{l}\right) + \dGT_{\vect{1},\vect{1}}\left(\P{n}{a,b}{l}\right) \notag \\
                     & = q^{l-a_1} \P{n-1}{a,b}{l} + q^{l-a_1+\ssum{b}-\ssum{a}-n} \P{n-1}{a,b}{l-1}, \displaybreak[1] \\
\label{eq:dGTe-Q}%
\dGTe[\Q{n}{a,b}{l}] & = \dGT_{\vect{0},\vect{0}}\left(\Q{n}{a,b}{l}\right) = q^{l-a_1} q^{\ssum{b} - \frac{l\abs{\bdy}}{2}} \Q{n-1}{a,b}{l}, \displaybreak[1] \\
\label{eq:dGTa-P}%
\dGTa[\P{n}{a,b}{l}] & = \dGT_{\vect{1},\vect{0}}\left(\P{n}{a,b}{l}\right) = q^{l-a_1-n} q^{\ssum{b} - \frac{l\abs{\bdy}}{2}} \P{n-1}{a+\vect{1},b}{l}, \\
\intertext{and}
\label{eq:dGTa-Q}%
\dGTa[\Q{n}{a,b}{l}] & = \dGT_{\vect{0},\vect{-1}}\left(\Q{n}{a,b}{l}\right) + \dGT_{\vect{1},\vect{0}}\left(\Q{n}{a,b}{l}\right) \notag \\
                     & = q^{l-a_1+\ssum{b} - \ssum{a} -\frac{n\abs{\bdy}}{2}} \Q{n-1}{a+\vect{1},b}{l+1} + q^{l-a_1-n} \Q{n-1}{a+\vect{1},b}{l}.
\end{align}

\chapter{Describing the kernel}
\label{sec:kernel}

Finally, in this chapter we'll use the methods developed to this
point to try to pin down the kernel of the representation functor.
In particular, Theorems \ref{thm:I=H}, \ref{thm:square-switch} and
\ref{thm:kekule}, stated in the three subsequent sections,
describe particular classes of relations amongst diagrams. Theorems \ref{thm:square-switch} and
\ref{thm:kekule} additionally rule out any other similar relations.

\section{The $I=H$ relations}
\begin{thm}
\label{thm:I=H}%
For each $a,b,c \in \{0, \dotsc, n\}$ with $a+b+c \leq n$, there's an identity 
\begin{align}
\label{eq:I=H}
\mathfig{0.1}{IH/I_abc} & = (-1)^{(n+1)a} \mathfig{0.1}{IH/H_abc}, \\
\intertext{and another}
\label{eq:I=H-outgoing}
\mathfig{0.1}{IH/I_out_abc} & = (-1)^{(n+1)a} \mathfig{0.1}{IH/H_out_abc}.
\end{align}
\end{thm}
The proof appears in \S \ref{sec:proofs-of-relations}.
\begin{rem}
It's somewhat unfortunate that there's a sign here in the first place, and that it depends unsymmetrically on the parameters. (See also \S \ref{sec:temperley-lieb} and \S \ref{sec:kim}, for a discussion of sign differences between my setup and previous work.) Replacing the vertices in Equation \eqref{eq:I=H} with standard ones,
we can equivalently write this identity in the following four ways (obtaining each by flipping tags, per Equation \eqref{eq:dual-of-a-tag}):
\begin{align}
\label{eq:I=H-equivalent}%
 \mathfig{0.1}{IH/I_abcd_l} & = (-1)^{(n+1)a} \mathfig{0.1}{IH/H_abcd_d} &
 \mathfig{0.1}{IH/I_abcd_r} & = (-1)^{(n+1)b} \mathfig{0.1}{IH/H_abcd_d} \displaybreak[1] \\
 \mathfig{0.1}{IH/I_abcd_r} & = (-1)^{(n+1)c} \mathfig{0.1}{IH/H_abcd_u} &
 \mathfig{0.1}{IH/I_abcd_l} & = (-1)^{(n+1)d} \mathfig{0.1}{IH/H_abcd_u}, \notag
\end{align}
where here $a+b+c+d=n$. These signs are special cases of the $6-j$ symbols for $\uqsl{n}$.
\end{rem}

\section{The square-switch relations}
\label{sec:square-switch}%
We can now state the `square switch' relations, which essentially say
that every $\mathcal{P}$-type square can be written as a linear
combination of $\mathcal{Q}$-type squares, as vice versa. (Refer back to \S \ref{sec:polygons} for the definitions of $\mathcal{P}$- and $\mathcal{Q}$-polygons.)

\begin{thm}
\label{thm:square-switch}%
The subspace ($\SS$ for `square switch')
\begin{align*}
 \SS^n_{a,b} =
    \begin{cases}
        \Span{\A}{\P{n}{a,b}{l} - \sum_{m=\max{a}}^{\min{b}} \qBinomial{n+\ssum{a}-\ssum{b}}{m+l-\ssum{b}} \Q{n}{a,b}{m}}_{l=\max{b}}^{\min{a}+n} & \text{if $n+\ssum{a} - \ssum{b} \geq 0$} \\
        \Span{\A}{\Q{n}{a,b}{l} - \sum_{m=\max{b}}^{n+\min{a}} \qBinomial{\ssum{b}-n-\ssum{a}}{m+l-\ssum{a}-n} \P{n}{a,b}{m}}_{l=\max{a}}^{\min{b}} & \text{if $n+\ssum{a} - \ssum{b} \leq 0$}
    \end{cases}
\end{align*}
of $\AP^n_{a,b} + \AQ^n_{a,b} \subset
\Hom{\SymCat^n}{\eset}{\Label(a,b)}$ lies in the kernel of the
representation functor.
\end{thm}
The proof appears in \S \ref{sec:proofs-of-relations}. Notice that
we're stating this result for all tuples $a$ and $b$, not just those
for which $\SS$ is non-trivial. This is a minor point, except that
it makes the base cases of our inductive proof easy.

\begin{rem}
In the case $n+\ssum{a}-\ssum{b}=0$, this becomes simply
$$SS^n_{a,b} = \Span{\A}{\P{n}{a,b}{l} - \Q{n}{a,b}{\ssum{b}-l}}_{l=\max{b}}^{\min{a}+n}.$$
\end{rem}

\begin{rem}
For $n+\ssum{a}-\ssum{b} \leq 0$, at the maximal or minimal values of $l$, we simply get
$\P{n}{a,b}{\max{b}} - \Q{n}{a,b}{\min{b}}$ (since the only terms
for which the $q$-binomial is nonzero are $m \geq \min{b}$) and
$\P{n}{a,b}{\min{a}+n} - \Q{n}{a,b}{\max{a}}$ (again, the
$q$-binomials vanish unless $m \leq \max{a}$). We already knew these
were identically zero in $\SymCat_n$, by the discussion at the end
of \S \ref{sec:polygons} on relations between $\mathcal{P}$- and
$\mathcal{Q}$-polygons. A similar statement applies when $n+\ssum{a}-\ssum{b} \leq 0$.
\end{rem}

\subsubsection{Loops, and let bigons be bygones}
Over on the representation theory side, any loops or bigon must become a multiple
of the identity, by Schur's lemma. We can see this as a special case
of the $\SS$ relations. (Alternatively, we'd be able to derive the
same formulas from Theorem \ref{thm:kekule} in the next
section.)

Consider the boundary flow labels $a=(0,0), b=(k,0)$, so $$\Label(a,b) = ((0,-),(0,+),(k,-),(k,+)).$$ Then
$\AP^n_{a,b} = \Span{\A}{\P{n}{a,b}{l}}_{l=k}^n$, while $\AQ^n_{a,b}
= \Span{\A}{\Q{n}{a,b}{0}}$. By \ref{thm:square-switch}, we have
$$\SS^n_{a,b} = \Span{\A}{\P{n}{a,b}{l} - \qBinomial{n-k}{l-k}
\Q{n}{a,b}{0}}.$$ Diagrammatically, this says that in the quotient by $\ker \Rep$,
\begin{align*}
 \mathfig{0.2}{PQ/bigonP} & = \qBinomial{n-k}{l-k} \mathfig{0.2}{PQ/bigonQ} \\
\intertext{or, removing trivial edges, and cancelling tags}
 \mathfig{0.115}{webs/bigon_nkl} & = \mathfig{0.1}{webs/bigon_kl} = \qBinomial{n-k}{l-k} \mathfig{0.024}{webs/strand_k}.
\end{align*}
The special case $k = 0$ evaluates loops:
\begin{equation*}
\mathfig{0.07}{webs/circle_l} = \qBinomial{n}{l}.
\end{equation*}

\section{The \Kekule relations}
\label{sec:kekule}%
Friedrich August \Kekule is generally thought to have been the first
to suggest the cyclic structure of the benzene molecule
\cite{wiki:Kekule}. A simplistic description of the benzene molecule
is as the quantum superposition of two mesomers, each with
alternating double bonds
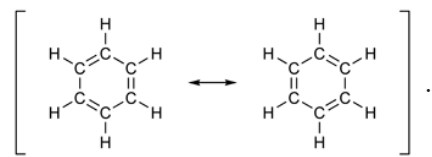
\ifpdf
    $$\mathfig{0.36}{benzene-mesomers}.$$
\else
    $$\mathfig{0.36}{benzene-mesomers_small}.$$
\fi The similarity to the relation between the two hexagonal
diagrams in $\SymCat_4/\ker$,
\begin{equation}
\label{eq:kekule}%
 \mathfig{0.1}{webs/sl_4/kekule1} + \mathfig{0.1}{webs/sl_4/kekule3} = \mathfig{0.1}{webs/sl_4/kekule2} + \mathfig{0.1}{webs/sl_4/kekule0}
\end{equation}
discovered by Kim \cite{math.QA/0310143}, prompted Kuperberg to suggest the name `\Kekule
relation', even though in my classification all of the relations for
$\uqsl{2}$ and $\uqsl{3}$ are also of this type!

In the following, we'll use the convenient notation $\sumhat{x} =
\ssum{x} - \max{x}$, $\sumtah{x} = \ssum{x} - \min{x}$. (With the
convention that $\sumhat{\eset} = \sumtah{\eset} = 0$.)

\begin{defn}
We'll say a set of flow labels $(a,b)$ is $n$-\emph{hexagonal} if
there's a sequence of six elements of $\Label(a,b)$, each
non-trivial for $n$, i.e. between $1$ and $n-1$ inclusive, which are
alternately incoming and outgoing edges.
\end{defn}
\begin{rem}
It's not enough to simply say there are $6$ non-trivial edges. For
example, the boundary labels
$$((5,+),(1,-),(0,+),(1,-),(0,+),(1,-),(0,+),(1,-),(0,+),(1,-),(0,+),(1,-))$$
seem to allow hexagonal webs, but modulo the $I=H$ relations, the
$\mathcal{P}$-polygons are all actually bigons in disguise.
\end{rem}

\begin{thm}
\label{thm:kekule}%
The subspaces $\APR^n_{a,b} \subset \AP^n_{a,b} \subset
\Hom{\SymCat^n}{\eset}{\Label(a,b)}$ and $\AQR^n_{a,b} \subset
\AQ^n_{a,b} \subset \Hom{\SymCat^n}{\eset}{\Label(a,b)}$ defined by
\begin{equation}
 \APR^n_{a,b} = \Span{\A}{ \sum_{k=-\sumhat{b}}^{-\sumtah{a} + 1} (-1)^{j+k} \qBinomial{j+k-\max{b}}{j-\ssum{b}} \qBinomial{\min{a} + n -j -k}{\ssum{a} + n -1 -j} \P{n}{a,b}{j+k}}_{j=\ssum{b}}^{\ssum{a} +n-1}
\end{equation}
\begin{multline}
 \AQR^n_{a,b} = \operatorname{span}_{\A}\left\{ \sum_{k=-\sumtah{a}}^{-\sumhat{b} + \frac{n\abs{\bdy}}{2} + 1} (-1)^{j+k} \qBinomial{j+k-\max{a}}{j-\ssum{a}} \times \right. \\
    \left. \rule{0mm}{10mm} \times \qBinomial{\min{b}-j-k}{\ssum{b}-n(\frac{\abs{\bdy}}{2}-1)-1-j} \Q{n}{a,b}{j+k} \right\}_{j=\ssum{a}}^{\ssum{b}-n(\frac{\abs{\bdy}}{2}-1)-1}
\end{multline}
are in the kernel of the representation functor. Moreover, when the
representation functor is restricted to $\AP^n_{a,b}$, its kernel is
exactly $\APR^n_{a,b}$, and similarly for $\AQ^n_{a,b}$ and
$\AQR^n_{a,b}$. Further, in the case that $(a,b)$ is $n$-hexagonal,
when the representation functor is restricted to $\AP^n_{a,b} +
\AQ^n_{a,b}$, its kernel is exactly $\APR^n_{a,b} + \AQR^n_{a,b}$.
\end{thm}%
The proof appears in \S \ref{sec:proofs-of-relations}.

\begin{rem}
It's worth being a little careful here; remember from Lemma
\ref{lem:P=Q} that sometimes $\AP^n_{a,b} \cap \AQ^n_{a,b} \neq 0$;
when $a$ or $b$ is constant, the intersection is $1$ dimensional,
and when both are constant, but not equal, it's $2$ dimensional.
Besides this, however, you can understand the second half of the
above Theorem as saying that there are no relations between
$\mathcal{P}$-type polygons, and $\mathcal{Q}$-type polygons, when
those polygons are at least hexagons. In particular, this
contradicts the conjecture of \cite{math.GT/0506403}, which expected
to find relations looking like the `square-switch' relations
described in Theorem \ref{thm:square-switch} amongst larger
polygons.
\end{rem}

The elements of $\APR^n_{a,b}$ described above have `breadth'
$\sumhat{b}-\sumtah{a}+2$; that is, they have that many terms. We
can generalise the definition of $n$-hexagonal to say that a pair
$(a,b)$ has $n$-circumference at least $2k$ if there's a sequence of $2k$ edges
in $\Label(a,b)$, which are nontrivial for $n$ and alternately incoming and
outgoing. Then
\begin{lem}
\label{lem:breadth-size}%
If $(a,b)$ has $n$-circumference at least $2k$, the
breadth of the relations in $\APR^n_{a,b}$ is at least $k+1$. In
particular, the only place that relations of breadth $2$ appear is
when $(a,b)$ has $n$-circumference no more than $2$.
\end{lem}
The proof appears in \S \ref{sec:proofs-of-relations}, after Lemma \ref{lem:breadth}.

\section{More about squares}
\begin{thm}
\label{thm:square-switch-complement}%
Further, when $n + \ssum{a} - \ssum{b} \geq 0$, the space
$\AQR^n_{a,b}$ of relations amongst the $\mathcal{Q}$-squares
becomes $0$, and in fact $\APR^n_{a,b} \subset \SS^n_{a,b}$.
(Conversely, when $n + \ssum{a} - \ssum{b} \leq 0$, $\APR^n_{a,b}=0$
and $\AQR^n_{a,b} \subset \SS^n_{a,b}$.) When $n+\ssum{a} - \ssum{b}
\geq 0$, defining
\begin{equation*}
 {\SS'}^n_{a,b} =
        \Span{\A}{\Q{n}{a,b}{m} - \sum_{l=n+\ssum{a}-m}^{n+\min{a}} (-1)^{m+l+n+\ssum{a}} \qBinomial{m+l-1-\ssum{b}}{m+l-n-\ssum{a}} \P{n}{a,b}{l} }_{m=\max{a}}^{\min{b}}
\end{equation*}
we find that $SS^n_{a,b} = \APR^n_{a,b} \directSum SS'^n_{a,b}$.
When $n+\ssum{a} - \ssum{b} \leq 0$, we define instead
\begin{equation*}
 {\SS'}^n_{a,b} =
        \Span{\A}{\P{n}{a,b}{m} - \sum_{l=\ssum{b}-m}^{\min{b}} (-1)^{m+l+\ssum{b}} \qBinomial{m+l-n-1-\ssum{a}}{m+l-\ssum{b}} \Q{n}{a,b}{l} }_{m=\max{b}}^{\min{a}+n}
\end{equation*}
and find that $\SS^n_{a,b} = \AQR^n_{a,b} \directSum SS'^n_{a,b}$.
\end{thm}

\begin{cor}
\label{cor:square-kernel}%
The kernel of $\Rep$ on $\AP^n_{a,b} + \AQ^n_{a,b}$, in the case of
squares, is exactly $SS^n_{a,b}$.
\end{cor}
\begin{proof}
This follows from Theorems \ref{thm:kekule} and \ref{thm:square-switch-complement}. Suppose $x \in \ker \Rep_n \cap (\AP^n_{a,b} + \AQ^n_{a,b})$,
we'll show that it's zero as an element
of $(\AP^n_{a,b} + \AQ^n_{a,b})/SS^n_{a,b}$. In fact, $(\AP^n_{a,b} + \AQ^n_{a,b})/SS^n_{a,b} = \AP^n_{a,b}/SS^n_{a,b} = \AP^n_{a,b}/\APR^n_{a,b}$, by Theorem \ref{thm:square-switch-complement},
and by Theorem \ref{thm:kekule}, $(\ker \Rep_n \cap \AP^n_{a,b})/\APR^n_{a,b} = 0$.
\end{proof}

\section{A conjecture}
I wish I could prove the following
\begin{conj}
\label{conjecture}%
The kernel of the representation functor is
generated, as a pivotal category ideal, by the elements described in
Theorems \ref{thm:I=H}, \ref{thm:square-switch} and
\ref{thm:kekule}.
\end{conj}
\begin{proof}[Evidence.]
First of all, we have Corollary \ref{cor:square-kernel}, which confirms that the
conjecture holds for diagrams with at most four boundary points.

Second, Theorem \ref{thm:kekule} tells us more than simply that
certain elements are in the kernel; it tells us exactly the kernel
of the representation functor when restricted to diagrams with a
single (hexagonal or bigger) polygonal face. Any relations not
generated by those discovered here must then be `more nonlocal', in
the sense that they involve diagrams with multiple polygons.

Finally, there's some amount of computer evidence. I have a
{\tt{Mathematica}} package, sadly unreleased at this point, which
can explicitly produce intertwiners between representations of
arbitrary quantum groups $U_q \mathfrak{g}$. A program which can
translate diagrams into instructions for taking tensor products and
compositions of elementary morphisms (i.e. an implementation of the
representation functor described in \S \ref{sec:representation-functor}) can then
automatically look for linear dependencies amongst diagrams. For
small values of $n$ (up to 8 or 9), I used such a program to look
for relations involving two adjacent hexagons, but found nothing
besides relations associated, via Theorem \ref{thm:kekule}, to one of
the individual hexagons.
\end{proof}

\section{Examples: $\uqsl{n}$, for $n=2,3,4$ and $5$.}
For each $n$, we can enumerate the finite list of examples of each of the three families of relations.
To enumerate $\APR^n$ relations, we find each set of flow labels
$(a,b)$ such that every element of $\Label(a,b)$ is between $1$ and
$n-1$ inclusive, and $\dim \APR^n_{a,b} = n+\ssum{a} -\ssum{b} > 0$.
When sets of flow labels $(a,b)$ and $(\rotl^k(a),\rotl^k(b))$
differ by a cyclic permutation, we'll only discuss the
lexicographically smaller one. There's no need to separately enumerate the $\AQR^n$ relations; they're just rotations of the $\APR^n$ relations.
To enumerate $\SS^n$ relations, we'll do this same, but weakening the inequality to $n+\ssum{a}-\ssum{b} \geq 0$. Recall that we don't need to look at the extreme
values of $l$ in the $\SS$ relations in Theorem \ref{thm:square-switch}; those relations automatically hold already in $\SymCat_n$.

For $n=2$ and $n=3$, there are no $\APR^n_{a,b}$ relations with $a,b$ of length $3$ or more; the relations of length $2$ or less follow from the $\SS^n$ relations.
We recover exactly the Temperley-Lieb loop relation for $n=2$, and Kuperberg's loop, bigon and square relations for $n=3$, from $\SS^3_{(0,0),(0,0)}, \SS^3_{(0,0),(0,1)}$ and $\SS^3_{(0,0),(1,1)}$ respectively.

For $n=4$, we have non-trivial $\APR^4_{a,b}$ spaces for
\begin{align*}
(a,b) & = \begin{cases}
            (\eset,\eset) & \text{when $k=0$} \\
            ((0),(1)),((0),(2)),\text{ or } ((0),(3)) & \text{when $k=1$} \\
            ((0,0),(1,1)),((0,0),(1,2)),\text{ or } ((0,0),(2,2)) & \text{when $k=2$} \\
          \end{cases}
\intertext{and for $k=3$}
(a,b) & = ((0,0,0),(1,1,1)).
\end{align*}
Thus we have relations involving loops:
\begin{align}
\APR^4_{\eset,\eset} & = \splitSpan{\A}{\qi{4} \P{4}{\eset,\eset}{0} - \P{4}{\eset,\eset}{1}, - \qi{3} \P{4}{\eset,\eset}{1} + \qi{2} \P{4}{\eset,\eset}{2},}{\qi{2} \P{4}{\eset,\eset}{2} - \qi{3} \P{4}{\eset,\eset}{3}, - \P{3}{\eset,\eset}{3} + \qi{4} \P{4}{\eset,\eset}{4}} \notag \\
                     & = \splitSpan{\A}{\qi{4} - \mathfig{0.048}{webs/circle_1}, - \qi{3} \mathfig{0.048}{webs/circle_1} + \qi{2} \mathfig{0.048}{webs/circle_2},}{\qi{2} \mathfig{0.048}{webs/circle_2} - \qi{3} \mathfig{0.048}{webs/circle_3}, - \mathfig{0.048}{webs/circle_3} + \qi{4}}, \label{eq:sl_4-loops} \\
\intertext{and bigons:}
\APR^4_{(0),(1)}     & = \splitSpan{\A}{-\qi{3} \P{4}{(0),(1)}{1} + \P{4}{(0),(1)}{2}, \qi{2} \P{4}{(0),(1)}{2} - \qi{2} \P{4}{(0),(1)}{3},}{- \P{4}{(0),(1)}{3} + \qi{3} \P{4}{(0),(1)}{4}} \notag \\
                     & = \Span{\A}{-\qi{3} \mathfig{0.02}{webs/strand_1} + \mathfig{0.06}{webs/bigon_1_2}, \qi{2} \mathfig{0.06}{webs/bigon_1_2} - \qi{2} \mathfig{0.06}{webs/bigon_1_3}} \label{eq:sl_4-bigon1}, \\
\APR^4_{(0),(2)}     & = \Span{\A}{\qi{2} \P{4}{(0),(2)}{2} - \P{4}{(0),(2)}{3}, - \P{4}{(0),(2)}{3} + \qi{2} \P{4}{(0),(2)}{4}} \notag \\
                     & = \Span{\A}{\qi{2} \mathfig{0.012}{webs/strand_2} - \mathfig{0.06}{webs/bigon_2_3}}, \label{eq:sl_4-bigon2} \\
\intertext{and}
\APR^4_{(0),(3)}     & = \Span{\A}{-\P{4}{(0),(3)}{3} + \P{4}{(0),(3)}{4}} = 0, \label{eq:sl_4-bigon3}
\end{align}
(there's a redundancy in each $APR^4_{(0),(k)}$ space, since $\P{4}{(0),(k)}{k} = \P{4}{(0),(k)}{4}$)
and then relations involving squares:
\begin{align}
\APR^4_{(0,0),(1,1)}     & = \splitSpan{\A}{-\qi{3} \P{4}{(0, 0),(1, 1)}{1} + \qi{2} \P{4}{(0, 0),(1, 1)}{2} - \P{4}{(0, 0),(1, 1)}{3},}{ \P{4}{(0, 0),(1, 1)}{2} - \qi{2} \P{4}{(0, 0),(1, 1)}{3} + \qi{3} \P{4}{(0, 0),(1, 1)}{4}} \notag \\
                         & = \splitSpan{\A}{-\qi{3} \mathfig{0.075}{webs/two_strands_horizontal} + \qi{2} \mathfig{0.075}{webs/sl_4/P_00_11_2} - \mathfig{0.075}{webs/sl_4/P_00_11_3},}{ \mathfig{0.075}{webs/sl_4/P_00_11_2} - \qi{2} \mathfig{0.075}{webs/sl_4/P_00_11_3} + \qi{3} \mathfig{0.075}{webs/two_strands_vertical}}, \label{eq:sl_4-squares11} \\
\APR^4_{(0,0),(1,2)}     & = \Span{\A}{\P{4}{(0, 0),(1, 2)}{2} - \P{4}{(0, 0),(1, 2)}{3} + \P{4}{(0, 0),(1, 2)}{4}} \notag \\
                         & = \Span{\A}{\mathfig{0.0375}{webs/sl_4/P_00_12_2} - \mathfig{0.075}{webs/sl_4/P_00_12_3} + \mathfig{0.075}{webs/sl_4/P_00_12_4}}, \label{eq:sl_4-squares12} \\
\intertext{and}
\APR^4_{(0,1),(2,2)}     & = \Span{\A}{\P{4}{(0, 1),(2, 2)}{2} - \P{4}{(0, 1),(2, 2)}{3} + \P{4}{(0, 1),(2, 2)}{4}} \notag \\
                         & = \Span{\A}{\mathfig{0.075}{webs/sl_4/P_01_22_2} - \mathfig{0.075}{webs/sl_4/P_01_22_3} + \mathfig{0.09}{webs/sl_4/P_01_22_4}}, \label{eq:sl_4-squares22} \\
\intertext{and finally the eponymous \Kekule relation, involving
hexagons:}
\APR^4_{(0,0,0),(1,1,1)} & = \splitSpan{\A}{-\P{4}{(0, 0, 0),(1, 1, 1)}{1} + \P{4}{(0, 0, 0),(1, 1, 1)}{2} - }{ - \P{4}{(0, 0, 0),(1, 1, 1)}{3} + \P{4}{(0, 0, 0),(1, 1, 1)}{4}} \notag \\
                         & = \Span{\A}{-\mathfig{0.11}{webs/sl_4/kekule0} + \mathfig{0.11}{webs/sl_4/kekule1} - \mathfig{0.11}{webs/sl_4/kekule2} + \mathfig{0.11}{webs/sl_4/kekule3}}.
\end{align}
Here we're telling some small lies; $\APR^4_{a,b}$ is always a
subspace of $\Hom{\SymCat^4}{\eset}{\Label(a,b)}$, but I've drawn
some of these diagrams as elements of other $\operatorname{Hom}$
spaces, via rotations.

The interesting $\SS^4$ spaces are
\begin{align*}
\SS^4_{(0, 0),(1, 1)} & = \splitSpan{\A}{\P{4}{(0, 0),(1, 1)}{2} - \Q{4}{(0, 0),(1, 1)}{0} - \qi{2} \Q{4}{(0, 0),(1, 1)}{1},}{\P{4}{(0, 0),(1, 1)}{3} - \qi{2} \Q{4}{(0, 0),(1, 1)}{0} - \Q{4}{(0, 0),(1, 1)}{1}} \\
                      & = \splitSpan{\A}{\mathfig{0.1}{webs/sl_4/P_00_11_2} - \mathfig{0.1}{webs/two_strands_vertical} - \qi{2} \mathfig{0.1}{webs/two_strands_horizontal},}{\mathfig{0.1}{webs/sl_4/P_00_11_3} - \qi{2} \mathfig{0.1}{webs/two_strands_vertical} - \mathfig{0.1}{webs/two_strands_horizontal}},  \displaybreak[1] \\
\SS^4_{(0, 0),(1, 2)} & = \Span{\A}{\P{4}{(0, 0),(1, 2)}{3} - \Q{4}{(0, 0),(1, 2)}{0} - \Q{4}{(0, 0),(1, 2)}{1}} \\
                      & = \Span{\A}{\mathfig{0.1}{webs/sl_4/P_00_12_3} - \mathfig{0.1}{webs/sl_4/Q_00_12_0} - \mathfig{0.051}{webs/sl_4/Q_00_12_1}},  \displaybreak[1] \\
\SS^4_{(0, 0),(2, 2)} & = \Span{\A}{\P{4}{(0, 0),(2, 2)}{3} - \Q{4}{(0, 0),(2, 2)}{1}} \\
                      & = \Span{\A}{\mathfig{0.1}{webs/sl_4/P_00_22_3} - \mathfig{0.1}{webs/sl_4/Q_00_22_1}}, \\
\intertext{and}
\SS^4_{(0, 1),(2, 2)} & = \Span{\A}{\P{4}{(0, 1),(2, 2)}{3} - \Q{4}{(0, 1),(2, 2)}{1} - \Q{4}{(0, 1),(2, 2)}{2}} \\
                      & = \Span{\A}{\mathfig{0.1}{webs/sl_4/P_01_22_3} - \mathfig{0.11}{webs/sl_4/Q_01_22_1} - \mathfig{0.1}{webs/sl_4/Q_01_22_2}}.
\end{align*}
Notice the redundancies here: $\SS^4_{(0, 0),(1, 1)} = \APR^4_{(0, 0),(1, 1)}$, $\SS^4_{(0, 0),(1, 2)} = \APR^4_{(0, 0),(1, 2)}$ and $\SS^4_{(0, 1),(2, 2)} = \APR^4_{(0, 1),(2, 2)}$,
but $\SS^4_{(0, 0),(2, 2)}$ is independent of any of the $\APR$ relations.

For $n=5$, we'll just look at non-trivial $\APR^5_{a,b}$ spaces where $a$ and $b$ are each of length at least $3$; we know
what the length $0$ and $1$ relations look like (loops and bigons), and the length $2$ relations all follow from $\SS$ relations. Thus
we have
\begin{align*}
(a,b) & = ((0,0,0),(1,1,1)),((0,0,0),(1,1,2)),((0,0,1),(1,2,2)),((0,1,1),(2,2,2)) \\
\intertext{or}
      & = ((0,0,0,0),(1,1,1,1)).
\end{align*}
Then
\begin{align}
\APR^5&_{(0,0,0),(1,1,1)} = \notag \\
        & = \splitSpan{\A}{-\qi{4} \P{5}{\vect{0},\vect{1}}{1} + \qi{3} \P{5}{\vect{0},\vect{1}}{2} - \qi{2} \P{5}{\vect{0},\vect{1}}{3} + \P{5}{\vect{0},\vect{1}}{4},}{ \P{5}{\vect{0},\vect{1}}{2} - \qi{2} \P{5}{\vect{0},\vect{1}}{3} + \qi{3} \P{5}{\vect{0},\vect{1}}{4} - \qi{4} \P{5}{\vect{0},\vect{1}}{5}} \notag \\
        & = \splitSpan{\A}{-\qi{4} \mathfig{0.11}{webs/sl_5/P_000_111_1} + \qi{3} \mathfig{0.11}{webs/sl_5/P_000_111_2} - \qi{2} \mathfig{0.11}{webs/sl_5/P_000_111_3} + \mathfig{0.11}{webs/sl_5/P_000_111_4},}{ \mathfig{0.11}{webs/sl_5/P_000_111_2} - \qi{2} \mathfig{0.11}{webs/sl_5/P_000_111_3} + \qi{3} \mathfig{0.11}{webs/sl_5/P_000_111_4} - \qi{4} \mathfig{0.11}{webs/sl_5/P_000_111_5}}
\end{align}
\begin{align}
\APR^5_{(0,0,0),(1,1,2)} & = \splitSpan{\A}{\P{5}{(0, 0, 0),(1, 1, 2)}{2} - \P{5}{(0, 0, 0),(1, 1, 2)}{3} +}{+ \P{5}{(0, 0, 0),(1, 1, 2)}{4} - \P{5}{(0, 0, 0),(1, 1, 2)}{5}} \notag \\
                         & = \Span{\A}{\mathfig{0.11}{webs/sl_5/P_000_112_2} - \mathfig{0.11}{webs/sl_5/P_000_112_3} + \mathfig{0.11}{webs/sl_5/P_000_112_4} - \mathfig{0.11}{webs/sl_5/P_000_112_5}} \displaybreak[1] \\
\APR^5_{(0,0,1),(1,2,2)} & = \splitSpan{\A}{\P{5}{(0, 0, 1),(1, 2, 2)}{2} - \P{5}{(0, 0, 1),(1, 2, 2)}{3} +}{+ \P{5}{(0, 0, 1),(1, 2, 2)}{4} - \P{5}{(0, 0, 1),(1, 2, 2)}{5}} \notag \\
                         & = \Span{\A}{\mathfig{0.11}{webs/sl_5/P_000_112_2} - \mathfig{0.11}{webs/sl_5/P_000_112_3} + \mathfig{0.11}{webs/sl_5/P_000_112_4} - \mathfig{0.11}{webs/sl_5/P_000_112_5}} \displaybreak[1] \\
\APR^5_{(0,1,1),(2,2,2)} & = \splitSpan{\A}{\P{5}{(0, 1, 1),(2, 2, 2)}{2} - \P{5}{(0, 1, 1),(2, 2, 2)}{3} +}{+ \P{5}{(0, 1, 1),(2, 2, 2)}{4} - \P{5}{(0, 1, 1),(2, 2, 2)}{5}} \notag \\
                         & = \Span{\A}{\mathfig{0.11}{webs/sl_5/P_011_222_2} - \mathfig{0.11}{webs/sl_5/P_011_222_3} + \mathfig{0.11}{webs/sl_5/P_011_222_4} - \mathfig{0.11}{webs/sl_5/P_011_222_5}}
\end{align}
\begin{align}
\APR^5_{(0,0,0,0),(1,1,1,1)} & = \splitSpan{\A}{-\P{5}{(0, 0, 0, 0),(1, 1, 1, 1)}{1} + \P{5}{(0, 0, 0, 0),(1, 1, 1, 1)}{2} - \P{5}{(0, 0, 0, 0),(1, 1, 1, 1)}{3} +}{+ \P{5}{(0, 0, 0, 0),(1, 1, 1, 1)}{4} - \P{5}{(0, 0, 0, 0),(1, 1, 1, 1)}{5}} \notag \\
                             & = \splitSpan{\A}{-\mathfig{0.12}{webs/sl_5/octagon_1} + \mathfig{0.12}{webs/sl_5/octagon_2} - \mathfig{0.12}{webs/sl_5/octagon_3} +}{+ \mathfig{0.12}{webs/sl_5/octagon_4} - \mathfig{0.12}{webs/sl_5/octagon_5}}
\end{align}

The interesting $\SS^5$ relations are
\begin{align*}
\SS^5_{(0, 0),(1, 1)} & = \splitSpanThree{\A}{\P{5}{(0, 0),(1, 1)}{2} - \Q{5}{(0, 0),(1, 1)}{0} - \qi{3} \Q{5}{(0, 0),(1, 1)}{1},}{\P{5}{(0, 0),(1, 1)}{3} - \qi{3} \Q{5}{(0, 0),(1, 1)}{0} - \qi{3} \Q{5}{(0, 0),(1, 1)}{1},}{\P{5}{(0, 0),(1, 1)}{4} - \qi{3} \Q{5}{(0, 0),(1, 1)}{0} - \Q{5}{(0, 0),(1, 1)}{1}} \\
                      & = \splitSpanThree{\A}{\mathfig{0.1}{webs/sl_4/P_00_11_2} - \mathfig{0.1}{webs/two_strands_vertical} - \qi{3} \mathfig{0.1}{webs/two_strands_horizontal},}{\mathfig{0.1}{webs/sl_4/P_00_11_3} - \qi{3} \mathfig{0.1}{webs/two_strands_vertical} - \qi{3} \mathfig{0.1}{webs/two_strands_horizontal},}{\mathfig{0.1}{webs/sl_5/P_00_11_4} - \qi{3} \mathfig{0.1}{webs/two_strands_vertical} - \mathfig{0.1}{webs/two_strands_horizontal}}, \displaybreak[1] \\
\SS^5_{(0, 0),(1, 2)} & = \splitSpan{\A}{\P{5}{(0, 0),(1, 2)}{3} - \Q{5}{(0, 0),(1, 2)}{0} - \qi{2} \Q{5}{(0, 0),(1, 2)}{1},}{\P{5}{(0, 0),(1, 2)}{4} - \qi{2} \Q{5}{(0, 0),(1, 2)}{0} - \Q{5}{(0, 0),(1, 2)}{1}} \\
                      & = \splitSpan{\A}{\mathfig{0.1}{webs/sl_4/P_00_12_3} - \mathfig{0.1}{webs/sl_4/Q_00_12_0} - \qi{2} \mathfig{0.051}{webs/sl_4/Q_00_12_1},}{\mathfig{0.1}{webs/sl_5/P_00_12_4} - \qi{2} \mathfig{0.1}{webs/sl_4/Q_00_12_0} - \mathfig{0.051}{webs/sl_4/Q_00_12_1}},  \displaybreak[1] \\
\SS^5_{(0, 0),(1, 3)} & = \Span{\A}{\P{5}{(0, 0),(1, 3)}{4} - \Q{5}{(0, 0),(1, 3)}{0} - \Q{5}{(0, 0),(1, 3)}{1}} \\
                      & = \Span{\A}{\mathfig{0.1}{webs/sl_5/P_00_13_4} - \mathfig{0.1}{webs/sl_5/Q_00_13_0} - \mathfig{0.06}{webs/sl_5/Q_00_13_1}},  \displaybreak[1] \\
\SS^5_{(0, 0),(2, 2)} & = \splitSpan{\A}{\P{5}{(0, 0),(2, 2)}{3} - \Q{5}{(0, 0),(2, 2)}{1} - \Q{5}{(0, 0),(2, 2)}{2},}{\P{5}{(0, 0),(2, 2)}{4} - \Q{5}{(0, 0),(2, 2)}{0} - \Q{5}{(0, 0),(2, 2)}{1}} \\
                      & = \splitSpan{\A}{\mathfig{0.1}{webs/sl_4/P_00_22_3} - \mathfig{0.1}{webs/sl_4/Q_00_22_1} - \mathfig{0.1}{webs/sl_5/Q_00_22_2},}{\mathfig{0.1}{webs/sl_5/P_00_22_4} - \mathfig{0.1}{webs/sl_5/Q_00_22_0} - \mathfig{0.1}{webs/sl_4/Q_00_22_1}},  \displaybreak[1] \\
\SS^5_{(0, 1),(2, 2)} & = \splitSpan{\A}{\P{5}{(0, 1),(2, 2)}{3} - \Q{5}{(0, 1),(2, 2)}{1} - \qi{2} \Q{5}{(0, 1),(2, 2)}{2},}{\P{5}{(0, 1),(2, 2)}{4} - \qi{2} \Q{5}{(0, 1),(2, 2)}{1} - \Q{5}{(0, 1),(2, 2)}{2}} \\
                      & = \splitSpan{\A}{\mathfig{0.1}{webs/sl_4/P_01_22_3} - \mathfig{0.11}{webs/sl_4/Q_01_22_1} - \qi{2} \mathfig{0.1}{webs/sl_4/Q_01_22_2},}{\mathfig{0.1}{webs/sl_5/P_01_22_4} - \qi{2} \mathfig{0.11}{webs/sl_4/Q_01_22_1} - \mathfig{0.1}{webs/sl_4/Q_01_22_2}},  \displaybreak[1] \\
\SS^5_{(0, 1),(2, 3)} & = \Span{\A}{\P{5}{(0, 1),(2, 3)}{4} - \Q{5}{(0, 1),(2, 3)}{1} - \Q{5}{(0, 1),(2, 3)}{2}} \\
                      & = \Span{\A}{\mathfig{0.1}{webs/sl_5/P_01_23_4} - \mathfig{0.11}{webs/sl_5/Q_01_23_1} - \mathfig{0.051}{webs/sl_5/Q_01_23_2}}  \displaybreak[1] \\
\SS^5_{(0, 1),(3, 2)} & = \Span{\A}{\P{5}{(0, 1),(3, 2)}{4} - \Q{5}{(0, 1),(3, 2)}{1} - \Q{5}{(0, 1),(3, 2)}{2}} \\
                      & = \Span{\A}{\mathfig{0.1}{webs/sl_5/P_01_32_4} - \mathfig{0.11}{webs/sl_5/Q_01_32_1} - \mathfig{0.051}{webs/sl_5/Q_01_32_2}}, \\
\intertext{and}
\SS^5_{(0, 2),(3, 3)} & = \Span{\A}{\P{5}{(0, 2),(3, 3)}{4} - \Q{5}{(0, 2),(3, 3)}{2} - \Q{5}{(0, 2),(3, 3)}{3}} \\
                      & = \Span{\A}{\mathfig{0.1}{webs/sl_5/P_02_33_4} - \mathfig{0.11}{webs/sl_5/Q_02_33_2} - \mathfig{0.1}{webs/sl_5/Q_02_33_3}}.
\end{align*}

\section{Proofs of Theorems \ref{thm:I=H}, \ref{thm:square-switch} and \ref{thm:kekule}}
\label{sec:proofs-of-relations}%
\subsection{The $I=H$ relations}
\begin{proof}[Proof of Theorem \ref{thm:I=H}.]
We'll just show the calculation for Equation \eqref{eq:I=H}, with
`merging' vertices; the calculation for Equation
\eqref{eq:I=H-outgoing} is identical. It's an easy induction,
beginning by calculating $\dGT$ of both sides, using the calculations from \S \ref{sec:dGT-calculations}.
\begin{align*}
\dGT&\left(\mathfig{0.1}{IH/I_abc}\right)
     = (-1)^c \mathfig{0.12}{IH/I_a-1} + (-1)^{n+b+c} q^{-a} \mathfig{0.12}{IH/I_b-1}  + \\
        & \qquad \qquad \qquad + (-1)^b q^{-a-b} \mathfig{0.12}{IH/I_c-1} + (-1)^b \mathfig{0.1}{IH/I_abc}, \\
\intertext{while} \dGT&\left(\mathfig{0.1}{IH/H_abc}\right)
     = (-1)^{n+a+c} \mathfig{0.12}{IH/H_a-1} + (-1)^{n+a+b+c} q^{-a} \mathfig{0.12}{IH/H_b-1} + \\
        & \qquad \qquad \qquad + (-1)^{a+b} q^{-a-b} \mathfig{0.12}{IH/H_c-1} + (-1)^{a+b} \mathfig{0.1}{IH/H_abc} \displaybreak[1] \\
    & = (-1)^{n(a-1)} (-1)^{n+a+c} \mathfig{0.12}{IH/I_a-1} + (-1)^{na} (-1)^{n+a+b+c} q^{-a} \mathfig{0.12}{IH/I_b-1} + \\
        & \qquad + (-1)^{na} (-1)^{a+b} q^{-a-b} \mathfig{0.12}{IH/I_c-1} + (-1)^{na} (-1)^{a+b} \mathfig{0.1}{IH/I_abc},
\end{align*}
which differs by exactly an overall factor of $(-1)^{(n+1)a}$.

For the base cases of the induction, we need consider the
possibility that $a, b$ or $c$ is $0$ or $a+b+c=n$. It's actually
more convenient to consider one of the equivalent relations in
Equation \eqref{eq:I=H-equivalent}, where the corresponding base
cases are $a,b,c$ or $d=0$. Then I claim
\begin{equation}
\label{eq:I=H-base-case}%
\mathfig{0.1}{IH/I_0bcd_l} = (-1)^{(n+1)a}
\mathfig{0.1}{IH/H_0bcd_d}
\end{equation}
in $\SymCat_n$ (and similarly in the other $3$ cases when instead
$b, c$ or $d$ is zero), and so the difference is automatically in
the kernel of $\Rep$. Each of the four variations of Equation
\eqref{eq:I=H-base-case} follows from the `degeneration' relations
in $\SymCat_n$, given in Equation \eqref{eq:degeneration}. If $a=0$,
we have
\begin{equation}
\mathfig{0.1}{IH/I_0bcd_l} = \mathfig{0.1}{IH/I_0bcd_l_2} =
\mathfig{0.1}{IH/H_0bcd_d_2} = (-1)^{(n+1)0}
\mathfig{0.1}{IH/H_0bcd_d}.
\end{equation}
The middle equality here uses `tag cancellation' from Equation
\eqref{eq:tags-cancel}. The other cases are similar, but use `tag
flipping' from Equation \eqref{eq:dual-of-a-tag}, producing the
desired signs.
\end{proof}

\subsection{The square-switch relations}
\begin{proof}[Proof of Theorem \ref{thm:square-switch}.]
This is quite straightforward, by induction.
We'll just do the $n+\ssum{a} - \ssum{b} \geq 0$ case; the other
follows immediately by rotation, by \S \ref{sec:rotations}.

The base cases are easy; since either each $b_i \geq a_i$, or $\AP^n_{a,b} = \AQ^n_{a,b} = 0$, when we're at $n=0$ the only interesting case
is $a = b = (0,0)$. Then $\SS^0_{(0,0),(0,0)}$ is simply the span of $\P{0}{(0,0),(0,0)}{0} - \Q{0}{(0,0),(0,0)}{0}$; this element is actually zero in $\SymCat_0$, by
Lemma \ref{lem:P=Q}, and so automatically in the kernel of $\Rep$.

For the inductive step, we simply show
that for $l=\max{b},\dotsc,\min{a}+n$ the component of
\begin{equation}
\label{eq:dGT-SS}%
\dGT[\P{n}{a,b}{l} - \sum_{m=\max{a}}^{\min{b}}
\qBinomial{n+\ssum{a}-\ssum{b}}{m+l-\ssum{b}} \Q{n}{a,b}{m}]
\end{equation}
in $\Hom{\SymCat^{n-1}}{\eset}{\Label(a+a',b+b')}$ lies in
$\SS^{n-1}_{a+a',b+b'}$ for each $a',b'$. The expression in Equation
\eqref{eq:dGT-SS} has one matrix for each subset of the four
boundary points. The only nonzero matrix entries correspond to those
subsets for which there is a reduction path with endpoints
coinciding with the subset. Thus, numbering the four boundary points
$1,2,3$ and $4$, we have $\dGT = \dGTe + \dGT_{\{1,2\}} +
\dGT_{\{1,4\}} + \dGT_{\{3,4\}} + \dGT_{\{3,2\}} + \dGTa$. In the
notation of \S \ref{sec:path-model},
\begin{align*}
\restrict{\dGT_{\{1,2\}}}{\AP} & = \dGT_{(1,1),(0,1)}, &
\restrict{\dGT_{\{1,2\}}}{\AQ} & = \dGT_{(0,0),(-1,0)}, \\
\restrict{\dGT_{\{1,4\}}}{\AP} & = \dGT_{(1,0),(0,0)}, &
\restrict{\dGT_{\{1,4\}}}{\AQ} & = \dGT_{(1,0),(0,0)}, \\
\restrict{\dGT_{\{3,4\}}}{\AP} & = \dGT_{(1,1),(1,0)}, &
\restrict{\dGT_{\{3,4\}}}{\AQ} & = \dGT_{(0,0),(0,-1)}, \\
\restrict{\dGT_{\{3,2\}}}{\AP} & = \dGT_{(0,1),(0,0)}, &
\restrict{\dGT_{\{3,2\}}}{\AQ} & = \dGT_{(0,1),(0,0)}.
\end{align*}
We now compute all of the components of Equation \eqref{eq:dGT-SS}.
First, use Equations \eqref{eq:dGTab-on-P} and \eqref{eq:dGTab-on-Q} to write down
\begin{align*}
\dGT_{\{1,2\}}\left(\P{n}{a,b}{l}\right) & = (-1)^{\ssum{a}} q^{-n-\ssum{a}+\ssum{b}} \P{n-1}{a+(1,1),b+(0,1)}{l}, \\
\dGT_{\{1,2\}}\left(\Q{n}{a,b}{l}\right) & = (-1)^{\ssum{a}} q^{-n-\ssum{a}+\ssum{b}} \Q{n-1}{a+(0,0),b+(-1,0)}{l} \\
                                         & = (-1)^{\ssum{a}} q^{-n-\ssum{a}+\ssum{b}} \Q{n-1}{a+(1,1),b+(0,1)}{l+1}, \displaybreak[0] \\
\dGT_{\{1,4\}}\left(\P{n}{a,b}{l}\right) & = q^{-n-a_1+b_2} \P{n-1}{a+(1,0),b+(0,0)}{l}, \\
\dGT_{\{1,4\}}\left(\Q{n}{a,b}{l}\right) & = q^{-n-a_1+b_2} \Q{n-1}{a+(1,0),b+(0,0)}{l}, \displaybreak[0] \\
\dGT_{\{3,4\}}\left(\P{n}{a,b}{l}\right) & = (-1)^{\ssum{a}} q^{-n-\ssum{a}+\ssum{b}} \P{n-1}{a+(1,1),b+(1,0)}{l}, \\
\dGT_{\{3,4\}}\left(\Q{n}{a,b}{l}\right) & = (-1)^{\ssum{a}} q^{-n-\ssum{a}+\ssum{b}} \Q{n-1}{a+(0,0),b+(0,-1)}{l} \\
                                         & = (-1)^{\ssum{a}} q^{-n-\ssum{a}+\ssum{b}} \Q{n-1}{a+(1,1),b+(1,0)}{l+1}, \\
\dGT_{\{3,2\}}\left(\P{n}{a,b}{l}\right) & = q^{-a_1+b_1} \P{n-1}{a+(0,1),b+(0,0)}{l}, \\
\intertext{and}
\dGT_{\{3,2\}}\left(\Q{n}{a,b}{l}\right) & = q^{-a_1+b_1} \Q{n-1}{a+(0,1),b+(0,0)}{l}
\end{align*}
and notice that the coefficients appearing are in each case
independent of whether we're acting on $\AP$ or $\AQ$, and moreover
that the coefficients are independent of the internal flow labels
$l$. We can use this to show that $\dGT_{\{i,j\}}\left(\P{n}{a,b}{l}
- \sum_{m=\max{a}}^{\min{b}}
\qBinomial{n+\ssum{a}-\ssum{b}}{m+l-\ssum{b}} \Q{n}{a,b}{m}\right)$
lies in $\SS^{n-1}$ for each of the four pairs $\{i,j\}$. For
example,
\begin{align*}
\dGT&_{\{1,2\}}\left(\P{n}{a,b}{l} - \sum_{m=\max{a}}^{\min{b}} \qBinomial{n+\ssum{a}-\ssum{b}}{m+l-\ssum{b}} \Q{n}{a,b}{m}\right) = \\
    & = (-1)^{\ssum{a}} q^{-n-\ssum{a}+\ssum{b}} \times \\
    & \qquad \times \left(\P{n-1}{a+(1,1),b+(0,1)}{l} - \sum_{m=\max{a}}^{\min{b}} \qBinomial{n+\ssum{a}-\ssum{b}}{m+l-\ssum{b}} \Q{n-1}{a+(1,1),b+(0,1)}{m+1}\right) \\
    & = (-1)^{\ssum{a}} q^{-n-\ssum{a}+\ssum{b}} \times \\
    & \qquad \times \left(\P{n-1}{a+(1,1),b+(0,1)}{l} - {} \right. \\
    & \qquad \qquad \left. - \sum_{m=\max{a}}^{\min{b}} \qBinomial{n-1+\ssum{(a+(1,1))}-\ssum{(b+(0,1))}}{m+1+l-\ssum{(b+(0,1))}} \Q{n-1}{a+(1,1),b+(0,1)}{m+1}\right) \\
\intertext{which becomes, upon reindexing the summation,}
    & \phantom{=} (-1)^{\ssum{a}} q^{-n-\ssum{a}+\ssum{b}} \times \\
    & \qquad \times \left(\P{n-1}{a+(1,1),b+(0,1)}{l} - {} \right. \\
    & \qquad \qquad \left. - \sum_{m=\max{a}+1}^{\min{b}+1} \qBinomial{n-1+\ssum{(a+(1,1))}-\ssum{(b+(0,1))}}{m+l-\ssum{(b+(0,1))}} \Q{n-1}{a+(1,1),b+(0,1)}{m}\right)
\end{align*}
Now, making use of the fact that if $\min{b}+1 > \min{(b+(0,1))}$, the last term vanishes since $\Q{n-1}{a+(1,1),b+(0,1)}{\min{b}+1} = 0$, we can rewrite the
summation limits and obtain
\begin{multline*}
\dGT_{\{1,2\}}\left(\P{n}{a,b}{l} - \sum_{m=\max{a}}^{\min{b}} \qBinomial{n+\ssum{a}-\ssum{b}}{m+l-\ssum{b}} \Q{n}{a,b}{m}\right) = \\
    \shoveleft{= (-1)^{\ssum{a}} q^{-n-\ssum{a}+\ssum{b}} \times \left(\P{n-1}{a+(1,1),b+(0,1)}{l} - {} \right.} \\
    \left. - \sum_{m=\max{(a+(1,1))}+1}^{\min{(b+(0,1))}} \qBinomial{n-1+\ssum{(a+(1,1))}-\ssum{(b+(0,1))}}{m+l-\ssum{(b+(0,1))}} \Q{n-1}{a+(1,1),b+(0,1)}{m}\right).
\end{multline*}
The parenthesised expression is exactly an element of the spanning set for $\SS^{n-1}_{a+(1,1),b+(0,1)}$.
It remains to show that $\dGTe$ and $\dGTa$ carry $\SS^n$
into $\SS^{n-1}$. Both are similar; here's the calculation
establishing the first.
\begin{align*}
\dGT_\eset&\left(\P{n}{a,b}{l} - \sum_{m=\max{a}}^{\min{b}} \qBinomial{n+\ssum{a}-\ssum{b}}{m+l-\ssum{b}} \Q{n}{a,b}{m}\right) = \\
    & = q^{-a_1}\left(q^l \P{n-1}{a,b}{l} + q^{l-n-\ssum{a}+\ssum{b}} \P{n-1}{a,b}{l-1} - \sum_{m=\max{a}}^{\min{b}} q^{-m} q^{\ssum{b}} \qBinomial{n+\ssum{a}-\ssum{b}}{m+l-\ssum{b}} \Q{n-1}{a,b}{m}\right), \\
\intertext{which, using the $q$-binomial identity $\qBinomial{s}{t}
= q^t \qBinomial{s-1}{t} + q^{t-s} \qBinomial{s-1}{t-1}$, we can
rewrite as}
    & = q^{-a_1}\left(q^l \P{n-1}{a,b}{l} - \sum_{m=\max{a}}^{\min{b}} q^l \qBinomial{n-1+\ssum{a}-\ssum{b}}{m+l-\ssum{b}} \Q{n-1}{a,b}{m} + \right.\\
    & \qquad \left. q^{l-n-\ssum{a}+\ssum{b}} \P{n-1}{a,b}{l-1} - \sum_{m=\max{a}}^{\min{b}} q^{l-n-\ssum{a}+\ssum{b}} \qBinomial{n-1+\ssum{a}-\ssum{b}}{m+l-1-\ssum{b}} \Q{n-1}{a,b}{m}\right),
\end{align*}
which is in $\SS^{n-1}_{a,b}$.
\end{proof}

\subsection{The \Kekule relations}
We now prove Theorem \ref{thm:kekule}, describing the relations
amongst polygonal diagrams. We begin with a slight reformulation of
the goal; it turns out that the orthogonal complement $\APR^\bot$
of the relations in the dual space $\AP^*$ is easier to describe.

\begin{lem}
\label{lem:APR-complement}%
The orthogonal complement ${\APR^n_{a,b}}^\bot$ in ${\AP^n_{a,b}}^*$
is spanned by the elements
$$e^{\mathcal{P},n}_{a,b;j^*} = \sum_{k^* = \ssum{b}}^{n+\ssum{a}}
\qBinomial{n + \ssum{a} - \ssum{b}}{k^* - \ssum{b}} (\P{n}{a,b}{k^*
+ j^*})^*$$ for $j^* = - \sumhat{b}, \dotsc, -\sumtah{a}$.
\end{lem}
\begin{proof}
First note that $\APR^n_{a,b} \subset \AP^n_{a,b}$ is an $(n + \ssum{a} -
\ssum{b})$-dimensional subspace of $\AP^n_{a,b}$; the spanning set we've described for it really is linearly independent in $\SymCat_n$.
It's a subspace of a $(\min{a} - \max{b} + n + 1)$-dimensional vector space, so we expect the orthogonal complement to be $(1 +
\sumhat{b} - \sumtah{a})$-dimensional, as claimed in the lemma. We
thus need only check that each $e^{\mathcal{P},n}_{a,b;j^*}$ annihilates each
$d^{\mathcal{P},n}_{a,b;j}$, for $j = \ssum{b}, \dotsc, \ssum{a} + n - 1$. Thus we
calculate
\begin{multline*}
e^{\mathcal{P},n}_{j^*}(d^{\mathcal{P},n}_j)
     = \sum_{k = \max(-\sumhat{b}, \ssum{b} + j^* - j)}^{\min(- \sumtah{a} + 1, j^* - j + n + \ssum{a})}
        (-1)^{j+k} \qBinomial{j+k-\max{b}}{j-\ssum{b}} \times \\ \times  \qBinomial{\min{a} + n - j - k}{\ssum{a} +n -1 -j} \qBinomial{n+\ssum{a} - \ssum{b}}{j-j^*+k-\ssum{b}}
\end{multline*}
and observe that the limits of the summation are actually
irrelevant; outside the limits at least one of the three quantum
binomial coefficients is zero.

Lemma \ref{lem:ed-zero} in the appendix \S \ref{sec:qbinomial} shows
that this $q$-binomial sum is zero.
\end{proof}
\begin{cor}
The orthogonal complement ${\AQR^n_{a,b}}^\bot$ in ${\AQ^n_{a,b}}^*$
is spanned by the elements
$$e^{\mathcal{Q},n}_{a,b;j^*} = \sum_{k^* = \ssum{a}}^{\ssum{b} - n \left(\frac{\abs{\bdy}}{2} - 1\right)} \qBinomial{\ssum{b} - \ssum{a} - n \left(\frac{\abs{\bdy}}{2} - 1\right)}{k^* - \ssum{a}} \Q{n}{a,b}{j^*+k^*}^\bot$$
for $j^* = - \sumtah{a}, \dotsc, -\sumhat{b} + n \left(\frac{\abs{\bdy}}{2} - 1\right)$.
\end{cor}
\begin{proof}
This follows immediately from the previous lemma, by rotation, as in \S \ref{sec:rotations}.
\end{proof}

The proof of Theorem \ref{thm:kekule} now proceeds inductively.
We'll assume that $$\ker \Rep \cap \AP^{n-1}_{a,b} = \APR^{n-1}_{a,b}$$ for all flows $a, b$.
We'll look at $\dGTe$ first; we know that anything in $\ker \Rep
\cap \AP^n_{a,b}$ must lie inside $\dGTe^{-1}(\APR^{n-1}_{a,b})$,
and we'll discover that in fact $\dGTe^{-1}(\APR^{n-1}_{a,b}) =
\APR^n_{a,b}$, in Lemma \ref{lem-dgt-est-inverse}. Subsequently, we
need to check, in Lemma \ref{lem-dgt-ab}, that
$\dGT_{a',b'}(\APR^n_{a,b}) \subset \APR^{n-1}_{a+a',b+b'}$ for each
allowed pair of boundary flow reduction patterns $a', b'$. Modulo
these lemmas, that proves the first part of Theorem \ref{thm:kekule} (at least the statements about $\APR$; Lemma \ref{lem:aqr} deals with the statements about $\AQR$). Later, Lemma \ref{lem:rest-of-kekule} deals with the other parts of Theorem \ref{thm:kekule}.

Explicit formulas for $\dGTe$ and $\dGT_{a',b'}$ are easy to
come by, from \S \ref{sec:path-model}:
\begin{align}
 \dGTe[\P{n}{a,b}{l}] & = q^{l-a_1} \left( \P{n-1}{a,b}{l} + q^{\ssum{b} - \ssum{a} - n} \P{n-1}{a,b}{l-1} \right) \label{eq:dGTe-P-2} \\ 
\intertext{and}
 \dGT_{a',b'}(\P{n}{a,b}{l}) & = \kappa_{n,a,b,a',b'} q^{-l(m(a',b') - 1)} \P{n-1}{a+a',b+b'}{l} \label{eq:dGT-ab} 
\end{align}
where here $\kappa_{n,a,b,a',b'}$ is some constant not depending on
$l$, which for now we don't need to keep track of, but can be reconstructed from Equation \eqref{eq:dGTab-on-P}, and $m(a',b')$ is
the number of components of the reduction path indexed by $a',b'$,
which is just $\ssum{a}' - \ssum{b}'$.

It's worth noting at this point that the $\mathcal{P}$-webs
appearing in the right hand sides of Equations \eqref{eq:dGTe-P-2} and
\eqref{eq:dGT-ab} are sometimes zero; the internal flow value might
be out of the allowable range.

In particular, the allowed values of $l$ for $\P{n}{a,b}{l} \in
\AP^n_{a,b}$ are $l = \max{b}, \dotsc, \min{a} + n$, and so the
target of $\dGTe$ is one dimension smaller than its source.
Moreover, it's easy to see from Equation \eqref{eq:dGTe-P-2} that
the kernel of $\dGTe$ is exactly one dimensional. Life is a little
more complicated for $\dGT_{a',b'}$; if $b'_i = 1$ for any $i$ where
$b_i = \max{b}$, then
\begin{equation*}
 \dGT_{a',b'}\left(\P{n}{a,b}{\max{b}}\right) = 0
\end{equation*}
and if $a'_i = 1$ for every $i$ where $a_i = \min{a}$, then
\begin{equation*}
 \dGT_{a',b'}\left(\P{n}{a,b}{\min{a} + n}\right) = 0
\end{equation*}
and otherwise $\dGT_{a',b'}$ is actually nonzero on each
$\mathcal{P}$-web. Thus the kernel of $\dGT_{a',b'}$ may be $0, 1$ or $2$
dimensional.

\begin{lem}
\label{lem-dgt-est-inverse}%
The only candidates for the kernel of the representation functor on
$\AP^n_{a,b}$ are in $\APR^n_{a,b}$, since
$$\dGTe^{-1}(\APR^{n-1}_{a,b}) = \APR^n_{a,b}.$$
\end{lem}
\begin{proof}
We've already noticed that the kernel of $\dGTe$ is $1$-dimensional,
so we actually only need to check that $\dGTe[\APR^n_{a,b}] \subset
\APR^{n-1}_{a,b}$, or, equivalently, that
$$\dGTe^*({\APR^{n-1}_{a,b}}^\bot) \subset {\APR^n_{a,b}}^\bot.$$
Easily, we obtain $$\dGTe^*(\P{n-1}{a,b}{l^*}^*) = q^{l^*-a_1}
\left( \P{n}{a,b}{l^*}^* + q^{\ssum{b} - \ssum{a} - n + 1}
\P{n}{a,b}{l^*+1}^* \right),$$ simply by taking duals in Equation
\eqref{eq:dGTe-P}. 
Then
\begin{align*} 
 \dGTe^* & (e^{\mathcal{P},n-1}_{a,b;j^*})
     = \dGTe^* \left( \sum_{k^* = \ssum{b}}^{n-1+\ssum{a}} \qBinomial{n - 1 + \ssum{a} - \ssum{b}}{k^* - \ssum{b}} (\P{n-1}{a,b}{k^*+j^*})^* \right) \\
    & = q^{-a_1} \sum_{k^* = \ssum{b}}^{n-1+\ssum{a}} q^{k^* + j^*} \qBinomial{n - 1 + \ssum{a} - \ssum{b}}{k^* - \ssum{b}} \left( \P{n}{a,b}{k^*+j^*}^* + q^{\ssum{b} - \ssum{a} - n + 1} \P{n}{a,b}{k^*+j^*+1}^* \right) \\
\intertext{and by reindexing the summation of the second terms, we can rewrite this as}
 \dGTe^* & (e^{\mathcal{P},n-1}_{a,b;j^*}) = q^{j^*-a_1+\ssum{b}} \sum_{k^* = \ssum{b}}^{n+\ssum{a}} \left( q^{k^* - \ssum{b}} \qBinomial{n - 1 + \ssum{a} - \ssum{b}}{k^* - \ssum{b}} + \right. \\
    & \qquad \qquad \qquad \qquad \left. q^{k^* - n - \ssum{a}} \qBinomial{n-1+\ssum{a}-\ssum{b}}{k^*-1-\ssum{b}} \right) \P{n}{a,b}{k^*+j^*}^* \\
    & = q^{j^*-a_1+\ssum{b}} \sum_{k^* = \ssum{b}}^{n+\ssum{a}} \qBinomial{n + \ssum{a} - \ssum{b}}{k^* - \ssum{b}} \P{n}{a,b}{k^*+j^*}^* \\
    & = q^{j^*-a_1+\ssum{b}} e^{\mathcal{P},n}_{a,b;j^*},
\end{align*}
tidily completing the proof.
\end{proof}
\begin{rem}
By rotation, we can also claim that $$\dGTa^{-1}(\AQR^{n-1}_{a,b}) =
\AQR^n_{a,b},$$
and that if you're in the kernel of $\Rep$, and in $\AQ^n_{a,b}$, you must also be in $\AQR^n_{a,b}$.
\end{rem}

\begin{lem}
\label{lem-dgt-ab}%
The candidate relations $\APR^n_{a,b}$ really are in the kernel of
the representation functor, since $$\dGT_{a',b'}(\APR^n_{a,b})
\subset \APR^{n-1}_{a+a',b+b'}$$ for each pair $a',b'$ (other than the pairs $a'=b'=\vect{0}$ and $a'=b'=\vect{1}$; the previous Lemma dealt with that matrix entry of $\dGT$).
\end{lem}
\begin{proof}
For this we need the quantum Vandermonde identity
\cite{wiki:q-Vandermonde},
$$\qBinomial{x+y}{z} = q^{yz} \sum_{i=0}^y q^{-(x+y)i}
\qBinomial{y}{i} \qBinomial{x}{z-i}.$$ Throughout, we'll write $m=
\ssum{a'} - \ssum{b'}$ for the number of components of the reduction
path. We want to show
\begin{multline}
\label{eq:dGT-on-e}%
\dGT_{a',b'}^*(e^{\mathcal{P},n-1}_{a+a',b+b';j^*}) =
\kappa_{n,a+a',b+b',a',b'} \times \\ \times \sum_{k^* = \ssum{b}}^{n+\ssum{a}+m-1}
q^{-(j^*+k^*+\ssum{b'})(m-1)}
\qBinomial{n-\ssum{a}-\ssum{b}+m-1}{k^*-\ssum{b}}
\P{n}{a,b}{k^*+j^*+\ssum{b}}^*
\end{multline}
is equal to
$$\sum_{i=0}^{m-1} X_i e^{\mathcal{P},n}_{a,b,j^* + \ssum{b} + i} = \sum_{i=0}^{m-1} X_i \sum_{k_* = \ssum{b}}^{n+\ssum{a}} \qBinomial{n+\ssum{a}-\ssum{b}}{k_* - \ssum{b}} \P{n}{a,b}{k_* + j^* + i + \ssum{b'}}^*
   \in {\APR^{n}_{a,b}}^\bot$$
for some coefficients $X_i$. Looking at the coefficient of
$\P{n}{a,b}{k^*+j^*+\ssum{b'}}^n$ in Equation \eqref{eq:dGT-on-e}
and applying the quantum Vandermonde identity with
$x=n+\ssum{a}-\ssum{b}$, $y=m-1$, and $z=k^*-\ssum{b}$, we obtain
\begin{align*}
\kappa_{n,a+a',b+b',a',b'} & q^{-(j^*+k^*+\ssum{b'})(m-1)} q^{(m-1)(k^*-\ssum{b})} \times \\
& \times \sum_{i=0}^{m-1} q^{-(n+\ssum{a}-\ssum{b} + m -1)}
\qBinomial{m-1}{i}\qBinomial{n+\ssum{a}-\ssum{b}}{k^* -\ssum{b}-i}
\end{align*}
Choosing $X_i = \kappa_{n,a+a',b+b',a',b'}
q^{(m-1)(-j^*-\ssum{b}-\ssum{b'})} q^{-(n+\ssum{a}-\ssum{b}+m-1)i}$,
(which, note, is independent of $k^*$, the index of the term we're
looking at) this is exactly the coefficient of
$\P{n}{a,b}{k^*+j^*+\ssum{b'}}^n$ in $\sum_{i=0}^{m-1} X_i
e^{\mathcal{P},n}_{a,b,j^* + \ssum{b} + i}$.
\end{proof}

\begin{lem}
\label{lem:aqr}%
The intersection of $\ker \Rep$ and $\AQ$ is exactly $\AQR$.
\end{lem}
\begin{proof}
Again, this is just by rotation, see \S \ref{sec:rotations}.
\end{proof}

The next two Lemmas, and the subsequent proof of Lemma \ref{lem:breadth-size}, are kinda hairy. Hold on tight!

\begin{lem}
For any $n$-hexagonal flows $(a,b)$ of length $k \geq 3$ such that
$\AP^n_{a,b} \neq 0$ one of the following must hold:
\begin{enumerate}
\item The flows $(a,b)$ are also $(n-1)$-hexagonal.
\item The flows $(a+\vect{1},b)$ are $(n-1)$-hexagonal.
\item There is some $(a',b') \in (\{0,1\}^k)^2$, so $(a+a',b+b')$ is
$(n-1)$-hexagonal and moreover $\dGT_{a',b'}$ maps $\AP^n_{a,b}$
faithfully into $\AP^{n-1}_{a+a',b+b'}$.
\item $n \leq 3$.
\end{enumerate}
\end{lem}
\begin{proof}
For $(a,b)$ not to be $(n-1)$-hexagonal, at least $k-5$ edge
labels must be $(n-1)$, and for $(a+\vect{1},b)$ not to be
$(n-1)$-hexagonal, at least $k-5$ edge labels must be $1$.
Already, this establishes the lemma for $k > 10$ and $n>2$
(and it's trivially true for $n\leq2$).

We thus need to deal with $k = 6, 8$ or $10$. If
$k = 10$, one of the first two alternatives hold, unless
$\Label(a,b)$ is some permutation of
$(1,1,1,1,1,n-1,n-1,n-1,n-1,n-1)$. It's easy enough to see that no
permutation is possible, given the condition that the incoming edges
and outgoing edges have the same sum. For $k = 6$ or $8$,
we can be sure that both $1$ and $n-1 \in \Label(a,b)$. Notice that
this implies $\APR^n_{a,b} = \AQR^n_{a,b} = 0$, and both
$\AP^n_{a,b}$ and $\AQ^n_{a,b}$ are at most $2$-dimensional, with
$\AP^n_{a,b}$ spanned by $\P{n}{a,b}{\max{b}}$ and
$\P{n}{a,b}{\min{a}+n}$ and $\AQ^n_{a,b}$ spanned by
$\Q{n}{a,b}{\max{a}}$ and $\Q{n}{a,b}{\min{b}}$.
Let's now assume that we never see (not even `cyclically') as contiguous
sequences of boundary edges labels either
\begin{subequations}
\label{eq:forbidden-sequences}
\begin{equation}
((n-1,-),(*\leq1,+),(*,-),(*\leq1,+),(*,-),\dotsc,(*,-),(*\leq1,+),(n-1,-))
\end{equation}
or
\begin{equation}
((n-1,+),(*\leq1,-),(*,+),(*\leq1,-),(*,+),\dotsc,(*,+),(*\leq1,-),(n-1,+)).
\end{equation}
\end{subequations}
We can then prescribe a reduction path on the $\mathcal{P}$-polygons,
which corresponds to a pair $(a',b')$ establishing the third
alternative above. Consider the reduction
path which starts at each incoming edge labeled $n-1$, and ends at
the next (heading counterclockwise) outgoing edge with a label
greater than $1$, and which additionally ends at each outgoing edge
labeled $n-1$, having started at the previous incoming edge with a
label greater than $1$. That this forms a valid reduction path (i.e.
there are no overlaps between the components described) follows
immediately from the assumption of the previous paragraph, that certain sequences do not appear. The pair
$(a',b')$ corresponding (via the correspondence discussed in \S \ref{sec:path-model})
to this reduction path then satisfies the first condition of the
third alternative above, that $(a+a',b+b')$ be $(n-1)$-hexagonal.
However, it doesn't obviously satisfy the second condition. The only
way that $\dGT_{a',b'}$ might not be faithful on $\AP^n_{a,b}$ is if
it kills $\P{n}{a,b}{\max{b}}$ or $\P{n}{a,b}{\min{a}+n}$ (remember
$\max{b} = \min{a}+n-1$ or $\min{a}+n$, since $0 < \dim \AP^n_{a,b}
\leq 2$). To kill $\P{n}{a,b}{\max{b}}$, the reduction path we've
described would have to traverse an internal edge labeled $0$.
However, it turns out there's `no room' for this. Consider some
component of the reduction path starting at an incoming edge labeled
$n-1$. It then looks like:
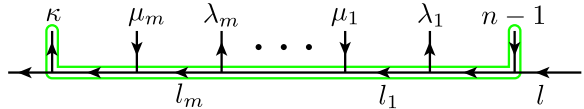
\begin{equation*}
\mathfig{0.5}{PQ/no_zero_reduction}
\end{equation*}
for some $0\leq l \leq n$, $1 < \kappa \leq n$, $\lambda_i = 0$ or $1$ and
$0\leq \mu_i \leq n$. If any internal $0$ edge gets reduced, there must be some $i$ so that the internal edge $l_i$ also gets reduced. However, $l_1 = l + n-1 -
\lambda_1 \geq l+n-2$, and generally $l_i = l + n - 1 +
\sum_{j=1}^{i-1} (-\lambda_j+\mu_i) - \lambda_i \geq
l+n+\numberof{\mu_i \geq 1}-i-1$. Thus for $l_i$ to be $0$, we must
have
\begin{equation}
\label{eq:zero-reduction-inequality}
m \geq i \geq n-1+\numberof{\mu_i \geq 1}.
\end{equation}
Now, if $k = 6$, none of the external edges can be zero,
since $(a,b)$ is $n$-hexagonal. Thus $\numberof{\mu_i \geq 1} = m$,
so as long as $n \geq 2$, the inequality of Equation
\eqref{eq:zero-reduction-inequality} cannot be satisfied. If, on the
other hand $k = 8$, we can only say $\numberof{\mu_i \geq
1} \geq m-1$ (since any non-adjacent pair of external edges being zero makes it
impossible to satisfy the `alternating' part of the condition for
being $n$-hexagon), so $m \geq n-1+m-1$, or $n \leq 2$. The argument
preventing $\dGT_{a',b'}$ from killing $\P{n}{a,b}{\min{a}+n}$ is
much the same.

Finally, we need to show that the appearance of one of the
`forbidden sequences' from Equation \eqref{eq:forbidden-sequences}
forces $n \leq 3$, so that the last alternative holds. Suppose that, for
some $m \geq 0$,
$((n-1,-),(\lambda_1,+),(\kappa_1,-),\dotsc,(\kappa_m,-),(\lambda_{m+1},+),(n-1,-)$
with $\lambda_i = 0$ or $1$ appears as a contiguous subsequence of
$\Label(a,b)$. The internal edge labels on either side of this
subsequence then differ by $2n-2+\sum_{i=1}^m \kappa_i -
\sum_{i=1}^{m+1} \lambda_i$. We know this quantity must be no more
than $n$, since $\AP^n_{a,b} \neq 0$, i.e., there is some internal
flow value $l$ so $\P{n}{a,b}{l} \neq 0$. However, $n-2+\sum_{i=1}^m
\kappa_i-\sum_{i=1}^{m+1} \lambda_i \geq n-2 + \numberof{\kappa_i
\neq 0} - \numberof{\lambda_i = 1} \geq n+\numberof{\kappa_i \neq
0}-m-3$. Thus we must have $n \leq m+3-\numberof{\kappa_i \neq 0}$.
If such a subsequence appears with $m=0$, this immediately forces
$n\leq 3$. If such a subsequence appears with $m>0$, either $n \leq
3$ or $\numberof{\kappa_i \neq 0}  < m$, so there's at least one
$\kappa_i = 0$. This certainly isn't possible for $k=3$; $(a,b)$
couldn't be $n$-hexagonal.

We're now left with a fairly finite list of cases to check for the $k=4$ case; for each
one, even in the presence of a `forbidden sequence' we'll explicitly
describe an ad hoc reduction path, so that the corresponding
$(a',b')$ pair satisfies the third alternative. We need to consider
both $m=1$ and $m=2$; there's no room for $m \geq 3$ if $k=4$. For
$m=1$, the boundary edges (up to a rotation) are
$$((n-1,-),(\lambda_1,+),(0,-),(\lambda_2,+),(n-1,-),(\mu_1,+),(\mu_2,-),(\mu_3,+))$$
with $\lambda_i = 0$ or $1$ and some $\mu_i$. In order for this to be $n$-hexagon, we
can't have $\mu_1 = n$, $\mu_2 = 0$, $\mu_3
= n$ or both $\lambda_i = 0$. Subject then to the condition that incoming and outgoing
labels have the same sum, we must have one of
\begin{align*}
\mu_2 & = 2 & \lambda_1 & = 1 & \lambda_2 & = 1 & \mu_1 & = n-1 & \mu_3 = n - 1, \\
\mu_2 & = 1 & \lambda_1 & = 1 & \lambda_2 & = 1 & \mu_1 & = n-2 & \mu_3 = n - 1, \\
\mu_2 & = 1 & \lambda_1 & = 1 & \lambda_2 & = 1 & \mu_1 & = n-1 & \mu_3 = n - 2, \\
\mu_2 & = 1 & \lambda_1 & = 0 & \lambda_2 & = 1 & \mu_1 & = n-1 & \mu_3 = n - 1, \\
\intertext{or}
\mu_2 & = 1 & \lambda_1 & = 1 & \lambda_2 & = 0 & \mu_1 & = n-1 & \mu_3 = n - 1.
\end{align*}
For the first three, however, there's an $(a',b')$ pair:
\begin{align*}
((1,0,1,1),(0,0,0,0)), \\
((1,1,1,0),(0,1,0,0)), \\
\intertext{or}
((1,0,1,1),(0,0,0,1)),
\end{align*}
and for the last two, we find that $\AP^n_{a,b} = 0$. For $m=2$, the boundary edges are
$$((n-1,-),(\lambda_1,+),(\kappa_1,-),(\lambda_2,+),(\kappa_2,-),(\lambda_3,+),(n-1,-),(\mu,+))$$
with $\lambda_i = 0$ or $1$ and at least one of the $\kappa_i = 0$. If both $\kappa_i = 0$, or if $\mu = n$, we can't achieve $n$-hexagon-ness. But then
the sum of the incoming labels is $2n-2+\kappa_1+\kappa_2 \geq 2n-1$, while the sum of the outgoing edges is $\lambda_1+\lambda_2+\lambda_3+\mu \leq n+2$.
If $2n-1 \leq n+2$, $n\leq3$ anyway, and we're done.

The argument for the forbidden sequence with opposite orientations is almost
identical.
\end{proof}

\begin{lem}
\label{lem:breadth}%
The intersection $\APR^n_{a,b} \cap \Span{\A}{\P{n}{a,b}{\max{b}}, \P{n}{a,b}{\max{b}+1}, \dotsc, \P{n}{a,b}{\ssum{b}-\sumtah{a}}}$ is zero.
\end{lem}
\begin{proof}
Say $x = \sum_{i=0}^{\sumhat{b}-\sumtah{a}} x_i \P{n}{a,b}{\max{b}+i}$, write $i_{\text{max}}$ for the greatest $i$ so $x_i \neq 0$. Define
$j^* = i_{\text{max}}- \ssum{b}$, and note that since $\max{b} \leq i_{\text{max}} \leq \ssum{b}-\sumtah{a}$ we must have $-\sumhat{b} \leq j^* \leq -\sumtah{a}$. Then apply
$e^{\mathcal{P},n}_{a,b;j^*}$, as defined in Lemma \ref{lem:APR-complement} to $x$, obtaining:
$$e^{\mathcal{P},n}_{a,b;j^*}(x) = x_{i_{\text{max}}} \neq 0,$$
so $x \notin \APR^n_{a,b}$.
\end{proof}

We now give the proof of an earlier lemma, as we're about to need it.
\begin{proof}[Proof of Lemma \ref{lem:breadth-size}.]
\newcommand{\imax}{{i_{\text{max}}}}
\newcommand{\imin}{{i_{\text{min}}}}
We can just prove the $\infty$-circumference case. Moreover, if
we're assuming $(a,b)$ has $\infty$-circumference $2k$, we can
further assume that $a$ and $b$ each actually have length $k$.
(Otherwise, there's some $i$ so $b_i = a_i$, or $b_i = a_{i+1}$; in
that case, `snip out' those entries from $a$ and $b$.) Now, pick
some $\imin$ so $a_\imin = \min{a}$, and then pick as $\imax$ the
first value of $i \geq \imin$ (in the cyclic sense) so that $b_\imax
= \max{b}$. If $\imax = \imin$, we're done; $\sumhat{b}-\sumtah{a} =
\sum_{i \neq \imax} b_i - a_i \geq k-1$. Otherwise, we'll show how
to modify $a$ and $b$, preserving the value of
$\sumhat{b}-\sumtah{a}$, to make $\imax - \imin$ (again, measured
cyclically) smaller. Specifically, for each $\imin \leq i < \imax$,
increase $b_i$ by $b_\imax - b_{\imax - 1}$, and for each $\imax < i
\leq \imax$, increase $a_i$ by the same. The resulting sequence is
still valid, in the sense that $b_i \geq a_i, a_{i+1}$ for each $i$.
Clearly, this doesn't change any of $\max{b}$, $\min{a}$, $\ssum{b}$ or
$\ssum{a}$, but it reduces the minimal value of $\imax-\imin$ by
one.
\end{proof}

We the previous three monstrosities tamed, we can finish of Theorem \ref{thm:kekule}.

\begin{lem}
\label{lem:rest-of-kekule}
If $(a,b)$ is a pair of flow labels of length at least $3$, and all
the external edge labels in $\Label(a,b)$ are between $1$ and $n-1$
inclusive, the intersection of $\ker \Rep$ and $\AP^n_{a,b} +
\AQ^n_{a,b}$ is $\APR^n_{a,b} + \AQR^n_{a,b}$.
\end{lem}
\begin{proof}
If $\AQ^n_{a,b} \neq 0$ but $\AP^n_{a,b} = 0$, first apply a rotation, as described in \S \ref{sec:rotations}.
We need to deal with the four alternatives of the previous lemma.

Suppose $n-1 \notin \Label(a,b)$, so we're in the first alternative.
Then after applying $\dGTe$, we're in a situation where
\ref{thm:kekule} holds in full; the kernel of $\Rep$ is just
$\APR^{n-1} + \AQR^{n-1}$. (The argument below works also if $1
\notin \Label(a,b)$ and we're in the second alternative, mutatis
mutandis, exchanging the roles of $\dGTe$ and $\dGTa$.) Thus suppose
$X \in \AP^n + \AQ^n$ is in $\ker \Rep$. Then $\dGTe[X] =
X_{\eset;\mathcal{P}} + X_{\eset;\mathcal{Q}}$, for some
$X_{\eset;\mathcal{P}} \in \APR^{n-1}$ and $X_{\eset;\mathcal{Q}}
\in \AQR^{n-1}$. Pick some $X_\mathcal{P} \in
\dGTe^{-1}(X_{\eset;\mathcal{P}})$, which we know by Lemma
\ref{lem-dgt-est-inverse} must also be an element of $\APR^n$. Now
$\dGTe[X-X_\mathcal{P}] = X_{\eset;\mathcal{Q}} \in \AQR^{n-1}$, so
$X-X_\mathcal{P}$ must lie in $\AQ^n$. We want to show it's actually
in $\AQR^n$. We know (via Lemma \ref{lem-dgt-ab}), that
$X-X_\mathcal{P} \in \ker \Rep$, so $\dGTa(X-X_\mathcal{P}) \in
\AQR^{n-1}$, by the remark following Lemma
\ref{lem-dgt-est-inverse}. Thus $X-X_\mathcal{P} \in \AQR^n$. The
decomposition $X = X_{\mathcal{P}} + (X-X_\mathcal{P})$ establishes
the desired result.

Suppose that there's some $(a',b')$ with
$\restrict{\dGT_{a',b'}}{\AP^n_{a,b}}$ faithful and $(a+a',b+b')$
$(n-1)$-hexagonal. Now, if $1, n-1 \notin \Label(a,b)$, as in the last
lemma we know $\AP^n_{a,b}$ is spanned by $\P{n}{a,b}{\max{b}}$ and
$\P{n}{a,b}{\min{a}+n}$ and $\AQ^n_{a,b}$  is spanned by
$\Q{n}{a,b}{\max{a}}$ and $\Q{n}{a,b}{\min{b}}$ (remember that in each case the pairs of elements might actually coincide). Thus
\begin{align*}
\dGT_{a',b'}\left(\P{n}{a,b}{\max{b}}\right)   & = (-1)^{\inner{b'}{a+\rotl[a]}} q^{\max{b}(\ssum{b'}-\ssum{a'}+1)} \times \\
                                               & \qquad \times q^{\inner{\rotl[a']}{b}-\inner{b'}{a}-a_1-n a_1'} \P{n-1}{a+a',b+b'}{\max{b}} \neq 0, \\
\dGT_{a',b'}\left(\P{n}{a,b}{\min{a}+n}\right) & = (-1)^{\inner{b'}{a+\rotl[a]}} q^{(\min{a}+n)(\ssum{b'}-\ssum{a'}+1)} \times \\
                                               & \qquad \times  q^{\inner{\rotl[a']}{b}-\inner{b'}{a}-a_1-n a_1'} \P{n-1}{a+a',b+b'}{\min{a}+n} \neq 0, \\
\dGT_{a',b'}\left(\Q{n}{a,b}{\max{a}}\right)   & = (-1)^{\inner{b'}{a+\rotl[a]}} q^{\max{a}(\ssum{a'}-\ssum{b'}-\frac{\abs{\bdy}}{2}+1)} \times \\
                                               & \qquad \times  q^{\ssum{b} - n\ssum{b'} +\inner{b'}{\rotl[a]}-\inner{a'}{b}-a_1-n a_1'} \Q{n-1}{a+a',b+b'}{\max{a}}, \\
\intertext{and}
\dGT_{a',b'}\left(\Q{n}{a,b}{\min{b}}\right)   & = (-1)^{\inner{b'}{a+\rotl[a]}} q^{\min{b}(\ssum{a'}-\ssum{b'}-\frac{\abs{\bdy}}{2}+1)} \times \\
                                               & \qquad \times  q^{\ssum{b} - n\ssum{b'} +\inner{b'}{\rotl[a]}-\inner{a'}{b}-a_1-n a_1'} \Q{n-1}{a+a',b+b'}{\min{b}}.
\end{align*}
Thus suppose $x_1 \P{n}{a,b}{\max{b}} + x_2 \P{n}{a,b}{\min{a}+n} +
y_1 \Q{n}{a,b}{\max{a}} + y_2 \Q{n}{a,b}{\min{b}}$ is in $\ker
\Rep_n$. Then $$x_1' \P{n-1}{a+a',b+b'}{\max{b}} + x_2'
\P{n-1}{a+a',b+b'}{\min{a}+n} + y_1' \Q{n-1}{a+a',b+b'}{\max{a}} +
y_2' \Q{n-1}{a+a',b+b'}{\min{b}}$$ is in $\ker \Rep_{n-1}$, where
$x_1'$ is a nonzero multiple of $x_1$, and $x_2'$ is a nonzero
multiple of $x_2$, by the formulas above. Since $(a+a',b+b')$ is
$(n-1)$-hexagonal, we can inductively assume $\ker \Rep_{n-1} = \APR^{n-1}_{a+a',b+b'} +
\AQR^{n-1}_{a+a',b+b'}$, so
$x_1' \P{n-1}{a+a',b+b'}{\max{b}} + x_2'
\P{n-1}{a+a',b+b'}{\min{a}+n} \in \APR^{n-1}_{a+a',b+b'}$. However,
by Lemma \ref{lem:breadth}, along with the knowledge that the relations in $\APR^{n-1}_{a+a',b+b'}$ have `breadth' at least $3$, as per Lemma \ref{lem:breadth-size},
this implies $x_1'=x_2'=0$, and thus
$x_1 = x_2 = 0$. Now we know $y_1 \Q{n}{a,b}{\max{a}} + y_2
\Q{n}{a,b}{\min{b}}$ is in $\ker \Rep_n$, which when restricted to
$\AQ^n_{a,b}$ is just $\AQR^n_{a,b}$, by Lemma \ref{lem:aqr}. Again
by Lemma \ref{lem:breadth}, $y_1 = y_2 = 0$.

The final alternative is $n\leq 3$, where the lemma is obvious from
the known complete descriptions of $\ker \Rep_2$ and $\ker \Rep_3$.
\end{proof}

\subsection{More about squares}
To prove Theorem \ref{thm:square-switch-complement}, that $\SS =
\APR \directSum \SS'$, we need to establish the following three
lemmas. (Trivially, $\APR \cap \SS' = 0$.) We'll in fact
only deal with the $n+\ssum{a}-\ssum{b}>0$ case of Theorem \ref{thm:square-switch-complement}; the $n+\ssum{a}-\ssum{b}=0$ is trivial, and, as usual, the $n+\ssum{a}-\ssum{b}<0$ case follows by rotation.

\begin{lem}
$\APR \subset \SS$.
\end{lem}
\begin{proof}
In the quotient $(\AP+\AQ)/\SS$,
\begin{multline*}
\sum_{k=-\sumhat{b}}^{\sumtah{a} + 1} (-1)^{j+k} \qBinomial{j+k-\max{b}}{j-\ssum{b}} \qBinomial{\min{a} + n -j -k}{\ssum{a} + n -1 -j} \P{n}{a,b}{j+k} = \\
\sum_{k=-\sumhat{b}}^{\sumtah{a} + 1} (-1)^{j+k} \qBinomial{j+k-\max{b}}{j-\ssum{b}} \qBinomial{\min{a} + n -j -k}{\ssum{a} + n -1 -j} \sum_{m=\max{a}}^{\min{b}} \qBinomial{n+\ssum{a}-\ssum{b}}{m+j+k-\ssum{b}} \Q{n}{a,b}{m}.
\end{multline*}
The coefficient of $\Q{n}{a,b}{m}$ in this is
\begin{equation*}
\sum_{k=-\sumhat{b}}^{\sumtah{a} + 1} \sum_{m=\max{a}}^{\min{b}} (-1)^{j+k} \qBinomial{j+k-\max{b}}{j-\ssum{b}} \qBinomial{\min{a} + n -j -k}{\ssum{a} + n -1 -j} \qBinomial{n+\ssum{a}-\ssum{b}}{m+j+k-\ssum{b}}.
\end{equation*}
Not only does this vanish, but each term indexed by a particular value of $m$ vanishes separately: replacing $m$ with $-j^*$, this is precisely the $q$-binomial identity of Lemma \ref{lem:ed-zero}.
\end{proof}

\begin{lem}
$\SS' \subset \SS$.
\end{lem}
\begin{proof}
We need to show that each element of the spanning set of $\SS'$
presented in Theorem \ref{thm:square-switch-complement} is in $\SS$.
In $\SS'/(\SS \cap \SS')$,
\begin{multline*}
\Q{n}{a,b}{m} - \sum_{l=n+\ssum{a}-m}^{n+\min{a}} (-1)^{m+l+n+\ssum{a}} \qBinomial{m+l-1-\ssum{b}}{m+l-n-\ssum{a}} \P{n}{a,b}{l} = \\
   = \Q{n}{a,b}{m} - \sum_{l=n+\ssum{a}-m}^{n+\min{a}} (-1)^{m+l+n+\ssum{a}} \qBinomial{m+l-1-\ssum{b}}{m+l-n-\ssum{a}} \sum_{m'=\max{a}}^{\min{b}} \qBinomial{n+\ssum{a}-\ssum{b}}{m+l-\ssum{b}} \Q{n}{a,b}{m'}.
\end{multline*}
The coefficient of $\Q{n}{a,b}{m'}$ here is
\begin{equation*}
\delta_{mm'} - (-1)^{m+n+\ssum{a}} \sum_{l=n+\ssum{a}-\min{b}}^{n+\min{a}} (-1)^l \qBinomial{m+l-1-\ssum{b}}{m+l-n-\ssum{a}} \qBinomial{n+\ssum{a}-\ssum{b}}{m'+l-\ssum{b}}.
\end{equation*}
That this is zero follows from the $q$-binomial identity proved in \S \ref{sec:qbinomial} as Lemma \ref{lem:SS'-SS}.
\end{proof}

\begin{lem}
$\SS \subset \SS' + \APR$.
\end{lem}
\begin{proof}
We take an element of the spanning set described for $\SS$, considered as an element of $(\AP + \AQ)/\SS'$. Using the relations from $\SS'$, we write
this as an element of $\AP/\SS'$, and check that this element is annihilated by each element of the spanning set for $\APR^\bot$ described in Lemma \ref{lem:APR-complement}. Thus,
\begin{multline*}
\P{n}{a,b}{l} - \sum_{m=\max{a}}^{\min{b}} \qBinomial{n+\ssum{a}-\ssum{b}}{m+l-\ssum{b}} \Q{n}{a,b}{m} = \\
 = \P{n}{a,b}{l} - \sum_{m=\max{a}}^{\min{b}} \qBinomial{n+\ssum{a}-\ssum{b}}{m+l-\ssum{b}} \sum_{l'=n+\ssum{a}-m}^{n+\min{a}} (-1)^{m+l'+n+\ssum{a}} \qBinomial{m+l'-1-\ssum{b}}{m+l'-n-\ssum{a}} \P{n}{a,b}{l'}.
\end{multline*}
Applying $\sum_{k^* = \ssum{b}}^{n+\ssum{a}} \qBinomial{n + \ssum{a} - \ssum{b}}{k^* - \ssum{b}} (\P{n}{a,b}{k^* + j^*})^*$ to this we get
\begin{multline*}
\qBinomial{n+\ssum{a}-\ssum{b}}{k^*-\ssum{b}} \delta_{l,k^*+j^*} - \\ \sum_{m=\max{a}}^{\min{b}} \sum_{l'=n+\ssum{a}-m}^{n+\min{a}} \sum_{k^* = \ssum{b}}^{n+\ssum{a}} (-1)^{m+l'+n+\ssum{a}} \qBinomial{n+\ssum{a}-\ssum{b}}{m+l-\ssum{b}} \\ \qBinomial{m+l'-1-\ssum{b}}{m+l'-n-\ssum{a}} \qBinomial{n + \ssum{a} - \ssum{b}}{k^* - \ssum{b}} \delta_{l',k^*+j^*}.
\end{multline*}
First noticing that the lower limit for $l'$ is redundant, a slight variation Lemma \ref{lem:SS'-SS} lets us evaluate the sum over $m$, obtaining
\begin{align*}
&\qBinomial{n+\ssum{a}-\ssum{b}}{k^*-\ssum{b}} \delta_{l,k^*+j^*} - \sum_{l'=-\infty}^{n+\min{a}} \sum_{k^* = \ssum{b}}^{n+\ssum{a}} \delta_{l,l'} \qBinomial{n + \ssum{a} - \ssum{b}}{k^* - \ssum{b}} \delta_{l',k^*+j^*} \\
&= \qBinomial{n+\ssum{a}-\ssum{b}}{k^*-\ssum{b}} \delta_{l,k^*+j^*} - \sum_{k^* = \ssum{b}}^{n+\ssum{a}} \qBinomial{n + \ssum{a} - \ssum{b}}{k^* - \ssum{b}} \delta_{l,k^*+j^*} \\
&= \qBinomial{n+\ssum{a}-\ssum{b}}{k^*-\ssum{b}} \delta_{l,k^*+j^*} - \qBinomial{n + \ssum{a} - \ssum{b}}{k^* - \ssum{b}} \delta_{l,k^*+j^*} \\
&= 0. \qedhere
\end{align*}
\end{proof}

\chapter{Relationships with previous work}
\label{sec:relationships}%
\section{The Temperley-Lieb category}
\label{sec:temperley-lieb}%
The usual Temperley-Lieb category (equivalently, the Kauffman bracket skein module) has unoriented strands, and a `loop value' of $-\qi{2}$. The category for $\uqsl{2}$ described here
has oriented strands, and a `loop value' of $\qi{2}$, along with orientation reversing `tags', which themselves can be flipped for a sign.
Nevertheless, the categories are equivalent, as described below (in parallel with a description of the relationship between my category for $\uqsl{4}$ and Kim's previous conjectured one) in \S \ref{sec:tags-and-orientations}.
The modification of Khovanov homology described by myself and Kevin Walker in \cite{math.GT/0701339}, producing a fully functorial invariant,
is based on a categorification of the quantum group skein module, rather than the Kauffman skein module.

\section{Kuperberg's spider for $\uqsl{3}$}
\label{sec:kuperberg}%
Kuperberg's work on $\uqsl{3}$ is stated in the language of spiders.
These are simply (strict) pivotal categories, with an alternative
set of operations emphasised; instead of composition and tensor
product, `stitch' and `join'. Further, in a spider the isomorphisms
such as $\Hom{}{a}{b\tensor c} \Iso \Hom{}{a \tensor c^*}{b}$ are
replaced with identifications.%
In this and the following sections, I'll somewhat freely mix
vocabularies.

In his work, the only edges that appear are edges labelled by $1$,
and there are no tags. As previously pointed out, the only trivalent
vertices in my construction at $n=3$ involve three edges labelled by
$1$ (here I'm not counting the flow vertices; remember they're
secretly defined in terms of the original vertices and tags). This
means that in any diagram, edges labelled by $2$ only appear either
at the boundary, connected to a tag, or in the middle of an internal
edge, with tags on either sides. Since tags cancel without signs in
the $n=3$ theory, we can ignore all the internal tags.

There's thus an easy equivalence between Kuperberg's $n=3$ pivotal
category (defined implicitly by his spider) and mine. (In fact, this
equivalence holds before or after quotienting out by the appropriate
relations.) In one direction, to his category, we send $(1,\pm) \To
(1,\pm)$ and $(2,\pm) \To (1,\mp)$ at the level of objects, and
`chop off' any external $2$-edges, and their associated tag, at the
level of morphisms.\footnote{This picturesque description needs a
patch for the identity $2$-edge; there we create a pair of tags
first, so we have something to chop off at either end.} The functor
the other direction is just the `inclusion'. Of the two composition
functors, the one on his category is actually the identity; the one
on mine is naturally isomorphic to the identity, via the tensor
natural transformation defined by
\begin{align*}
\phi_{(1,+)} & = \mathfig{0.02}{webs/strand_1} & \phi_{(1,-)} & = \mathfig{0.02}{webs/strand_1_down} \\
\phi_{(2,+)} & = \mathfig{0.032}{webs/sl_3/dt2} & \phi_{(2,-)} & = \mathfig{0.032}{webs/sl_3/da2} \\
\end{align*}

\section{Kim's proposed spider for $\uqsl{4}$}
\label{sec:kim}%
In the final chapter of his Ph.D. thesis, Kim \cite{math.QA/0310143}
conjectured relations for the $\uqsl{4}$ spider. These relations agree exactly
with mine at $n=4$. He discovered his relations by calculating the dimensions
of some small $\uqsl{4}$ invariant spaces. As soon as it's possible to write
down more diagrams with a given boundary than the dimension of the
corresponding space in the representation theory, there must be linear
combinations of these diagrams in the kernel. Consistency conditions coming
from composing with fixed diagrams which reduce the size of the boundary
enabled him to pin down all the coefficients, although, as with my work, he's
unable to show that he's found generators for the kernel.

\section{Tags and orientations}
\label{sec:tags-and-orientations}
In this section I'll describe an equivalence between my categories
and the categories previously described for $n=2$, in \S
\ref{sec:temperley-lieb} and for $n=4$, in \S \ref{sec:kim}. The
description will also encompass the equivalence described for $n=3$
in \S \ref{sec:kuperberg}; I'm doing these separately because
further complications arise when $n$ is even.

The usual spider for $\uqsl{2}$, the Temperley-Lieb category, has
unoriented edges. Similarly, Kim's proposed spider for $\uqsl{4}$
does not specify orientations on the `thick' edges, that is, those
edges labelled by $2$. Moreover, as in Kuperberg's $\uqsl{3}$
spider, only a subset of the edge labels I use appear; he only has
edges labelled $1$ and $2$.

Nevertheless, those spiders are equivalent to the spiders described
here. First, the issue of the edge label $3$ being disallowed is
treated exactly as above in \S \ref{sec:kuperberg}; we notice that
at $n=4$, there are no vertices involving edges labelled $3$, so we
can remove any internal $3$-edges by cancelling tags, and, at the
cost of keeping track of a natural isomorphism, remove external
$3$-edges too. There's a second issue, however, caused by the
unoriented edges. The category equivalences we define will have to
add and remove orientation data, and tags.

Thus we define two functors, $\iota$, which decorates spider
diagrams with `complete orientation data', and $\pi$, which forgets
orientations and tags (entirely for $n=2$, and only on the $2$ edges
for $n=4$). The forgetful functor also `chops off' any $3$ edges (in
the $n=4$ case), just like the functor in \S \ref{sec:kuperberg}.
The decorating functor $\iota$ simply fixes an up-to-isotopy
representative of a diagram, and orients each unoriented edge up the
page, placing tags at critical points of unoriented edges as
follows:
\begin{align*}
\iota\left(\mathfig{0.1}{pairings/p_unoriented}\right)  = \mathfig{0.1}{pairings/p_tag_middle}\\
\iota\left(\mathfig{0.1}{pairings/cp_unoriented}\right) = \mathfig{0.1}{pairings/cp_tag_middle}.
\end{align*}
It's well defined because opposite tags cancel. The composition $\pi
\compose \iota$ is clearly the identity. The other composition
$\iota \compose \pi$ isn't quite, but is natural isomorphic to the
identity, via the tag maps.

\section{Murakami, Ohtsuki and Yamada's trivalent
graph invariant}%
\label{sec:MOY}%
In \cite{MR2000a:57023}, Murakami, Ohtsuki and Yamada (MOY,
hereafter) define an invariant of closed knotted trivalent graphs,
which includes as a special case the HOMFLYPT
polynomial.\footnote{There's also a `cheat sheet', containing a
terse summary of their construction, at
\url{http://katlas.math.toronto.edu/drorbn/index.php?title=Image:The_MOY_Invariant_-_Feb_2006.jpg}.}
Their graphs carry oriented edge labels, and the trivalent vertices
are exactly as my `flow vertices'. (They don't have `tags'.) They
don't make any explicit connection with $\uqsl{n}$ representation
theory, although this is certainly their motivation. In fact, they
say:
\begin{quote}
We also note that our graph invariant may be obtained (not checked
yet) by direct computations of the universal $R$-matrix. But the
advantage of our definition is that it does not require any
knowledge of quantum groups nor representation theory.
\end{quote}

One of their closed graphs can be interpreted in my diagrammatic
category $\SymCat_n$; pushing it over into the representation theory
$\Rep \uqsl{n}$, it must then evaluate to a number. Presumably, this
number must be a multiple of their evaluation (modulo replacing
their $q$ with my $q^2$), with the coefficient depending on the
vertices appearing in the diagram, but not their connectivity. Thus
given a suitable closed trivalent graph $D$, the two evaluations
would be related via
\begin{equation}
\left\langle D \right\rangle_{\text{MOY}} = \lambda(D) \Rep_n\left(D\right).
\end{equation}
Knowing this evaluation coefficient $\lambda(D)$ explicitly would be
nice. You might approach it by either `localising' the MOY
formulation\footnote{If you're interested in trying, perhaps ask
Dror Bar-Natan or myself, although we don't have that much to say;
the localising step is easy enough.}, or deriving a recursive
version of the MOY evaluation function, writing the evaluation of a
graph for $n$ in terms of the evaluation of slightly modified graphs
for $n-1$. In particular, there's an evaluation function in my
category, obtained by branching all the way down to $n=0$, which is
of almost the same form as their evaluation function. It's a sum
over multiple reduction paths, such that each edge is traversed by
as many reductions paths as its label. Moreover, each reduction path
comes with an index, indicating which step of the branching process
it is `applied' to the diagram at. This index corresponds exactly to
the labels in $\mathcal{N}$ in MOY, after a linear change of
variable. Writing down the details of this should produce a formula
for $\lambda(D)$.


Modulo the translation described in the previous paragraph, their
(unproved) Proposition A.10 is presumably equivalent to the
$n+\ssum{a}-\ssum{b} \geq 0$ case of Theorem
\ref{thm:square-switch}.

\section{Jeong and Kim on $\uqsl{n}$}
In \cite{math.GT/0506403}, Jeong and Kim independently discovered a
result analogous to both cases of Theorem \ref{thm:square-switch},
using a dimension counting argument to show that there must be such
relations, and finding coefficients by gluing on other small webs.
(Of course, they published first, and have priority on that
theorem.) They never describe an explicit map from trivalent webs to
the representation theory of $\uqsl{n}$, however; they posit
relations for loops and bigons, and $I=H$ relations, which differ
from ours up to signs, and use these to show that the relations of
Theorem \ref{thm:square-switch} must hold in order to get the
dimensions of spaces of diagrams right.

They later make a conjecture which says (in my language) that even
for $2k$-gons with $k \geq 3$, each $\mathcal{P}$-$2k$-gon can be
written in terms of $\mathcal{Q}$-$2k$-gons and smaller polygons. My
results show this conjecture is false.

\section{Other work}
\label{sec:other-work}%
Sikora, in \cite{MR2171796}, defines an invariant of oriented
$n$-valent braided ribbon graphs. His graphs do not have labels on
the edges, and allow braidings of edges. The invariant is defined by
some local relations, including some normalisations, a `traditional
skein relation' for the braiding,
$$\mathfig{0.6}{sikora/skein}$$
and a relation expressing a pair of $n$-valent vertices as a linear combination of braids
$$\mathfig{0.65}{sikora/symmetriser}.$$
Easily, every closed graph evaluates to a number: the second
relation above lets you remove all vertices, resulting in a linear
combination of links, which can be evaluated via the first relation.
He further explains the connection with Murakami, Ohtsuki and
Yamada's work, giving a formula for their evaluation in terms of his
invariant.

Although his work does not expressly use the language of a category
of diagrams, it's straightforward to make the translation. He
implicitly defines a braided tensor category, with objects generated
by a single object, the unlabeled strand, and morphisms generated by
caps, cups, crossings, and the $n$-valent vertices.

This category, with no more relations than he explicitly gives,
ought to be equivalent to the full subcategory of $\Rep \uqsl{n}$
generated by tensor powers of the standard representation. Note that
this isn't the same as the category $\FundRep \uqsl{n}$ used here;
it has even fewer objects, although again its Karoubi envelope is
the entire representation category. However, I don't think that this
equivalence is obvious, at least with the currently available
results. Certainly he proves that there is a functor to this
representation category (by explicitly construct $\uqsl{n}$
equivariant tensors for the braiding, coming from the $R$-matrix,
and for the $n$-valent vertices). Even though he additionally proves
that there are no nontrivial quotients of his invariant, this does
not prove that the functor to the representation theory is faithful.
Essentially, there's no reason why `open' diagrams, such as appear
in the $\operatorname{Hom}$ spaces of the category, shouldn't have
further relations amongst them. Any way of closing up such a
relation would have to result in a linear combination of closed
diagrams which evaluated to zero simply using the initially
specified relations. Alternatively, we could think of this as a
question about`nondegeneracy': there's a pairing on diagrams, giving
by gluing together pairs of diagrams with the same, but oppositely
oriented, boundaries. It's $\Complex(q)$-valued, since every closed
diagram can be evaluated, but it may be degenerate. That is, there
might be elements of the kernel of this pairing which do not follow
from any of the local relations. See \cite[\S 3.3]{math.GT/0612754}
for a discussion of a similar issue in $\csl{3}$ Khovanov homology.

Perhaps modulo some normalisation issues, one can write down
functors between his category and mine, at least before imposing
relations. In one direction, send an edge labeled $k$ in my category
to $k$ parallel ribbons in Sikora's category, and send each vertex,
with edges labeled $a,b$ and $c$ to either the incoming or outgoing
$n$-valent vertex in Sikora's category. In the other direction, send
a ribbon to an edge labeled $1$, an $n$-valent vertex to a tree of
trivalent vertices, with $n$ leaves labeled $1$ (up to a sign it
doesn't matter, by the $I=H$ relations, which tree we use), and a
crossing to the appropriate linear combination of the two
irreducible diagrams with boundary $((1,+),(1,+),(1,-),(1,-))$.

This suggests an obvious question: are the elements of the kernel of
the representation functor I've described generated by Sikora's
relations? They need not be, if the diagrammatic pairing on Sikora's
category is degenerate.

Yokota's work in \cite{MR1427678} is along similar lines as
Sikora's, but directed towards giving an explicit construction of
quantum $SU(n)$ invariants of $3$-manifolds.

\chapter{Future directions}
Sad though it is to say, I've barely scratched the surface here.
There are a number of obvious future directions.
\begin{itemize}
\item Prove Conjecture \ref{conjecture}; that the
whole kernel of $\Rep$ is generated by the relations I've described
via tensor product and composition.
\item Try to find a `confluent' set of relations, according to the definition of
\cite{math.QA/0609832}. There's no obvious basis of diagrams modulo the relations I've given, for $n \geq 4$; this is essentially just saying that
the relations I've presented are not confluent.
\item Find an evaluation algorithm; that is, describe how to use
the relations given here to evaluate any closed web. The authors of
\cite{math.GT/0506403} claim to do so, using only the $I=H$ and
`square-switch' relations, although I have to admit not being able
to follow their proof.
\item Complete the discussion from \S \ref{sec:MOY}, explicitly giving the relationship between the evaluation of a closed diagram in my category,
and its evaluation in the MOY theory.
\item In fact, the first \Kekule relation, appearing in Equation \eqref{eq:kekule} at
$n=4$, follows from the `square switch' and `bigon' relations, as
follows. First, applying Equation \eqref{eq:sl_4-squares12} across the middle of a hexagon, we see
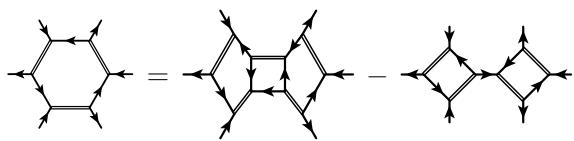
\begin{align*}
\mathfig{0.11}{webs/sl_4/kekule1_db} & = \mathfig{0.15}{webs/sl_4/calculation_1a} - \mathfig{0.15}{webs/sl_4/calculation_1b} \\
\intertext{and then, using the $I=H$ relation to switch two of the edges in the first term, and Equation \eqref{eq:sl_4-squares11} to resolve each of the squares in the second term,}
    & = \mathfig{0.15}{webs/sl_4/calculation_2a} - \left(\qi{2}^2 \mathfig{0.08}{webs/sl_4/kekule0} + \qi{2} \mathfig{0.08}{webs/sl_4/calculation_2b} + \qi{2} \mathfig{0.08}{webs/sl_4/calculation_2c} + \mathfig{0.08}{webs/sl_4/kekule3}\right). \\
\intertext{Next, we resolve the remaining squares using Equation \eqref{eq:sl_4-squares12}, obtaining}
\mathfig{0.11}{webs/sl_4/kekule1_db} & = \left(\mathfig{0.11}{webs/sl_4/kekule2_db} + \mathfig{0.1}{webs/sl_4/calculation_3a} + \mathfig{0.1}{webs/sl_4/calculation_3b} + \mathfig{0.1}{webs/sl_4/calculation_3c}\right) - \\
    & \qquad \qquad - \left(\qi{2}^2 \mathfig{0.08}{webs/sl_4/kekule0} + \qi{2} \mathfig{0.08}{webs/sl_4/calculation_2b} + \qi{2} \mathfig{0.08}{webs/sl_4/calculation_2c} + \mathfig{0.08}{webs/sl_4/kekule3}\right) \\
\intertext{and then resolve the newly created squares using Equation \eqref{eq:sl_4-squares11} again, and resolve the newly created bigons using Equation \eqref{eq:sl_4-bigon1}}
    & = \mathfig{0.11}{webs/sl_4/kekule2_db} + \left(1 + \qi{3} - \qi{2}^2 + 1\right)\mathfig{0.08}{webs/sl_4/kekule0} - \mathfig{0.08}{webs/sl_4/kekule3} \\
    & = \mathfig{0.11}{webs/sl_4/kekule2_db} + \mathfig{0.09}{webs/sl_4/kekule0} - \mathfig{0.08}{webs/sl_4/kekule3}.
\end{align*}
Is there any evidence that this continues to happen? Perhaps all the \Kekule relations follow in a similar manner from `square-switch' relations?
\item Try to do the same thing for the other simple Lie algebras.
Each Lie algebra has a set of fundamental representations,
corresponding to the nodes of its Dynkin diagram. You should first
look for intertwiners amongst these fundamental representations,
until you're sure (probably by a variation of the Schur-Weyl duality
proof given here) that they generate the entire representation
theory. After that, describe some relations, and perhaps prove that
you have all of them. I have a quite useful computer program for
discovering relations, which is general enough to attempt this
problem for any simple Lie algebra. On the other hand, the methods
of proof used here, based on the multiplicity free branching rule
for $\csl{n}$, will need considerable modification. For the notion of Gel`fand-Tsetlin basis in the orthogonal and symplectic groups, see \cite{MR1010814}. 
\item Describe the representation theory of the restricted quantum
group at a root unity in terms of diagrams; in particular find a
diagrammatic expression for the idempotent in $\tensor_a
\V{a}{n}^{m_a}$ which picks out the irrep of high weight $\sum_a m_a
\lambda_a$, for each high weight. Kim solved this problem for
$\uqsl{3}$ in his Ph.D. thesis \cite{MR2221529}. It would be nice to
start simply by describing the `lowest root of unity' quotients, in
which only the fundamental representations survive. After this, one might learn how to use these categories of trivalent graphs to
give a discussion, parallel to that in \cite{MR1710999}, of the modular tensor categories appearing at roots of unity.
\item Categorify everything in sight! Bar-Natan's work \cite{math.GT/0606318} on local
Khovanov homology provides a categorification of the $\uqsl{2}$
theory. Khovanov's work \cite{MR2100691} on a foam model for
$\csl{3}$ link homology categorifies the $\uqsl{3}$ theory, although
it's only made explicitly local in later work by Ari Nieh and myself
\cite{math.GT/0612754}. Finding an alternative to the matrix
factorisation method \cite{math.QA/0401268,math.QA/0505056} of
categorifying the $\uqsl{n}$ knot invariants, based explicitly on a
categorification of the $\uqsl{n}$ spiders, is a very tempting, and
perhaps achievable, goal! 
\end{itemize}
If you have answers to any of these problems, I'd love to hear about
them; if you have partial answers, or even just enthusiasm for the
problems, I'd love to work on them with you!

\bibliographystyle{plain}
\bibliography{bibliography/bibliography}

\appendix
\chapter{Appendices}
\section{Boring $q$-binomial identities}
\label{sec:qbinomial}%
In this section we'll prove some $q$-binomial identities needed in
various proofs in the body of the thesis. The identities are
sufficiently complicated that I've been unable to find combinatorial
arguments for them. Instead, we'll make use of the Zeilberger
algorithm, described in \cite{MR1379802}. In particular, we'll use
the {\tt{Mathematica`}} implementation in \cite{MR1395420, schorn},
and the implementation of the $q$-analogue of the Zeilberger
algorithm described in \cite{MR1448687,riese}.

If you're unhappy with the prospect of these identities being proved
with computer help, you should read Donald Knuth's foreword in
\cite{MR1379802}
\begin{quote}
Science is what we understand well enough to explain to a computer.
Art is everything else we do. During the past several years an
important part of mathematics has been transformed from an Art to a
Science: No longer do we need to get a brilliant insight in order to
evaluate sums of binomial coefficients, and many similar formulas
that arise frequently in practice; we can now follow a mechanical
procedure and discover the answers quite systematically.
\end{quote}

For each use made of the the $q$-Zeilberger algorithm, I've included
a {\tt{Mathematica}} notebook showing the calculation. These
are available in the {\tt{/code/}} subdirectory, after you've
downloaded the source of this thesis from the arXiv. 

Unfortunately, the conventions I've used for $q$-binomials don't
quite agree with those used in the implementation of the
$q$-Zeilberger algorithm. They are related via $\qBinomial{n}{k} =
q^{k(k-n)}\qBinomial[q^2]{n}{k}^\text{qZeil}$, although in my uses
of their algorithms, I've replaced $q^2$ with $q$ throughout; this
apparently produces cleaner looking results.

\begin{lem}
\label{lem:ed-zero}
For $\ssum{b} \leq j \leq \ssum{a}+n-1$
and $-\sumhat{b} \leq j^* \leq -\sumtah{a}$
\begin{multline}
\label{eq:triple-binomial}%
\sum_{k = \max(-\sumhat{b}, \ssum{b} + j^* - j)}^{\min(- \sumtah{a} + 1, j^* - j + n + \ssum{a})}
        (-1)^{j+k} \qBinomial{j+k-\max{b}}{j-\ssum{b}} \times \\ \times \qBinomial{\min{a} + n - j - k}{\ssum{a} +n -1 -j} \qBinomial{n+\ssum{a} - \ssum{b}}{j-j^*+k-\ssum{b}}
        = 0.
\end{multline}
\end{lem}
\begin{proof}
Writing $\Sigma_n$ for the sum above (suppressing the dependence on
$\ssum{a}, \ssum{b}, \min{a}, \max{b}, j$ and $j^*$), the
$q$-Zeilberger algorithm reports (see
{\tt{/code/triple-identity.nb/}}) that
\begin{equation*}
\Sigma_n = \frac{q^{(\sumtah{a}+\ssum{b}-j+j^*-1)}(1-q^{(2+2n+2\min{a}-2\max{b})})}{1-q^{(2\ssum{a}+2\sumhat{b}-2j+2j^*+2n)}} \Sigma_{n-1}
\end{equation*}
as long as
\begin{equation}
\label{eq:triple-identity-inequalities}%
n+\ssum{a}-\ssum{b} \neq 0 \qquad \qquad \text{and} \qquad \qquad
n+\ssum{a}+\sumhat{b}-j+j^* \neq 0.
\end{equation}
In the case that $n+\ssum{a}-\ssum{b} = 0$, the top of the third $q$-binomial in Equation \eqref{eq:triple-binomial} is zero, so the sum
is automatically zero unless $k = \ssum{b}+j^*-j$, in which case the sum collapses to
\begin{equation*}
(-1)^{\ssum{b}+j^*} \qBinomial{\ssum{b}+j^*}{j-\ssum{b}} \qBinomial{-\sumtah{a}-j^*}{\ssum{b}-1-j}.
\end{equation*}
However, since $j \geq \ssum{b}$, the bottom of the second $q$-binomial here is less than zero, so the whole expression vanishes.

In the case that $n+\ssum{a}+\sumhat{b}-j+j^* = 0$, we can replace $n$ everywhere in Equation \eqref{eq:triple-binomial}, obtaining
\begin{equation*}
\sum_k (-1)^{j+k} \qBinomial{j+k-\max{b}}{j-\ssum{b}} \qBinomial{-\sumtah{a} - \sumhat{b} - j^* - k}{-\sumhat{b}-j^*-1} \qBinomial{-\ssum{b}-\sumhat{b}+j-j^*}{j-j^*-k-\ssum{b}}.
\end{equation*}
Again, the second $q$-binomial vanishes, since $j^* \geq -\sumhat{b}$.

Now that we're sure the identity holds when either of the
inequalities of Equation \eqref{eq:triple-identity-inequalities} are
broken, we can easily establish the identity when they're not; if it
holds for some value of $n$, it must hold for $n+1$ (notice that the
restriction on $j$ in the hypothesis of the lemma doesn't depend on
$n$). Amusingly, the easiest starting case for the induction is not
$n=0$, but $n$ sufficiently large and negative, even though this
makes no sense in the original representation theoretic context! In
particular, at $n \leq - \ssum{a} + \ssum{b} - 1$, the third
binomial in Equation \eqref{eq:triple-binomial} is always zero, so
the sum vanishes.
\end{proof}

For the next lemma, we'll need $q$-Pochhammer symbols,
\begin{equation*}
\qPoch{a}{q}{k} =
\begin{cases}
\prod_{j=0}^{k-1} (1-a q^j) & \text{if $k>0$,} \\
1 & \text{if $k=0$,} \\
\prod_{j=1}^{\abs{k}} (1-a q^{-j})^{-1} & \text{if $k<0$.} \\
\end{cases}
\end{equation*}
In particular, note that $\qPoch{q}{q}{k<0}^{-1}=0$, and $\qPoch{q}{q}{k>0} \neq 0$.

\begin{lem}
\label{lem:SS'-SS}%
For any $n, m, m' \in \Integer$, and $a=(a_1, a_2), b=(b_1, b_2) \in \Integer^2$, with $n+\ssum{a}-\ssum{b}>0$
\begin{equation*}
\label{eq:SS-identity}%
\delta_{mm'} - (-1)^{m+n+\ssum{a}} \sum_{l=n+\ssum{a}-\min{b}}^{n+\min{a}} (-1)^l \qBinomial{m+l-1-\ssum{b}}{m+l-n-\ssum{a}} \qBinomial{n+\ssum{a}-\ssum{b}}{m'+l-\ssum{b}} = 0.
\end{equation*}
\end{lem}
\begin{proof}
First, notice that the expression is invariant under the
transformation adding some integer to each of $m, m', l, a_1, a_2, b_1, b_2$.
We'll take advantage of this to assume $\ssum{b}+2m-1,
\ssum{b}+m+m'-1$, $2\ssum{b}+2m-1-n-\ssum{a}$, $\min{a}+m+n$, $\min{a}+m'+n$, and $\ssum{b}-\max{a}+m$ are all positive,
and $-2\ssum{b}-m-m'+n-\ssum{a}$ and $\max{a}-\ssum{b}-m'-1$ are both negative.

Again, writing $\Sigma_n$ for the left hand side of the above
expression, the $q$-Zeilberger algorithm reports (see
{\tt{/code/SS-identity.nb/}}) that as long as $m \neq m'$,
\begin{align*}
\Sigma_n & = (-1)^\bullet q^\bullet
            \frac{\qPoch{q}{q}{\ssum{b}+2m-1}}
                 {(q^m-q^{m'})\qPoch{q}{q}{\ssum{b}+m+m'-1} \qPoch{q}{q}{2\ssum{b}+2m-1-n-\ssum{a}} \qPoch{q}{q}{-2\ssum{b}-m-m'+n-\ssum{a}}} + \\
         & \quad + (-1)^\bullet q^\bullet
            \frac{\qPoch{q}{q}{\min{a}+m+n}}
                 {(q^m-q^{m'})\qPoch{q}{q}{\min{a}+m'+n} \qPoch{q}{q}{\ssum{b}-\max{a}+m} \qPoch{q}{q}{\max{a}-\ssum{b}-m'-1}}.
\end{align*}
This is exactly zero, using the inequalities described in the previous paragraph, and the definition of the $q$-Pochhammer symbol.

When $m=m'$, the $q$-Zeilberger algorithm reports that
\begin{equation}
\Sigma_n = (-1)^\bullet q^\bullet \frac{(q^{n+\ssum{a}-\ssum{b}}-1)}{(q^{\min{a}+m+n}-1)\qPoch{q}{q}{\ssum{b}-\max{a}+m}\qPoch{q}{q}{-\ssum{b}+\max{a}-m}} + \Sigma_{n-1}.
\end{equation}
The first term is zero, since $-\ssum{b}+\max{a}-m<0$, so $\qPoch{q}{q}{-\ssum{b}+\max{a}-m}^{-1} = 0$.

Thus we just need to finish off the case $m=m'$, $n+\ssum{a}-\ssum{b}=1$, where the left hand side of Equation \eqref{eq:SS-identity} reduces to
\begin{equation*}
\delta_{mm'} - (-1)^{m+\ssum{b}+1} \sum_{l=\sumtah{b}+1}^{\ssum{b}-\sumtah{a}+1} (-1)^l \qBinomial{m+l-1-\ssum{b}}{m+l-1-\ssum{b}} \qBinomial{1}{m'+l-\ssum{b}}.
\end{equation*}
There are then two terms in the summation, $l=\ssum{b}-m'$ and $l=\ssum{b}-m'+1$, so we obtain
\begin{equation*}
\delta_{mm'} - (-1)^{m-m'+1} \qBinomial{m-m'-1}{m-m'-1} - (-1)^{m-m'} \qBinomial{m-m'}{m-m'}.
\end{equation*}
If $m=m'$, the first $q$-binomial vanishes, but the second is $1$, while if $m > m'$, both $q$-binomials are equal to $1$, and cancel, and if $m < m'$, both $q$-binomials vanish.
\end{proof}

\section{The $\uqsl{n}$ spider cheat sheet.}
\subsection{The diagrammatic category, \S \ref{sec:quotients}.}

\subsubsection{Generators.}
\begin{gather*}
\mathfig{0.072}{pairings/pl_a}, \qquad \mathfig{0.072}{pairings/pr_a}, \qquad \mathfig{0.072}{pairings/cpl_a}, \qquad \mathfig{0.072}{pairings/cpr_a} \\
\intertext{for each $a = 1, \dotsc, n-1$,}
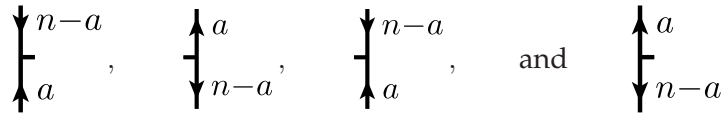
\mathfig{0.078}{duals/dtr}, \qquad
\mathfig{0.078}{duals/dal}, \qquad
\mathfig{0.078}{duals/dtl}, \qquad \text{and} \qquad
\mathfig{0.078}{duals/dar} \\
\intertext{for each $a = 1, \dotsc, n-1$, and}
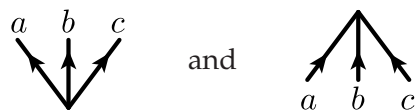
\mathfig{0.1}{vertices/v_outup_0} \qquad \text{and} \qquad
\mathfig{0.1}{vertices/v_inup_0}
\end{gather*}
for $0 \leq a,b,c \leq n$, with $a+b+c=n$.

\subsubsection{Relations: planar isotopy.}
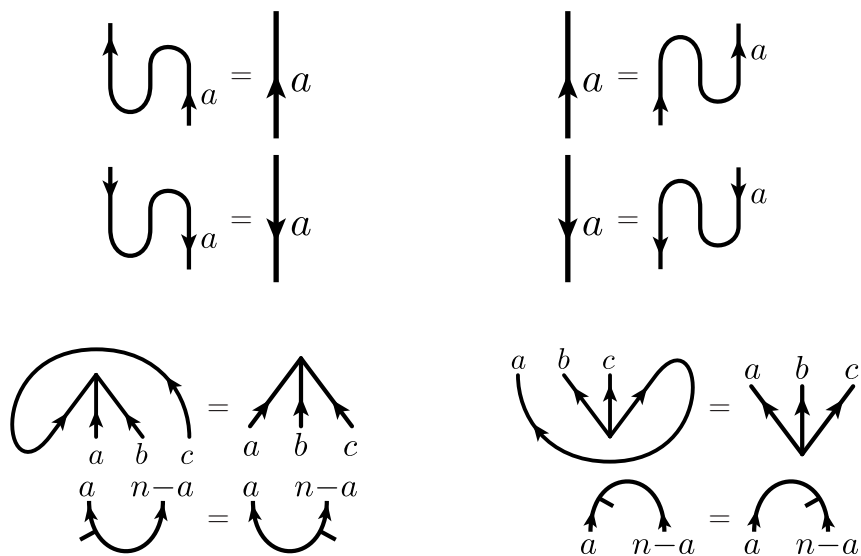
\begin{align*}
\mathfig{0.1}{isotopy/strand_left} & = \mathfig{0.04}{isotopy/strand} &
\mathfig{0.04}{isotopy/strand} & = \mathfig{0.1}{isotopy/strand_right} \displaybreak[1] \\
\mathfig{0.1}{isotopy/strand_down_left} & = \mathfig{0.04}{isotopy/strand_down} &
\mathfig{0.04}{isotopy/strand_down} & = \mathfig{0.1}{isotopy/strand_down_right}
\end{align*}
\begin{align*}
 \mathfig{0.16}{isotopy/v_in_rotated} & = \mathfig{0.1}{isotopy/v_in} & \mathfig{0.16}{isotopy/v_out_rotated} & = \mathfig{0.1}{isotopy/v_out} \displaybreak[1] \\
 \mathfig{0.1}{pairings/cp_tag_left} & = \mathfig{0.1}{pairings/cp_tag_right} &
 \mathfig{0.1}{pairings/p_tag_left} & = \mathfig{0.1}{pairings/p_tag_right}
\end{align*}
\subsubsection{Relations: (anti-)symmetric duality.}
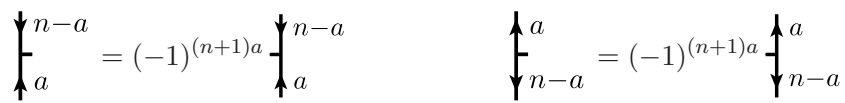
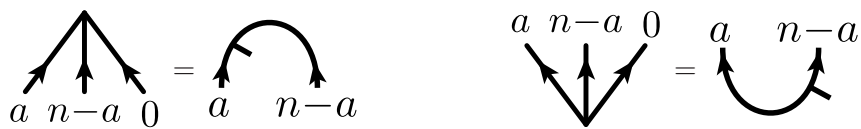
\begin{align*}
 \mathfig{0.065}{duals/dtr} & = (-1)^{(n+1)a} \mathfig{0.065}{duals/dtl} & \mathfig{0.065}{duals/dar} & = (-1)^{(n+1)a} \mathfig{0.065}{duals/dal}
\end{align*}
\subsubsection{Relations: degeneration.}
\begin{align*}
 \mathfig{0.13}{vertices/v_inup_degenerate}  & = \mathfig{0.13}{pairings/p_tag_left} &
 \mathfig{0.13}{vertices/v_outup_degenerate} & = \mathfig{0.13}{pairings/cp_tag_right}
\end{align*}

\subsection{The map from diagrams to representation theory, \S \ref{sec:representation-functor}.}
\begin{align*}
\Rep\left(\mathfig{0.1}{vertices/v_outup_0}\right) & = \mathfig{0.1}{morphisms/vout_abc} &
\Rep\left(\mathfig{0.1}{vertices/v_inup_0}\right) & = \mathfig{0.1}{morphisms/vin_abc} \\
\Rep\left(\mathfig{0.065}{duals/dtr}\right) & = \mathfig{0.1}{morphisms/d} &
\Rep\left(\mathfig{0.065}{duals/dal}\right) & = \mathfig{0.1}{morphisms/dinv} \displaybreak[1] \\
\Rep\left(\mathfig{0.065}{duals/dtl}\right) & = (-1)^{(n+1)a} \mathfig{0.1}{morphisms/d} &
\Rep\left(\mathfig{0.065}{duals/dar}\right) & = (-1)^{(n+1)a} \mathfig{0.1}{morphisms/dinv} \displaybreak[1] \\
\Rep\left(\mathfig{0.06}{pairings/pl_a}\right)  & = \mathfig{0.06}{morphisms/pl_a} &
\Rep\left(\mathfig{0.06}{pairings/cpl_a}\right) & = \mathfig{0.06}{morphisms/cpl_a}  \displaybreak[1]\\
\Rep\left(\mathfig{0.06}{pairings/pr_a}\right)  & = \mathfig{0.06}{morphisms/pr_a} &
\Rep\left(\mathfig{0.06}{pairings/cpr_a}\right) & = \mathfig{0.06}{morphisms/cpr_a} \\
\end{align*}

\subsection{The diagrammatic Gel`fand-Tsetlin functor, \S \ref{sec:dGT}.}
\subsubsection{Pairings and copairings.}
\begin{align}
 \dGT'\left(\mathfig{0.066}{pairings/pl_a}\right) & = \mathfig{0.066}{pairings/pl_a} + \mathfig{0.099}{pairings/pl_am} &
 \dGT'\left(\mathfig{0.066}{pairings/pr_a}\right) & = q^{n-a} \mathfig{0.066}{pairings/pr_a} + q^{-a} \mathfig{0.099}{pairings/pr_am} \notag \\
 \dGT'\left(\mathfig{0.066}{pairings/cpl_a}\right) & = \mathfig{0.066}{pairings/cpl_a} + \mathfig{0.099}{pairings/cpl_am} &
 \dGT'\left(\mathfig{0.066}{pairings/cpr_a}\right) & = q^{a-n} \mathfig{0.066}{pairings/cpr_a} + q^a \mathfig{0.099}{pairings/cpr_am},
 \tag{\ref{eq:dGT-on-pairings}}%
\end{align}

\subsubsection{Identifications with duals.}
\begin{align}
 \dGT'\left(\mathfig{0.065}{duals/dtr}\right) & = \mathfig{0.076}{duals/dtrb} + (-1)^a q^{-a} \mathfig{0.088}{duals/dtrt} \notag \\
 \dGT'\left(\mathfig{0.065}{duals/dtl}\right) & = (-1)^{n+a} \mathfig{0.076}{duals/dtlb} + q^{-a} \mathfig{0.088}{duals/dtlt} \notag \displaybreak[1] \\
 \dGT'\left(\mathfig{0.065}{duals/dar}\right) & = q^a \mathfig{0.088}{duals/darb} + (-1)^{n+a} \mathfig{0.076}{duals/dart} \notag \\
 \dGT'\left(\mathfig{0.065}{duals/dal}\right) & = (-1)^a q^a \mathfig{0.088}{duals/dalb} + \mathfig{0.076}{duals/dalt}. \tag{\ref{eq:dGT-on-tags}}%
\end{align}

\subsubsection{Curls.}
\begin{align*}
 \dGT\left(\mathfig{0.04}{curls/crc_a}\right) & = q^{n-a} \mathfig{0.07}{curls/crc_am} + q^{-a} \mathfig{0.04}{curls/crc_a} &
 \dGT\left(\mathfig{0.04}{curls/cra_a}\right) & = q^{a-n} \mathfig{0.07}{curls/cra_am} + q^a    \mathfig{0.04}{curls/cra_a}
\end{align*}

\subsubsection{Trivalent vertices.}
\begin{align*}
 \dGT'\left(\mathfig{0.1}{vertices/v_outup_0}\right) & =
     (-1)^c q^{b+c} \mathfig{0.1}{vertices/v_outup_1} + (-1)^a q^c \mathfig{0.09}{vertices/v_outup_2} + (-1)^b \mathfig{0.1}{vertices/v_outup_3}, \\
 \dGT'\left(\mathfig{0.1}{vertices/v_inup_0}\right) & =
     (-1)^c \mathfig{0.1}{vertices/v_inup_1} + (-1)^a q^{-a} \mathfig{0.095}{vertices/v_inup_2} + (-1)^b q^{-a-b} \mathfig{0.1}{vertices/v_inup_3} \\
\end{align*}

\begin{align*}
 \dGT\left(\mathfig{0.12}{vertices/sideways_vertex1_0}\right) & = (-1)^a \left(
    \mathfig{0.12}{vertices/sideways_vertex1_0} +
    (-1)^{n+b} \mathfig{0.155}{vertices/sideways_vertex1_1} +
    q^{-a} \mathfig{0.135}{vertices/sideways_vertex1_2}
 \right) \displaybreak[1] \\
 \dGT\left(\mathfig{0.12}{vertices/sideways_vertex2_0}\right) & = (-1)^a \left(
    \mathfig{0.12}{vertices/sideways_vertex2_0} +
    (-1)^{n+b} q^b \mathfig{0.155}{vertices/sideways_vertex2_1} +
    \mathfig{0.145}{vertices/sideways_vertex2_2}
 \right)
\end{align*}

\begin{align*}
 \dGT\left(\mathfig{0.12}{vertices/sideways_vertex3_0}\right) & = (-1)^a \left(
    \mathfig{0.12}{vertices/sideways_vertex3_0} +
    (-1)^{n+b} q^b \mathfig{0.12}{vertices/sideways_vertex3_1} +
    \mathfig{0.12}{vertices/sideways_vertex3_2}
 \right) \displaybreak[1] \\
 \dGT\left(\mathfig{0.12}{vertices/sideways_vertex4_0}\right) & = (-1)^a \left(
    q^a \mathfig{0.12}{vertices/sideways_vertex4_0} +
    \mathfig{0.12}{vertices/sideways_vertex4_1} +
    (-1)^{n+b} q^{a-n} \mathfig{0.12}{vertices/sideways_vertex4_2}
 \right)
\end{align*}

\subsubsection{Larger webs.}
\begin{multline}
\dGT\left(\mathfig{0.25}{PQ/P_section_0}\right) =
    \mathfig{0.25}{PQ/P_section_0} + {} \\
  + q^{-a} \mathfig{0.26}{PQ/P_section_1} + \mathfig{0.26}{PQ/P_section_2} + {} \\
  + q^{-a} \mathfig{0.275}{PQ/P_section_3} + (-1)^{b+c} q^b \mathfig{0.275}{PQ/P_section_4} \tag{\ref{eq:dGT-P-section}}
\end{multline}
\begin{multline}
\dGT\left(\mathfig{0.18}{PQ/Q_section_0}\right) =
    q^b \mathfig{0.18}{PQ/Q_section_0} + {} \\
  + \mathfig{0.18}{PQ/Q_section_1} + \mathfig{0.18}{PQ/Q_section_2} + {} \\
  + q^b \mathfig{0.18}{PQ/Q_section_3} + (-1)^{a+c} q^{b+c-n} \mathfig{0.18}{PQ/Q_section_4}  \tag{\ref{eq:dGT-Q-section}}
\end{multline}

\subsubsection{$\mathcal{P}$'s and $\mathcal{Q}$'s, \S \ref{sec:polygons}.}
\begin{align*}
 \dGT_{a',b'}(\P{n}{a,b}{l}) & = (-1)^{b' \cdot (a+\rotl{a})} q^{l(\ssum{b'}-\ssum{a'}+1)} q^{\rotl{a'}\cdot b - b' \cdot a - a_1 - n a_1'} \P{n-1}{a+a',b+b'}{l} \\
 \dGT_{a',b'}(\Q{n}{a,b}{l}) & = (-1)^{b' \cdot (a+\rotl{a})} q^{l(\ssum{a'}-\ssum{b'} - \frac{\abs{\bdy}}{2} + 1)} q^{\ssum{b}+n\ssum{b'}+b' \cdot \rotl{a} - a' \cdot b - a_1 - n a'_1} \Q{n-1}{a+a',b+b'}{l} \\
\intertext{and}
 \dGTe(\P{n}{a,b}{l}) & = \dGT_{\vect{0},\vect{0}}(\P{n}{a,b}{l}) + \dGT_{\vect{1},\vect{1}}(\P{n}{a,b}{l}) \\
                           & = q^{l-a_1} \left( \P{n-1}{a,b}{l} + q^{-n-\ssum{a}+\ssum{b}} \P{n-1}{a,b}{l-1} \right) \\
 \dGTe(\Q{n}{a,b}{l}) & = \dGT_{\vect{0},\vect{0}}(\Q{n}{a,b}{l}) \\
                           & = q^{l-a_1} q^{\ssum{b} - \frac{l \abs{\bdy}}{2}} \Q{n-1}{a,b}{l} \displaybreak[1] \\
 \dGTa(\P{n}{a,b}{l}) & = \dGT_{\vect{1}, \vect{0}}(\P{n}{a,b}{l}) \\
                           & = q^{l-a_1} q^{\ssum{b} - \frac{l \abs{\bdy}}{2} -n } \P{n-1}{a+\vect{1},b}{l} \\
 \dGTa(\Q{n}{a,b}{l}) & = \dGT_{\vect{0},\vect{-1}}(\Q{n}{a,b}{l}) + \dGT_{\vect{1},\vect{0}}(\Q{n}{a,b}{l}) \\
                           & = q^{l-a_1} \left( q^{\ssum{b} - \ssum{a}-\frac{n\abs{\bdy}}{2}} \Q{n-1}{a+\vect{1},b}{l+1} + q^{-n} \Q{n-1}{a+\vect{1},b}{l} \right)
\end{align*}

\subsection{Conjectural complete set of relations, \S \ref{sec:kernel}.}
\subsubsection{The $I=H$ relations}
\begin{align*}
\mathfig{0.1}{IH/I_abc} & = (-1)^{(n+1)a} \mathfig{0.1}{IH/H_abc}, \\
\mathfig{0.1}{IH/I_out_abc} & = (-1)^{(n+1)a} \mathfig{0.1}{IH/H_out_abc}.
\end{align*}

\subsubsection{The `square-switching' relations}
\begin{equation*}
 \SS^n_{a,b} =
    \begin{cases}
        \Span{\A}{\P{n}{a,b}{l} - \sum_{m=\max{a}}^{\min{b}} \qBinomial{n+\ssum{a}-\ssum{b}}{m+l-\ssum{b}} \Q{n}{a,b}{m}}_{l=\max{b}}^{\min{a}+n} & \text{if $n+\ssum{a} - \ssum{b} \geq 0$} \\
        \Span{\A}{\Q{n}{a,b}{l} - \sum_{m=\max{b}}^{n+\min{a}} \qBinomial{\ssum{b}-n-\ssum{a}}{m+l-\ssum{a}-n} \P{n}{a,b}{m}}_{l=\max{a}}^{\min{b}} & \text{if $n+\ssum{a} - \ssum{b} \leq 0$}
    \end{cases}
\end{equation*}

\subsubsection{The \Kekule relations}
\begin{multline}
 \APR^n_{a,b} = \operatorname{span}_{\A}\left\{ d^{\mathcal{P},n}_{a,b;j} = \sum_{k=-\sumhat{b}}^{-\sumtah{a} + 1} (-1)^{j+k} \qBinomial{j+k-\max{b}}{j-\ssum{b}} \times \right. \\ \left. \rule{0mm}{10mm} \times \qBinomial{\min{a} + n -j -k}{\ssum{a} + n -1 -j} \P{n}{a,b}{j+k}\right\}_{j=\ssum{b}}^{\ssum{a} +n-1}
\end{multline}
\begin{multline}
 \AQR^n_{a,b} = \operatorname{span}_{\A}\left\{  d^{\mathcal{P},n}_{a,b;j} = \sum_{k=-\sumtah{a}}^{-\sumhat{b} + \frac{n\abs{\bdy}}{2} + 1} (-1)^{j+k} \qBinomial{j+k-\max{a}}{j-\ssum{a}}  \times \right. \\ \left. \rule{0mm}{10mm} \times \qBinomial{\min{b}-j-k}{\ssum{b}-n(\frac{\abs{\bdy}}{2}-1)-1-j} \Q{n}{a,b}{j+k} \right\}_{j=\ssum{a}}^{\ssum{b}-n(\frac{\abs{\bdy}}{2}-1)-1}
\end{multline}
with orthogonal complements
\begin{multline}
 {\APR^n_{a,b}}^\bot = \operatorname{span}_{\A}\left\{  e^{\mathcal{P},n}_{a,b;j^*} = \sum_{k^* = \ssum{b}}^{n+\ssum{a}} \qBinomial{n+\ssum{a}-\ssum{b}}{k^*-\ssum{b}}  \P{n}{a,b}{j^*+k^*}^\bot \right\}_{j^* = - \sumhat{b}}^{-\sumtah{a}} \\
 \shoveleft{
     {\AQR^n_{a,b}}^\bot = \operatorname{span}_{\A}\left\{  e^{\mathcal{Q},n}_{a,b;j^*} = \sum_{k^* = \ssum{a}}^{\ssum{b} - n \left(\frac{\abs{\bdy}}{2} - 1\right)} \qBinomial{\ssum{b} - \ssum{a} - n \left(\frac{\abs{\bdy}}{2} - 1\right)}{k^* - \ssum{a}}  \times \right.} \\
 \left. \times \rule{0mm}{10mm} \Q{n}{a,b}{j^*+k^*}^\bot \right\}_{j^* = - \sumtah{a}}^{-\sumhat{b} + n \left(\frac{\abs{\bdy}}{2} - 1\right)}
\end{multline}

\end{document}